\documentclass[a4paper,10pt]{amsart}
\usepackage{amsmath,amsfonts,amssymb, amsthm, xypic}
\usepackage{epsfig}
\usepackage[all]{xy}

\renewcommand{\deg}{\mathsf{deg}}
\DeclareMathOperator{\pr}{\mathsf{pr}}
\DeclareMathOperator{\rank}{\mathsf{rank}}
\DeclareMathOperator{\Aut}{\mathsf{Aut}}
\DeclareMathOperator{\Spec}{\mathrm{Spec}}
\DeclareMathOperator{\jet}{\mathsf{jet}}
\DeclareMathOperator{\ev}{\mathsf{ev}}
\DeclareMathOperator{\id}{\mathsf{id}}
\renewcommand{\span}{\mathsf{span}}
\DeclareMathOperator{\cC}{\mathsf{C}}

\def\Hom{\mathrm{Hom}}
\def\Z{\mathbb{Z}}
\def\Q{\mathbb{Q}}
\def\E{X}
\def\S{\mathcal{S}}
\def\Y{\mathcal{Y}}
\def\Ext{\mathrm{Ext}}
\def\H{\mathbf{H}}
\def\U{\mathbf{U}}
\def\UU{\boldsymbol{\mathcal{E}}}
\def\qed{$\hfill \checkmark$}
\def\N{\mathbb{N}}
\def\tto{\twoheadrightarrow}
\def\lto{\longrightarrow}
\def\a{\alpha}
\def\b{\beta}
\def\q{\mathbf{q}}
\def\qlb{\overline{\mathbb{Q}_l}}
\def\C{\mathbb{C}}

\def\x{\mathbf{x}}
\def\y{\mathbf{y}}
\def\z{\mathbf{z}}
\def\p{\mathbf{p}}

\def\kk{\boldsymbol{k}}
\def\KK{K}
\def\Rb{\mathbf{R}}
\def\Kb{\mathbf{K}}
\def\ll{\boldsymbol{l}}
\def\LL{\boldsymbol{\Lambda} \hspace{-.06in} \boldsymbol{\Lambda}}
\def\LLambda{\boldsymbol{\Lambda}}
\def\ZZ{\mathbf{Z}}
\def\s{\sigma}
\def\bs{\bar{\sigma}}
\def\vv{\nu}

\newtheorem{theo}{\bf{Theorem}}[section]
\newtheorem{lem}[theo]{Lemma}
\newtheorem{cor}[theo]{Corollary}
\newtheorem{prop}[theo]{Proposition}

\numberwithin{equation}{section}

\newcommand{\kA}{\mathcal{A}}
\newcommand{\kB}{\mathcal{B}}
\newcommand{\kC}{\mathcal{C}}
\newcommand{\kD}{\mathcal{D}}

\newcommand{\kE}{\mathcal{E}}
\newcommand{\kF}{\mathcal{F}}
\newcommand{\kG}{\mathcal{G}}
\newcommand{\kH}{\mathcal{H}}
\newcommand{\kO}{\mathcal{O}}
\newcommand{\kL}{\mathcal{L}}
\newcommand{\kP}{\mathcal{P}}
\newcommand{\kK}{\mathcal{K}}
\newcommand{\kM}{\mathcal{M}}
\newcommand{\kN}{\mathcal{N}}

\title{On the Hall algebra of an elliptic curve, I}
\author{Igor Burban* and Olivier Schiffmann$^\dag$}

\thanks{
\noindent
\\
*
Mathematisches Institut,
Friedrich-Wilhelms-Universit\"at Bonn,
Endenicher Allee 60,
D-53115, Bonn, Germany,
e-mail:\; \texttt{burban@math.uni-bonn.de} \\
$\dag$ Institut Math\'ematique de Jussieu, Paris 6, 175 rue du Chevaleret, 75013 Paris,  France,
\\ e-mail:\;\texttt{olive@math.jussieu.fr}}

\begin{document}
\maketitle

{{\flushright{{
\quote{\flushright{\textit{Forests may fall,\\
But not the dusk they shield.}\\
H.P. Lovecraft \\
}}}}}}

\tableofcontents
\setlength{\unitlength}{10pt}

\centerline{I\footnotesize{NTRODUCTION}}

\vspace{.2in}

\paragraph{}
 Among the oldest and still most fundamental objects in representation theory and
combinatorics are the
rings of symmetric polynomials
$$\LLambda^+= \C[x_1, x_2, \ldots]^{\mathfrak{S}_{\infty}}:=\underset{\longleftarrow}{\text{Lim}}\;
\C[x_1, \ldots, x_r]^{\mathfrak{S}_r},$$
and symmetric Laurent polynomials
$$\LLambda= \C[x_1^{\pm 1}, x_2^{\pm 1}, \ldots]^{\mathfrak{S}_{\infty}}.$$
These rings admit  numerous algebraic and geometric realizations, but one of
the historically first constructions, dating  to the work of
Steinitz in 1900   completed later by Hall, was given in terms of what is now called
the classical Hall algebra $\mathbf{H}$ (see \cite{Mac}, Chapter II
). This algebra has a basis consisting of isomorphism classes of
abelian $q$-groups, where $q$  is a fixed prime power,  and the structure
constants are defined by counting extensions between such abelian
groups. In fact, these structure constants are polynomials in $q$, and we
can therefore  consider $\mathbf{H}$ as a $\C[q^{\pm 1}]$-algebra.
A theorem of Steinitz and Hall provides an isomorphism $\mathbf{H} \simeq
\LLambda^+_q=\C[q^{\pm 1}][x_1, x_2,\ldots]^{\mathfrak{S}_{\infty}}$. Under
this isomorphism, the natural basis of $\mathbf{H}$ (resp. the natural
scalar product) is mapped to the basis of Hall-Littlewood polynomials
(resp. the Hall-Littlewood scalar product). In addition, Zelevinsky
\cite{Zel} endowed $\LLambda_q^+$ with a structure of a cocommutative
Hopf algebra  and the whole algebra $\LLambda_q=\LLambda \otimes \C[q^{\pm
1}]$ can  be recovered from $\LLambda^+_q$ by the Drinfeld double
construction. This Hopf algebra structure is also
intrinsically defined by means of  the
Hall algebra.

\vspace{.1in}

One aim of the present work is to initiate a similar approach for the rings of \textit{diagonal} symmetric polynomials
$$\LL^{++}=\C[x_1, x_2,\ldots, y_1, y_2, \ldots]^{\mathfrak{S}_{\infty}},
\qquad \LL^+=\C[x_1^{\pm 1}, x_2^{\pm 1},\ldots, y_1, y_2, \ldots]^{\mathfrak{S}_{\infty}}$$
and
$$\LL=\C[x^{\pm 1}_1, x^{\pm 1}_2, \ldots, y_1^{\pm 1},
y_2^{\pm 1}, \ldots]^{\mathfrak{S}_{\infty}},$$
with  $\mathfrak{S}_{\infty}$ acting simultaneously on the variables $x_i$
and $y_i$, based on the category of coherent sheaves on an elliptic curve.
These rings have recently attracted a lot of attention due to
its close relations to Macdonald's  polynomials and double affine Hecke algebras.

\vspace{.1in}

To any abelian category
$\mathcal{A}$ defined over a finite field $\kk = \mathbb{F}_q$ and
satisfying certain finiteness conditions one  can  attach an associative algebra
$\H_\mathcal{A}$ defined over the field $\mathbb{Q}(v)$, $v=\sqrt{q}^{-1}$
called the \emph{Hall algebra} of the category $\mathcal{A}$.
As a $\mathbb{Q}(v)$--vector space
$\H_\mathcal{A}$
has a basis parameterized by isomorphism classes of objects of $\mathcal{A}$ and
its structure constants are expressed via the number of extensions between the
objects of $\mathcal{A}$.
The interest in this construction grew considerably after
Ringel studied in \cite{Ri} the Hall algebra of the category of representations
of an arbitrary quiver $\vec{Q}$ and showed that it contains the positive part
$\U^+_v(\mathfrak{g})$ of the quantized enveloping algebra of the Kac-Moody
algebra $\mathfrak{g}$ associated to $\vec{Q}$.

In a similar direction, Kapranov considered
in \cite{Kap} a natural subalgebra $\mathbf{H}^{sph}_{\E}$ of the Hall algebra
$\mathbf{H}_{\E}$ of the category of coherent sheaves
$Coh({\E})$ on  a smooth projective curve ${\E}$ defined over a finite
field $\kk$. This \textit{spherical Hall algebra}
$\mathbf{H}^{sph}_{\E}$ plays an important  role in the Langlands program for the function field of $\E$
because  it can  be interpreted as the algebra of (everywhere unramified,
principal) Eisenstein series for $GL(n)$ for all $n$, with the product coming from the parabolic
induction functor.
In the case ${\E}=\mathbb{P}^1$ the algebra $\mathbf{H}^{sph}_{\E}$ is isomorphic
to the positive part of the quantum loop algebra $U_v(\mathfrak{L sl}_2)$ (see \cite{Kap} and also \cite{BK}).
In higher genus, Kapranov defined a surjective map from another  algebra ${\mathbb{U}}^+_{\E}$
(defined by generators and relations) to $\mathbf{H}^{sph}_{\E}$. Unfortunately, this map has a nontrivial kernel,
and it is not known how to describe it explicitly.

\vspace{.1in}

In this paper, we study in details  the  Hall algebra $\mathbf{H}_\E$
of an elliptic curve  $\E$  defined over $\kk$
and a certain subalgebra $\mathbf{U}^+_{\E}$ of ${\mathbf{H}}_\E$ which turns out to coincide with 
the spherical Hall algebra $\H^{sph}_{\E}$
of Kapranov. We show that $\mathbf{U}^+_{{\E}}$ is naturally a deformation of the ring of
\textit{diagonal} symmetric polynomials
$$\LL^+:=\C[x^{\pm 1}_1,x^{\pm 1}_2, \ldots, y_1,y_2,\ldots]^{\mathfrak{S}_{\infty}}.$$

In Theorem~\ref{T:main2} we provide an explicit description of the bialgebra $\mathbf{U}^+_{\E}$
by generators and relations. It is neither commutative, nor cocommutative.
In order to obtain a more symmetric and canonical
object, we consider the Drinfeld double $\mathbf{U}_\E$ of $\mathbf{U}^+_{\E}$,
which is now a deformation of the ring
$\LL=\C[x^{\pm 1}_1, x^{\pm 1}_2,
\ldots, y_1^{\pm 1},  y_2^{\pm 1}, \ldots]^{\mathfrak{S}_{\infty}}.$
We prove (Theorem~\ref{P:braidact}) that the group of exact auto-equivalences of
the derived category $D^b\bigl(Coh({\E})\bigr)$ naturally acts  on $\mathbf{U}_\E$ by \textit{algebra}
automorphisms, yielding an action of ${SL}(2,\Z)$  on $\mathbf{U}_\E$. In
Section 5 we construct a natural ``monomial'' basis
of $\mathbf{U}^+_{\E}$ (resp. of $\mathbf{U}_\E$) indexed by the
set of finite  \textit{convex}
paths in the region $(\Z^2)^+=\bigl\{(p,q) \in \Z^2\;|\; p \geq 1 \; \text{or}\;
p=0, q \geq 0\bigr\}$
(resp. in $\Z^2$). This basis is equivariant with
respect to the ${SL}(2,\Z)$-action.

\vspace{.1in}

We show that the structure constants of $\mathbf{U}_{\E}$  are Laurent polynomials in
$\sigma^{1/2}$ and $\bar{\sigma}^{1/2}$,  where $\sigma, \bar{\sigma}$ are  the Frobenius
eigenvalues on the $l$-adic cohomology group
$H^1({\E}_{\overline{\kk}}, \overline{\mathbb{Q}_l})$ (observe that $v=(\sigma \bar{\sigma})^{-1/2}$).
This allows us to consider $\U_{\E}$ as a $\C\bigl[\sigma^{\pm 1/2}, \, 
\bar{\sigma}^{\pm 1/2}\bigr]$-algebra. More precisely,
we introduce a generic version $\UU_{\Rb}$ of the Hall algebras $\U^+_\E$, which is defined over the ring
$\Rb = \C\bigl[\sigma^{\pm 1/2},\bar{\sigma}^{\pm 1/2}\bigr]$,  where $\sigma, \bar{\sigma}$ are now
formal parameters and which specializes
to all the algebras $\U_{\E}$.
Moreover, for the values
$\sigma=\bar{\sigma} = 1$ one gets the ring
$$(\UU_{\Rb})_{|\sigma=\bar{\sigma}=1} \simeq \LL=\C\bigl[x^{\pm 1}_1, x^{\pm 1}_2,
\ldots, y_1^{\pm 1},  y_2^{\pm 1}, \ldots\bigr]^{\mathfrak{S}_{\infty}}$$
of diagonal symmetric polynomials and $\UU_{\Rb}$ is a flat deformation of $\LL$.
We show that as in the case of $\U_{\E}$,
the algebra $\UU_{\Rb}$ has a monomial basis, a triangular decomposition, and carries an action of ${SL}(2,\Z)$ by automorphisms.

\vspace{.1in}

A  very interesting two-parameter deformation of the ring
$$\LL_n=\C\bigl[x_1^{\pm 1}, \ldots,
x_n^{\pm 1}, y_1^{\pm 1}, \ldots , y_n^{\pm 1}\bigr]^{\mathfrak{S}_n}$$
is provided by the \textit{spherical} double affine Hecke algebra
(DAHA) $\mathbf{S}\ddot{\mathbf{H}}_n$ of type $\mathfrak{gl}(n)$ (see \cite{Cherednik}).
In a joint  work \cite{SV} of the second-named author with E.~Vasserot
it is shown that there are surjective homomorphisms $\UU_{\Rb}
\tto \mathbf{S} \ddot{\mathbf{H}}_n$ for any
positive integer $n$, so that $\UU_{\Rb}$ may be thought of as the ``stable limit'' $\mathbf{S}\ddot{\H}_{\infty}$ of the type $A$ spherical DAHA.
In the companion paper \cite{S2}, we shall use a geometric version of the Hall algebra to construct
certain ``canonical bases'' of $\UU_{\Rb}$, which may be thought of as some ''double'' analogues of Kazhdan-Lusztig polynomials of type $A$.

\vspace{.1in}

The elliptic Hall algebra $\UU_{\Rb}$ has recently found applications in the geometric construction of Macdonald polynomials via Eisenstein series (see \cite{SV}), and in the computation of convolution algebras in the equivariant K-theory of Hilbert schemes of $\mathbb{A}^2$ and of the commuting variety (see \cite{SV2}).

\vspace{.1in}

Let us now briefly describe the content of this paper. After recalling
Atiyah's classification of coherent sheaves on an elliptic curve ${\E}$ and
the structure of the group of exact
auto-equivalences of the derived category $D^b\bigl(Coh({\E})\bigr)$ in Section~1, we introduce, following Ringel and Green, the Hall bialgebra
${\mathbf{H}}_{\E}$ of the category $Coh({\E})$ in Section~2. In
Section~3 we deal  with the Drinfeld double $\mathbf{D}{\mathbf{H}}_{\E}$ of
${\mathbf{H}}_{\E}$
and constructs an embedding of the group of exact auto-equivalences
of $D^b\bigl(Coh(\E)\bigr)$ into $\Aut(\mathbf{D}{\mathbf{H}}_{\E})$. The subalgebra
$\mathbf{U}_\E$ of $\mathbf{D}{\mathbf{H}}_{\E}$ we are interested in is
defined in Section~4.
The main theorem of this article, describing ${\mathbf{U}}_\E$ by
generators and relations is proven in Section~5.  Section 6 contains various important
properties of $\U_{\E}$ (integral form, central extension, etc). In the last Section 7
sum up main properties of the algebra $\U_{\E}$ proven in this article. Appendix A is devoted
a discussion of Fourier-Mukai transforms for elliptic curves defined over finite fields,
whereas in Appendix B we prove some basic properties of the Drinfeld double of a 
topological bialgebra.

\vspace{.1in}

\section{Coherent sheaves on elliptic curves}

\vspace{.1in}

\paragraph{\textbf{1.1.}}
Let $\kk$ be any field. Throughout the paper $\E$ denotes a smooth elliptic
curve defined over $\kk$, that is,
$\E$ is a smooth projective curve  of genus one having
a rational point. Note, that by Weil's inequality  in the case
of a finite field $\kk = \mathbb{F}_q$
we have  $\bigl||\E(\kk)| - (q+1)\bigr| \le 2\sqrt{q}$, hence  any
smooth projective curve of genus one has such a point.
We denote by $Coh(\E)$ its category of coherent sheaves.
Let us first outline, following Atiyah, the classification of coherent
sheaves on elliptic curves
(in \cite{A} it is assumed that $\kk$ is algebraically closed, but the proof can be applied for
an arbitrary field $\kk$). Recall that the slope of a sheaf
$\mathcal{F} \in Coh(\E)$ is
$\mu(\mathcal{F})=\deg(\mathcal{F})/\rank(\mathcal{F})$,
and that a sheaf $\mathcal{F}$ is \textit{semi-stable} (resp. \textit{stable}) if for any subsheaf $\mathcal{G} \subset \mathcal{F}$
we have $\mu(\mathcal{G}) \leq \mu(\mathcal{F})$ (resp. $\mu(\mathcal{G}) < \mu(\mathcal{F})$). The full subcategory ${\cC}_{\mu}$
of $Coh(\E)$ consisting of all semi-stable sheaves of a fixed
slope $\mu \in \mathbb{Q} \cup \{\infty\}$ is abelian, artinian and closed under
extensions. Moreover, if $\mathcal{F}, \mathcal{G}$ are semi-stable with $\mu(\mathcal{F}) < \mu(\mathcal{G})$ then
$\Hom(\mathcal{G},\mathcal{F})=\Ext(\mathcal{F},\mathcal{G})=0$. Any sheaf $\mathcal{F}$ possesses a unique filtration
(the Harder-Narasimhan filtration, or HN filtration)
$$0=\mathcal{F}^{r+1} \subset \mathcal{F}^{r} \subset \cdots \subset \mathcal{F}^1=\mathcal{F}$$
for which $\mathcal{F}^i/\mathcal{F}^{i+1}$ is semi-stable of slope, say
$\mu_i$, and $\mu_1 < \cdots  <\mu_r$. Observe that ${\cC}_{\infty}$
is just the category of torsion sheaves, and hence is equivalent to the
product category $\prod_{x } \mathcal{T}or_x$, where $x$ runs through the
set of closed points of $\E$  and  $\mathcal{T}or_x$ denotes the
category of torsion sheaves supported at $x$. Since $\mathcal{T}or_x$
is equivalent to the category of finite length modules over the local ring $R_x$ of the point $x$,
 there is a unique
simple sheaf $\kO_{x}$ in $\mathcal{T}or_x$.

\begin{theo}[\cite{A}]\label{T:Ati} The following holds~:
\begin{enumerate}
\item[i)] the HN filtration of any coherent sheaf
splits (non-canonically). In particular,  any indecomposable coherent sheaf is semi-stable,
\item[ii)] the set of stable sheaves of slope $\mu$ is the set of simple objects
of ${\cC}_{\mu}$,
\item[iii)] there are canonical exact equivalences of abelian categories
$\epsilon_{\nu, \mu}:{\cC}_{\mu}
\stackrel{\sim}{\lto} {\cC}_{\nu}$ for any $\mu, \nu \in \mathbb{Q}
\cup \{\infty\}$.
\end{enumerate}
\end{theo}

\vspace{.1in}

The Grothendieck group $K_0\bigl(Coh(\E)\bigr)$ of $Coh(\E)$ is equipped with the Euler bilinear form
$
\langle\; \,,\;\rangle~:\; K_0\bigl(Coh(\E)\bigr) \otimes K_0\bigl(Coh(\E)\bigr) \lto \mathbb{Z}
$ defined by the formula
$$
\overline{\mathcal{F}} \otimes \overline{\mathcal{G}} \mapsto \mathrm{dim}\;\Hom(\mathcal{F},\mathcal{G}) - \mathrm{dim}\;\Ext(\mathcal{F},\mathcal{G}).
$$
There is a natural map
$
K_0\bigl(Coh(\E)\bigr) \lto K_0'\bigl(Coh(\E)\bigr):=\mathbb{Z}^2,
$
given by
$$
\overline{\mathcal{F}} \mapsto \bigl(\rank(\mathcal{F}), \deg(\mathcal{F})\bigr)
$$
whose kernel coincides with the radical of the form $\langle\; \,,\;\rangle$.
As we shall be mainly interested in the class of a sheaf in the numerical $K$--group
 $K_0'\bigl(Coh(\E)\bigr)$, we also
denote by $\overline{\mathcal{F}}$ the pair
$\bigl(\rank(\mathcal{F}),\deg(\mathcal{F})\bigr)$.
By the Riemann-Roch formula one has
$$\bigl\langle (r_1,d_1),(r_2,d_2) \bigr\rangle=r_1d_2-r_2d_1.$$
In particular, the Euler form is skew-symmetric in our case.

\vspace{.2in}

\paragraph{\textbf{1.2.}} Let $D^b\bigl(Coh({\E})\bigr)$ stand for  the bounded derived
category of coherent sheaves on $\E$. As $Coh(\E)$ has  global dimension
one, the structure of  $D^b\bigl(Coh({\E})\bigr)$ is very simple to describe:
any object of this category is isomorphic to its cohomology,
i.e.~$\mathcal{F}^\bullet \simeq \bigoplus_n H^n(\mathcal{F}^\bullet)[-n]$.

We also consider the so-called root category  $\mathcal{R}_{\E}=D^b\bigl(Coh({\E})\bigr)/[2]$,
where $[1]$ is the shift functor. This category can be described  as follows
\begin{enumerate}
\item $\mathrm{Ob}(\mathcal{R}_{\E}) = \bigl\{\kF^\pm| \kF \in \mathrm{Ob}(Coh_\E)\bigr\}$
\item $\mathrm{Hom}_{\mathcal{R}_\E}(\kF^\pm, \kG^\pm) = \mathrm{Hom}_\E(\kF,\kG)$ and
$\mathrm{Hom}_{\mathcal{R}_\E}(\kF^\pm, \kG^\mp) = \mathrm{Ext}_\E^1(\kF,\kG)$.
\end{enumerate}

\medskip
\noindent
The category $\mathcal{R}_\E$ is triangulated and there is a canonical exact functor
$$\Psi:D^b\bigl(Coh({\E})\bigr) \lto \mathcal{R}_\E$$
 inducing a group isomorphism $K_0(\E) \lto K_0(\mathcal{R}_{\E})$.
Since the shift $[2]$ preserves the Euler form
$\langle\;\,,\;\rangle$, we can define a morphism
$K_0(\mathcal{R}_{\E}) \lto K_0'\bigl(Coh(\E)\bigr)$, mapping $\overline{\kF^\pm}$ to the class $\pm \overline\kF$.
Moreover, one can view the root category $\mathcal{R}_\E$ as the category of two-periodic complexes with
the functor $\Psi$ being a  Galois  covering functor in the sense of Gabriel,
see \cite{X2per} for further details.

Next, let us consider auto-equivalences of triangulated categories
$D^b\bigl(Coh({\E})\bigr)$ and
$\mathcal{R}_{\E}$.
Let $\kE$ be a \emph{spherical object} in  the derived category $D^b\bigl(Coh({\E})\bigr)$,
i.e.~an object
satisfying
$\mathrm{Hom}(\kE, \kE) = \mathrm{Hom}\bigl(\kE, \kE[1]\bigr) = \kk$. For example the
structure sheaf the curve $\kO$ or
the structure
sheaf of a $\kk$-rational point $\kO_{x_0}$. Seidel and Thomas considered in
\cite{ST} the functor
$$
T_\kE: D^b(Coh_\E) \lto D^b(Coh_\E)
$$
defined by
$ T_\kE(\kF) = \mathrm{cone}\bigl({\mathrm{RHom}}(\kE, \kF)\overset{\kk}\otimes \kE \stackrel{ev}\lto \kF\bigr).
$
The functor $T_\kE$  is exact and if coherent sheaves $\kE$ and $\kF$  satisfy the condition
 $\mathrm{Ext}^1(\kE,\kF) = 0$, then
$T_\kE(\kF)$ is quasi-isomorphic to the complex $$\bigl(\mathrm{Hom}(\kE,\kF) \overset{\kk}\otimes
\kE \stackrel{ev}\lto
\kF\bigr) = \bigl(\kE^n \stackrel{ev}\lto \kF\bigr),$$
where $n = \mathrm{dim} \,\mathrm{Hom}(\kE, \kF)$. On the level of $K_0\bigl(Coh(\E)\bigr)$ the functor
$T_\kE$ induces the group homomorphism
 $t_\kE: K_0\bigl(Coh(\E)\bigr) \lto K_0\bigl(Coh(\E)\bigr)$,  given by
$$
\gamma \mapsto  \gamma - \langle \kE, \gamma\rangle\overline{\kE},
$$
where $\langle\; \,\,,\;\rangle$ denotes  the Euler form on $K_0\bigl(Coh(\E)\bigr)$.

Let $x_0$ be a rational point of $\E$.
In the basis $\bigl\{ \overline{\kO}, \overline{\kO}_{x_0}\bigr\}$ of the numerical $K$-group
$K'_0\bigl(Coh(\E)\bigr)$, the twist functors
$T_\kO$, $T_{\kO_{x_0}}$ and the shift $[1]$  induce linear transformations given by the  matrices
$$t_{\mathcal{O}}=\left(\begin{matrix} 1 & -1\\
0&1\end{matrix} \right), \qquad
t_{\kO_{x_0}}=\left( \begin{matrix} 1 & 0 \\
1 & 1 \end{matrix} \right), \qquad
t_{[1]}=\left(\begin{matrix} -1 & 0
\\ 0 & -1 \end{matrix} \right).$$
Observe that for any $\kk$-rational point  $x_0$
the equivalence $T_{\kO_{x_0}}$ preserves $Coh(\E)$ and is simply
given by $\mathcal{F} \mapsto \mathcal{F} \otimes \mathcal{O}(x_0)$, see
\cite[formula (3.11)]{ST}.

Due to  \cite[Proposition 2.10]{ST} the functor $T_\kE$ is an equivalence of
categories for any spherical object $\kE$ and by
\cite[Lemma 3.2]{ST} it
is isomorphic to a   Fourier-Mukai transform   with the kernel
$\mathrm{cone}(\kE^\vee \boxtimes \kE \lto \kO_\Delta) \in D^b\bigl(Coh(\E \times \E)\bigr)$.
Moreover, by   \cite[Proposition 2.13]{ST} we have the following braid group relation:
$$
T_{\kO_{x_0}} T_\kO T_{\kO_{x_0}} \cong T_\kO T_{\kO_{x_0}} T_\kO.
$$

\begin{prop}[see \cite{Muk,ST}]\label{Mukai}
Let $\Phi := T_{\kO_{x_0}} T_\kO T_{\kO_{x_0}}$, then $\Phi^2 \cong i^*[1]$, where $i$ is an involution
of $\E$ preserving $x_0$.  Moreover,  for the duality functor  $D =  {R}{\mathcal Hom}(-,\kO)$
we have an isomorphism
$$
D \circ \Phi \cong i^* \circ [1] \circ \Phi \circ D.
$$
\end{prop}

\noindent
\textit{Proof}. The braid group relation between $T_\kO$ and $T_{\kO_{x_0}}$ was proven in \cite{ST}
without any restrictions on the base field. However, in the proof of two other isomorphisms, given in
\cite{Muk} the assumption for  $\kk$ to be algebraically closed was used.
We refer to Appendix A for a  proof in the case of an arbitrary field.

\vspace{.15in}

From the above relations one deduces that the group
generated by $T_{\mathcal{O}}, T_{\mathcal{O}_{x_0}}$ and $[1]$ is the universal
covering  $\widetilde{SL}(2,\Z)$ of $SL(2,\Z)$ given by a  central  extension of $SL(2,\Z)$ by
$\Z$. Since in ${\Aut}(\mathcal{R}_{\E})$ we have $[1]^2 \simeq
\id$, the action of the group $\bigl\langle T_\kO, T_{\mathcal{O}_{x_0}}, [1]\bigr\rangle$ on
the root category
$\mathcal{R}_{\E}$
breaks up
to the action of  $\widehat{SL}(2,\Z)$, where  $\widehat{SL}(2,\Z)$ is a two-fold covering of
$SL(2,\Z)$.
That all may be summed up in the following commutative diagram:
$$\xymatrix{
\widetilde{SL}(2,\Z) \ar@{->>}[d] \ar@{^{(}->}[r] & \Aut\bigl(D^b\bigl(Coh({\E})\bigr)\bigr) \ar[d]\\
\widehat{SL}(2,\Z) \ar@{^{(}->}[r] \ar@{->>}[rd] & {\Aut}(\mathcal{R}_{\E}) \ar@{->>}[d]\\
& SL(2,\Z) = {\Aut}\bigl(K'_0(Coh(\E))\bigr)}
$$

\vspace{.05in}

For any $\nu \in \mathbb{Q} \cup \{\pm \infty\}$ denote by
$Coh_{\leq \nu}$ (resp. $Coh_{> \nu}$) the full subcategory
of $Coh(\E)$ consisting of sheaves all of whose indecomposable (=
semi-stable)
constituents have slope at most $\nu$ (resp. strictly greater then $\nu$).
Next, let $Coh^{\nu}(\E)$ be the full subcategory of $D^b\bigl(Coh({\E})\bigr)$
whose objects consist of direct sums $\mathcal{F} \oplus \mathcal{G}[1]$
where $\mathcal{F} \in Coh_{> \nu} , \mathcal{G} \in Coh_{\leq \nu}$.
This has the structure of an
abelian category as the heart of the t-structure
on $D^b\bigl(Coh({\E})\bigr)$ associated to the torsion pair $(Coh_{>\nu}, Coh_{\leq \nu})$.
One can view the category $Coh^\nu(\E)$ as a full
subcategory of the root category $\mathcal{R}_{\E}$.

\vspace{.05in}

For a spherical  sheaf $\mathcal{E}$ of class $(r,d) \in K_0'\bigl(Coh(\E)\bigr)$ and slope $\mu=\frac{d}{r}$
the auto-equivalence $T_{\mathcal{E}}$ establishes an equivalence
between   $Coh(\E)$ and  $Coh^{\nu}(\E)$, where $\nu=-\infty$ if
$\mu=\infty$ and $\nu=\mu-\frac{1}{r^2}$ if $\mu \neq \infty$. More generally, if
$\hat{\gamma} \in \widehat{SL}(2,\Z)$ is a lift of $\gamma \in SL(2,\Z)$ then $\hat{\gamma}$
sends $Coh(\E)$ to $Coh^\nu(\E)$ where $\nu=\frac{p'}{q'}$, and $(q',p')=\gamma(0,-1)$. Finally,
each equivalence $\epsilon_{\nu,\mu}$ in Atiyah's Theorem \ref{T:Ati}
can  be obtained as the restriction to ${\cC}_{\mu}$ of one of the
above auto-equivalences of $D^b\bigl(Coh({\E})\bigr)$ and $\mathcal{R}_{\E}$.
We can visualize the structure of the category $\mathcal{R}_{\E}$  by the following picture, where
$Coh(\E)^+ = Coh(\E)$ and $Coh(\E)^- = Coh(\E)[1]$.

\centerline{
\begin{picture}(36,16)
\thicklines
\put(7,7){\vector(0,1){7}}
\thinlines
\put(12.2,6){$\rank$}
\put(7,7){\line(0,-1){7}}
\put(7,7){\vector(1,0){7}}
\put(7,14.2){$\deg$}
\put(7,7){\line(-1,0){7}}
\put(7,7){\line(1,2){3}}
\put(10.2,13){$\mathcal{F}$}
\put(7,7){\line(-2,-1){4}}
\put(3,4.2){$\mathcal{G}[1]$}
\put(3,5){\circle*{.2}}
\qbezier(8,7)(7.9,7.5)(7.5,8)
\put(8.1,7.5){$\mu(\mathcal{F})$}
\put(10,10){$Coh(\E)^+$}
\put(1,10){$Coh(\E)^-$}
\qbezier(5.4,7)(5.5,6.7)(5.7,6.35)
\put(2.8,6){$\mu(\mathcal{G})$}
\put(10,13){\circle*{.2}}
\put(16,6.5){\huge{$\stackrel{\hat{\gamma}}{\longrightarrow}$}}
\put(28,7){\vector(0,1){7}}
\put(28,7){\line(0,-1){7}}
\put(28,7){\vector(1,0){7}}
\put(28,7){\line(-1,0){7}}
\put(28,14.2){$\deg$}
\put(33.2,6){$\rank$}
\thicklines
\put(28,7){\line(-1,1){7}}
\thinlines
\put(28,7){\line(1,-1){7}}
\qbezier(29.5,7)(29.35,6.5)(29,6)
\put(29.4,6.2){$\nu$}
\put(28.5,12){$Coh^{\nu}(\E)$}
\put(21.5,2){$Coh^{\nu}(\E)[1]$}
\put(32,7.5){$Coh_{> \nu}$}
\put(23,12){$Coh_{\leq \nu}[1]$}
\end{picture}}

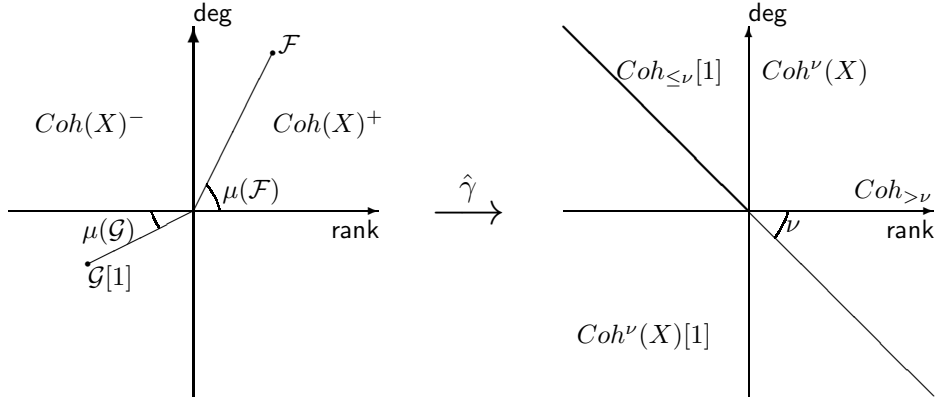
\centerline{Figure 1. The root category $\mathcal{R}_{\E}$ and its auto-equivalences}

\vspace{.2in}

\section{Hall algebra of an elliptic curve}

\vspace{.1in}

\paragraph{\textbf{2.1.}} From now on we assume that $\kk=\mathbb{F}_q$ is a
finite field,
fix a square root $v$ of $q^{-1}$ and  work over the quadratic field extension
$\KK = \Q(\sqrt{q})=\Q(v)$.
Note that $Coh(\E)$ is a hereditary abelian category.
Consider the free $\KK$-module $\mathbf{H}_{\E}$
with linear basis $\bigl\{[\mathcal{F}]\bigr\}$ where $\mathcal{F}$ runs through the set
 of isomorphism classes of objects in
$Coh(\E)$. There is a natural $\Z^2$-grading on $\mathbf{H}_{\E}$ given by
$\mathbf{H}_{\E}[\a]=\bigoplus_{\overline{\mathcal{F}}=\a} \KK [\mathcal{F}]$.
To a triple $(\mathcal{F},\mathcal{G},\mathcal{H})$ of coherent sheaves we associate the finite set
$\mathcal{P}_{\mathcal{F},\mathcal{G}}^\mathcal{H}$ of exact sequences
$0 \to \kG \to \kH \to \kF \to 0$. Next, we set 
${P}_{\mathcal{F},\mathcal{G}}^\mathcal{H}= \# \mathcal{P}_{\mathcal{F},\mathcal{G}}^\mathcal{H}$
and ${F}_{\mathcal{F},\mathcal{G}}^\mathcal{H}= \frac{\displaystyle {P}_{\mathcal{F},\mathcal{G}}^\mathcal{H}}{\displaystyle a_{\mathcal{F}}
a_{\mathcal{G}}}$, where  $a_{\mathcal{K}}=
\#\mathrm{Aut}(\mathcal{K})$ for a coherent sheaf  $\mathcal{K}$. As in \cite{Ri} we now define an associative product
on $\mathbf{H}_{\E}$ by the formula
\begin{equation}\label{E:Hallprod}
[\mathcal{F}] \cdot [\mathcal{G}]=v^{-\langle \mathcal{F},\mathcal{G} \rangle}
\sum_{\mathcal{H}} {F}_{\mathcal{F},\mathcal{G}}^\mathcal{H} [\mathcal{H}],
\end{equation}
and, following \cite{Green}, a coassociative coproduct
\begin{equation}\label{E:Hallcoprod}
\Delta\bigl([\mathcal{H}]\bigr)=\sum_{\mathcal{F},\mathcal{G}} v^{-\langle \mathcal{F},\mathcal{G}\rangle}\frac{{P}_{\mathcal{F},\mathcal{G}}^\mathcal{H}}{a_{\mathcal{H}}}
 [\mathcal{F}] \otimes [\mathcal{G}].
\end{equation}
(note that we are using the opposite of the algebra and coalgebra structures considered in \cite{Kap}).
The counit $\varepsilon: \H_\E \lto \KK$ is defined as follows
$$
\varepsilon\bigl([\kF]\bigr) =
\left\{
\begin{array}{ccc}
1 & \textrm{\,if\,} & \kF \cong 0 \\
0 & \textrm{\,if\,} & \kF \not\cong 0. \\
\end{array}
\right.
$$
Finally, the bilinear form given  by
$$\bigl([\mathcal{F}],[\mathcal{G}]\bigr)=\delta_{\mathcal{F},\mathcal{G}}
\frac{1}{a_{\mathcal{F}}}$$
is a non-degenerate Hopf pairing on ${\mathbf{H}}_{\E}$, i.e.~we have
$(ab,c) = \bigl(a \otimes b, {\Delta}(c)\bigr)$ for any $a,b,c \in {\mathbf{H}}_{\E}$ (see \cite{Green}).

\vspace{.2in}

\paragraph{\textbf{2.2.}} The comultiplication $\Delta$ only takes value
 in a certain completion of $\H_{\E}
\otimes \H_{\E}$ (the sum on the right-hand side of (\ref{E:Hallcoprod}) is infinite unless $\mathcal{H}$ is a torsion sheaf).
Note also that  the space
 $\mathbf{H}_{\E}[\a]$ is infinite dimensional for $\a = (r,d) \in \mathbb{Z}^2,  r>0$.

We denote  $(\Z^2)^+ = \bigl\{(q,p) \in \Z^2\;|\; q \geq 1\; \text{or}\; q=0, p \geq 0\bigr\}$
and for a given class $\a \in  (\mathbb{Z}^2)^+$
define $\H_{\E}^{\not\ge m}[\a] = \span \bigl\{[\kF]| \overline{\kF} = \a \,\,\mbox{and} \,\,
\kF \not\in Coh_{\ge m}\bigr\}$
and $\H_{\E}^{\ge m}[\a] = \span\bigl\{[\kF]| \overline{\kF} = \a \,\,
\mbox{and}\,\,
\kF \in Coh_{\ge m}\bigr\}$.

\begin{lem}\label{L:freqused}
For any class $\a \in (\mathbb{Z}^2)^+$ and any integer $m$ the vector space $\H_{\E}^{\ge m}[\a]$
is finite-dimensional.
\end{lem}

\noindent
\textit{Proof}. Note that for any $m\in \mathbb{Z}$  there are only
finitely many elements
$\a_1,\a_2,\dots,\a_t$ of $(\mathbb{Z}^2)^+$ such that $m \le  \mu(\a_1) < \dots < \mu(\a_t)$ and
$\a_1 + \a_2 + \dots  + \a_t = \a$. Moreover, it follows from the Atiyah's classification that for any
class $\beta \in (\mathbb{Z}^2)^+$ there are only finitely many semi-stable coherent sheaves of class
$\beta$. Since any coherent sheaf on an elliptic curve splits into a direct sum of semi-stable ones, the claim
easily follows.
\qed

\medskip
For any integer $m$ we have a surjective linear map of vector spaces $\jet_m: \H_\E[\a] \rightarrow
\H_\E^{\ge m}[\a]$ inducing an isomorphism $
\pi_m: \H_{\E}[\a]/\H_{\E}^{\not\ge m}[\a] \rightarrow \H_{\E}^{\ge m}[\a].
$
For any $m \le n$  the canonical embedding
$\H_{\E}^{\not\ge m}[\a] \rightarrow  \H_{\E}^{\not\ge n}[\a]$ induces a commutative diagram

\begin{tabular}{p{2.5cm}c}
&
\xymatrix
{ \H_{\E}[\a]/\H_{\E}^{\not\ge n}[\a] \ar[rrr]^{\pi_n} & & & \H_{\E}^{\ge n}[\a] \\
\H_{\E}[\a]/\H_{\E}^{\not\ge m}[\a]\ar[u] \ar[rrr]^{\pi_m} & & & \H_{\E}^{\ge m}[\a]
\ar[u]_{\varphi_{m,n}} \\
}
\end{tabular}

\noindent
Obviously, $\bigl(\H_{\E}^{\ge n}[\a], \varphi_{m,n}\bigr)$ forms a projective system, and we can define
\begin{equation}
\widehat{\H}_{\E}[\a] := \lim_{\underset{n}{\longleftarrow}}\bigl(\H_{\E}^{\ge n}[\a]\bigr).
\end{equation}

\noindent
One can view $\widehat{\H}_{\E}[\a]$  as the set of infinite sums
$\bigl\{\sum a_{\kF} [\kF]| a_{\kF} \in \KK,
\overline{\kF} = \a\bigr\}$.
For the sake of convenience we also denote by $\jet_n$
the canonical morphism $\widehat{\H}_{\E}[\a] \rightarrow  \H_{\E}^{\ge n}[\a]$.  By the universal property of the projective limit there is   an (injective)  linear map $\H_\E[\alpha] \rightarrow \widehat{\H}_\E[\alpha]$, and
 since the surjection
$\H_\E[\a] \rightarrow  \H_\E^{\not\ge n}[\a]$ splits, we may  consider  $\H_\E^{\not\ge n}[\a]$ as a subspace of
$\widehat{\H}_{\E}[\a]$ via the inclusion $\H_\E^{\not\ge n}[\a] \rightarrow  \H_\E[\alpha]
\lto \widehat{\H}_{\E}[\a]$.
So, the projection  $\jet_n: \widehat{\H}_{\E}[\a] \rightarrow \H_{\E}^{\ge n}[\a]$ is an idempotent morphism and if we
 denote
$r_n = 1 - \jet_n$, then any element $h \in \widehat{\H}_{\E}[\a]$ can be written as $\jet_n(h) + r_n(h)$, where
$\jet_n(h) \in \H_{\E}^{\ge n}[\a]$ and $\jet_n\bigl(r_n(h)\bigr) = 0$.
Using this formalism, the space $\H_{\E}[\a]$ viewed as a subset of  $\widehat{\H}_{\E}[\a]$ can be identified
 with the set  of  those sequences
$h = (h_n)$ for which $r_n(h_n) = 0$ for
$n \gg 0$.

\medskip

So, we
define $\widehat{\H}_{\E} := \bigoplus\limits_{\a \in (\mathbb{Z}^2)^+}
\widehat{\H}_{\E}[\a]$.
In a similar way, for $\a, \b \in (\mathbb{Z}^2)^+$ the sequence of vector spaces
$\bigl(\H_{\E}^{\ge n}[\a] \otimes \H_{\E}^{\ge m}[\beta]\bigr) =
\bigl(\H_{\E}[\a]/\H_{\E}^{\not\ge n}[\a] \otimes \H_{\E}[\b]/\H_{\E}^{\not\ge n}[\b]\bigr)$
forms a projective system and we put
\begin{equation}
\H_{\E}[\a] \widehat{\otimes}  \H_{\E}[\beta] :=
\lim_{\underset{n,m}{\longleftarrow}}(\H_{\E}^{\ge n}[\a] \otimes \H_{\E}^{\ge m}[\beta]).
\end{equation}
In this case as well
$
\H_{\E}[\a] \widehat{\otimes}  \H_{\E}[\beta]$ can be identified with the set of infinite sums
$ \bigl\{
\sum\limits_{\kF, \kG} b_{\kF, \kG} [\kF]\otimes [\kG]| \overline{\kF} = \alpha,
\overline{\kG} = \beta, \;
b_{\kF, \kG} \in \KK\bigr\}.
$
For $\gamma \in (\mathbb{Z}^2)^+$ we set
\begin{equation}
(\H_{\E} \widehat{\otimes}  \H_{\E})[\gamma] :=
\prod\limits_{\substack{\alpha +\beta = \gamma \\ \alpha, \beta \in (\mathbb{Z}^2)^+}}
\H_{\E}[\a] \widehat{\otimes}  \H_{\E}[\beta]
\end{equation}
and finally
\begin{equation}
\H_{\E} \widehat{\otimes}  \H_{\E} := \bigoplus\limits_{\gamma \in (\mathbb{Z}^2)^+}
\H_{\E} \widehat{\otimes}  \H_{\E}[\gamma].
\end{equation}

\begin{prop}\label{P:topbial}
In
 the notation as above the following properties hold
\begin{enumerate}
\item  $\widehat{\H}_{\E}$ and $\H_{\E} \widehat{\otimes}  \H_{\E}$ are associative algebras;
\item  the comultiplication
$\Delta: \H_\E \lto \H_{\E} \widehat{\otimes}  \H_{\E}$ is a ring homomorphism  and extends to a map
$\Delta: \widehat{\H}_\E \lto \H_{\E} \widehat{\otimes}  \H_{\E}$;
\item  let
$\Delta_{\a,\b}: \H_{\E}[\a+\b] \to \H_{\E}[\a] \widehat{\otimes} \H_{\E}[\b]$ stand for the
$(\a,\beta)$-component of $\Delta$, then
 $\Delta_{\a,\b}\bigl(\H_{\E}[\a+\b]\bigr) \subset \H_{\E}[\a]
\otimes \H_{\E}[\b]$.
\end{enumerate}
\end{prop}

\noindent
\textit{Proof.} Let us show that the composition map
$
\widehat{\H}_{\E}[\alpha] \otimes \widehat{\H}_{\E}[\beta] \stackrel{m}\lto
\widehat{\H}_{\E}[\alpha +\beta]
$  given by the rule
$\bigl(\sum a_{\kH}[\kH]\bigr) \otimes \bigl(\sum b_\kG [\kG]\bigr) \mapsto
\bigl(\sum a_\kH b_\kG [\kH][\kG]\bigr)$
is well-defined. Indeed,  for a fixed  coherent sheaf $\kF$ of class
$\overline\kF = \alpha + \beta$ there are finitely many exact sequences
$$
0 \lto \kG \lto \kF \lto \kH \lto 0
$$
such that $\overline\kH = \alpha$ and $\overline\kG = \beta$. To see this,  let
$\kF = \bigoplus\limits_{i = 1}^n \kF_i$ and
$\kH = \bigoplus\limits_{j = 1}^m \kH_j$ be
the  splittings of $\kF$ and $\kH$ into a direct sum of semi-stable objects, then the existence of an epimorphism
$\kF \twoheadrightarrow \kH$ implies the conditions
 $\rank(\kH) \le \rank(\kF)$ and
$\mu(\kH_j) \ge  \min\bigl\{\mu(\kF_i)| 1 \le i \le n \bigr\}$  for all
$1 \le j \le m$. Hence it follows that the degrees of all sheaves  $\kH_j$ are bounded below and
as $\sum\limits_{j=1}^m \deg(\kH_j) = \deg(\alpha)$,  they are also bounded above. By Atiyah's classification, there
are finitely many semi-stable sheaves of a given class and hence there are finitely many sheaves $\kH$ of class
$\alpha$ which are quotients of $\kF$. In the same way, there are only finitely many subsheaves of $\kF$ of class
$\beta$. This means that only finitely many sheaves from $\H_\E[\a]$ and  $\H_\E[\b]$ contribute
to the element $[\kF]$ from $\widehat{\H}_\E[\a+b]$, which shows that
the map  $m$ is well-defined.

\medskip
\noindent
In a similar fashion,  one  deals  with $\H_{\E} \widehat\otimes \H_{\E}$.
In this case
the map
\begin{equation}
\begin{split}
\prod\limits_{\alpha_1 + \beta_1 = \gamma_1}\hspace{-.1in}\H_{\E}[\alpha_1]
\widehat{\otimes}
\H_{\E}[\beta_1] \hspace{.1in} \otimes
\prod\limits_{\alpha_2 + \beta_2 = \gamma_1}\hspace{-.1in}\H_{\E}[\alpha_2]
\widehat{\otimes}
\H_{\E}[\beta_2] &\\ \stackrel{m}\lto
\prod_{\substack{\a_2 + \b_2 = \gamma_2 \\ \a_1 + \b_1 = \gamma_1}}\hspace{-.1in}
\H_{\E}[\alpha_1 + \alpha_2] \widehat{\otimes}
\H_{\E}[\beta_1 + \beta_2]
\end{split}
\end{equation}
 is convergent since for a given
$[\kF] \otimes [\kG] \in \H_{\E}[\gamma_1] \otimes
\H_{\E}[\gamma_2]$ there are finitely many surjective morphisms
$\kF \twoheadrightarrow \kM$ and $\kG \twoheadrightarrow \kN$, where $\kM$ and $\kN$ are
coherent sheaves  satisfying
$\overline\kM + \overline\kN = \gamma_1$.
The proof that $\Delta_{\a,\b}(\H_\E)[\a + \b] \subseteq \H_\E[\a] \otimes  \H_\E[\b]$
is completely analogous.

To see that $\Delta$ is a ring homomorphism, fix a pair of
tuples
$(\a, \b, \a', \b')$ and $(\gamma,\gamma',\delta,\delta')$ of elements of $K_0'(Coh(\E))$ satisfying
$$\gamma+\gamma'=\a,\quad \delta+\delta'=\b,\quad \gamma+\delta=\a',\quad \gamma'+\delta'=\beta'$$
and put
$$c_{\underline{\gamma}, \underline{\delta}}=(m \otimes m) \circ P_{23}
\circ (\Delta_{\gamma,\gamma'} \otimes \Delta_{\delta,\delta'}):
\H_{\E}[\a] \otimes \H_{\E}[\b] \to \H_{\E}[\a'] \otimes \H_{\E}[\b'],$$
where $P_{23}$ is the operator of permutation of the second and third
components. For any tuples of sheaves $(\mathcal{F},
\mathcal{G},\mathcal{H},\mathcal{K})$ such that $\overline{\mathcal{F}}=\a,
\overline{\mathcal{G}}=\b, \overline{\mathcal{H}}=\a', \overline{\mathcal{K}}=\b'$
let $c_{\underline{\gamma},\underline{\delta}}(\mathcal{F}, \mathcal{G},\mathcal{H},\mathcal{K})$
be the coefficient of
$[\mathcal{H}] \otimes [\mathcal{K}]$ in $c_{\underline{\gamma},\underline{\delta}}([\mathcal{F}]
\otimes [\mathcal{G}])$. It is easy to see that
$c_{\underline{\gamma},
\underline{\delta}}(\mathcal{F}, \mathcal{G},\mathcal{H},\mathcal{K})=0$ for all but finitely many tuples
$(\gamma,\gamma',\delta,\delta')$.

\begin{lem}\cite{Green}\label{L:Green}
The map
$$(m \circ \Delta)_{\a,\b}^{\a',\b'}=\bigoplus_{\underline{\gamma},\underline{\delta}} c_{\underline{\gamma},
\underline{\delta}}: \H_{\E}[\a] \otimes \H_{\E}[\b] \to
\H_{\E}[\a'] \otimes \H_{\E}[\b'].$$
satisfies the equality  $(m \circ \Delta)_{\a,\b}^{\a',\b'}=\Delta_{\a',\b'} \circ m$.
\end{lem}

\noindent
\emph{Note to the proof of Lemma}.
As in the case of quivers,  this result is equivalent to the following formula. Let $\kF, \kG, \kM, \kN$
be arbitrary coherent sheaves on $\E$. Then the following equality of Hall numbers is true:
\begin{equation}
\sum\limits_{\kH}\frac{P^{\kH}_{\kM, \kN} \cdot P^{\kH}_{\kF, \kG}}{a_\kH} =
\sum\limits_{\kA, \kB,\kC, \kD} q^{-\langle \kA, \kD\rangle}
\frac{P^{\kM}_{\kA, \kB} \cdot P^{\kN}_{\kC, \kD} \cdot P^{\kF}_{\kA, \kC} \cdot P^{\kG}_{\kB, \kD}}{a_\kA
\cdot a_\kB \cdot a_\kC \cdot a_\kD}.
\end{equation}
This  formula  can be  proved by essentially the same computation as in \cite{Green}
(a more detailed
proof in \cite{RingelGreen},
in which all arguments involving the dimension vector are replaced by the corresponding ones
involving $K_0'\bigl(Coh(\E)\bigr)$, can be
applied in our case literally). See also \cite{SLectures}.
\qed

Observe that since the Euler form $\langle \; \,,\;\rangle$ is antisymmetric it is not necessary
to twist the
multiplication in ${\mathbf{H}}_{\E} \otimes{\mathbf{H}}_{\E}$ as it was done in \cite{Green}.
Lemma \ref{L:Green} implies that the linear map $\Delta$ is a ring
homomorphism. The Proposition \ref{P:topbial} is proven.
\qed

\vspace{.1in}

A graded algebra with a graded coproduct satisfying the properties of
Proposition \ref{P:topbial} will be called a \textit{topological} bialgebra.
The next important lemma  says   that
the linear maps
$$m: \widehat{\H}_\E[\a] \otimes \widehat{\H}_\E[\b] \lto \widehat{\H}_\E[\a + \b],
\hspace{1cm} \Delta_{\a,\b}:  \widehat{\H}_\E[\a +\b] \lto \H_\E[\a] \widehat\otimes \H_\E[\b]$$
and
$$m:  (\H_\E[\a_1] \widehat\otimes \H_\E[\b_1]) \otimes (\H_\E[\a_2] \widehat\otimes \H_\E[\b_2])
\lto \H_\E[\a_1 + \b_1] \widehat\otimes \H_\E[a_2 + \b_2]$$
are \emph{continuous}. Recall that
for an element
$a \in \widehat{\H}_\E[\a]$ and $n\in \Z$ we can
 write $a = \jet_n(a) + r_n(a)$, where $\jet_n(a) \in \widehat{\H}_\E^{\ge n}$ and
$\jet_n\bigl((r_n(a)\bigr) = 0$.

\begin{lem}\label{L:cont}
For any two classes $\a,\b \in (\mathbb{Z}^2)^+$ and any $m \in \mathbb{Z}$
there exists another integer $n$ such that for any
$a \in \widehat{\H}_\E[\a]$ and $b\in \widehat{\H}[\b]$ we have
$\jet_m(ab) =
\jet_m\bigl(\jet_n(a) \jet_n(b)\bigr)$.
Similarly,  for any  pair of integers $ m,n$ there exists another pair
$k,l$ such that for all elements
$f \in \H_\E[\a_1] \widehat\otimes \H_\E[\b_1]$,
$g \in \H_\E[\a_2] \widehat\otimes \H_\E[\b_2]$  we have
$$\jet_{m,n}(fg) = \jet_{m,n}\bigl(\jet_{k,l}(f) \jet_{k,l}(g)\bigr).
$$
Finally, for any pair of integers $m,n$  there exists
$k$ such that for and any
$a \in \widehat\H_\E[\a + \b]$ we have
$\jet_{m,n}\bigl(\Delta_{\a, \b}(a)\bigr) =
\Delta_{\a, \b}\bigl(\jet_k(a)\bigr)$.
\end{lem}

\noindent
\textit{Proof}.
For any coherent sheaf $\kH$ of class  $\a +\b$  there are only finitely many
sheaves $\kF$ of class $\alpha$ such that there is a surjection $\kH \twoheadrightarrow \kF$. Hence,
we have  a finite number of exact sequences $0 \lto \kG \lto \kH \lto \kF \lto 0$ with
$\overline\kF = \a$ and $\overline\kG = \b$.
Since the vector space $\H_\E^{\ge m}[\a + \b]$ is  finite-dimensional, we see that
there exists $n$ such that
$r_n(a)$ and $r_n(b)$ do  not contribute to $\jet_m(ab)$.
The proof of two other statements is
completely analogous.
\qed

\medskip
\noindent
Later, we shall need the following
property of the Hopf pairing in $\H_\E$.

\begin{lem}\label{R:hopfp}
Let $\sum x_i' \otimes x''_i \in \H_\E \widehat\otimes \H_\E[\gamma]$ and $y \in \H_\E[\gamma]$ and suppose that
$\sum (x_i' x''_i, y) < \infty$. Then
$$
\sum\limits_i  (x_i' x''_i, y) =  \sum\limits_{i,j} (x'_i, y_j^{(1)})(x''_i, y_j^{(2)})
$$
where $\sum\limits_{j} y_j^{(1)} \otimes y_j^{(2)} = \Delta(y)$.
\end{lem}

\vspace{.1in}

\paragraph{\textbf{2.4.}} There exists a natural ``PBW-type'' decomposition for $\mathbf{H}_{\E}$.
For any $\mu \in \Q \cup \{\infty\}$ we consider the subspace
${\H}_{\E}^{(\mu)}  \subset \H_{\E}$ linearly spanned by classes
$\bigl\{[\mathcal{F}]\;|\;\mathcal{F} \in {\cC}_{\mu}\bigr\}$. Since the category
${\cC}_{\mu}$ is stable under extensions, ${\H}_{\E}^{(\mu)}$ is a
subalgebra of $\H_{\E}$ (but not a subbialgebra~!). The exact equivalence
$\epsilon_{\mu_1,\mu_2}$ defined in Theorem~\ref{T:Ati} gives rise to an
algebra isomorphism $\epsilon_{\mu_1,\mu_2}:~{\H}_{\E}^{(\mu_2)}
\stackrel{\sim}{\to} {\H}_{\E}^{(\mu_1)}$.
Let $\vec{\bigotimes\limits_{\mu}} \mathbf{H}^{(\mu)}_{\E}$ stand for the (restricted) tensor product of spaces $\mathbf{H}^{(\mu)}_{\E}$
with $\mu \in \mathbb{Q} \cup \{\infty\}$, ordered from left to right in increasing order, i.e.~for the vector space spanned by
elements of the form $a_{\mu_1} \otimes \cdots \otimes a_{\mu_r}$ with $a_{\mu_i} \in \mathbf{H}_{\E}^{(\mu_i)}$ and
$\mu_1 < \cdots < \mu_r$.

\begin{lem}\label{L:PBW} The multiplication map
 $m:\; \vec{\bigotimes\limits_{\mu}} \mathbf{H}^{(\mu)}_{\E}  \to
{\mathbf{H}}_{\E}$ is an isomorphism.\end{lem}

\noindent
\textit{Proof.} As the spaces $\mathrm{Ext}(\mathcal{F},\mathcal{G})$ vanish for
 $\mathcal{F} \in {\cC}_{\mu}$,
$\mathcal{G}_{\nu}\in {\cC}_{\nu}$ and $\mu < \nu$, we have, up to a power of $v$,
$[\mathcal{F}_1] \cdot [\mathcal{F}_2] \cdots [\mathcal{F}_r]=[\mathcal{F}_1 \oplus \cdots \oplus \mathcal{F}_r]$ if $\mathcal{F}_i
\in {\cC}_{\mu_i}$ and $\mu_1 < \ldots < \mu_r$. Any sheaf can  be decomposed into a direct sum of semi-stable summands, and these
are determined up to isomorphism. The statement easily follows.\qed

\vspace{.1in}

Let $\cC[\mu_1,\mu_2]$ be the full subcategory of sheaves whose HN decomposition only contains slopes $\mu \in [\mu_1, \mu_2]$.
This category is exact  and in particular stable under
extensions. Moreover, we have the following remark.

\vspace{.1in}

\addtocounter{theo}{1}
\paragraph{\textbf{Remark \thetheo}} For any $\mu_1\leq \mu_2$ the Hall algebra of the exact category
$\cC[\mu_1,\mu_2]$ is a subalgebra of $\mathbf{H}_{\E}$, isomorphic to
$\vec{\bigotimes\limits_{\mu_1 \leq \mu \leq \mu_2}} \mathbf{H}_{\E}^{(\mu)}$.

\vspace{.1in}
\noindent
We conclude this section by the following proposition.

\begin{prop}\label{P:Reineke}
 Consider the algebra
 $$
 T = \KK\langle X, Y^\pm\rangle/(Y^\pm Y^\mp = 1, \; XY^\pm  = v^{\pm 2}  Y^\pm X).$$
  Then there exists  a
$\KK$--linear
algebra homomorphism $\chi: \mathbf{H}_{\E} \to T$, called Reineke's character,  given by the formula
$$\chi\bigl([\kF]\bigr) = q^{-\a \b} \frac{X^\a Y^\b}{a_\kF},$$ where $(\a, \b) = \bigl(\rank(\kF), \deg(\kF)\bigr) \in \mathbb{Z}^2$.
\end{prop}

\noindent
\emph{Proof}. Let  $
\Gamma  \subset K_0'\bigl(Coh(\E)\bigr)$ be the semi-group
generated by the images of classes of coherent
sheaves on $\E$.  Following Reineke \cite{Reineke},  consider the associative algebra
$$
\KK\bigl(\Gamma, \langle\;\, , \;\rangle\bigr) =
\Bigl\{\sum_{\gamma \in \Gamma} a_\gamma t^\gamma | a_\gamma \in \KK \Bigr\}
$$
where the multiplication is given by the rule
$t^\alpha t^\beta = v^{\langle \alpha, \beta\rangle} t^{\alpha + \beta}$. Let $X = t^{\overline\kO}$ and
$Y^\pm = t^{\pm \overline{\kO}_{x_0}}$. Then we have:
 $Y^\pm Y^\mp = 1$ and
$t^{\overline\kO + \overline{\kO}_{x_0}} = v YX = v^{-1}XY,$ hence
$\KK\bigl(\Gamma, \langle\;\, , \;\rangle\bigr) \cong T$.
Finally,  by  \cite[Lemma 6.1]{Reineke} the linear map  $\chi:\mathbf{H}_{\E} \to T$ mapping $[\kF]$ to
$\frac{\displaystyle t^{\overline\kF}}{\displaystyle a_\kF} = q^{-\a \b}
\frac{\displaystyle X^\a Y^\b}{\displaystyle a_\kF}$ is an algebra homomorphism. \qed

\section{Drinfeld double of ${\mathbf{H}}_\E$}

\vspace{.1in}

\paragraph{\textbf{3.1.}} As in the case of quivers, it is natural to consider the Drinfeld double of the bialgebra
${\mathbf{H}}_{\E}$. This is what we do in this Section.

\begin{lem}\label{L:relH}
$\mathbf{H}_{\E}$ is isomorphic to the $\KK$-algebra generated by the collection of elements $\{x_{\mathcal{F}}\;|\mathcal{F}\;
\text{is\;semi-stable}\;\}$ subject to the set of relations
\begin{equation}\label{E:rel0}
x_{\mathcal{F}} \cdot x_{\mathcal{G}}=v^{-\langle \mathcal{F},\mathcal{G}\rangle} \sum_{\mathcal{H}}
{F}_{\mathcal{F},\mathcal{G}}^{\mathcal{H}} \underline{x}_{\mathcal{H}}, \qquad \forall\; \mathcal{F},\mathcal{G}\;
\text{semi-stable}
\end{equation}
where by definition $\underline{x}_{\mathcal{H}}=v^{\sum_{i<j}\langle \mathcal{H}_i,\mathcal{H}_{j}\rangle} x_{\mathcal{H}_1}
\cdots x_{\mathcal{H}_r}$ if $\mathcal{H}=\mathcal{H}_1 \oplus \cdots \oplus \mathcal{H}_r$ with all $\mathcal{H}_i$ being 
semi-stable and $\mu(\mathcal{H}_1) < \cdots < \mu(\mathcal{H}_r)$.
\end{lem}

\noindent
\textit{Proof.} Let $\mathbf{G}$ be the algebra defined above. By construction, there is a morphism $\phi: \mathbf{G} \to \mathbf{H}_{\E}$,
 which is surjective by virtue of Lemma~\ref{L:PBW} (we have $\phi(\underline{x}_{\mathcal{H}})=[\mathcal{H}]$).
Let $\mathbf{G}' \subset \mathbf{G}$ denote the linear span of elements $\underline{x}_{\mathcal{H}}$ for $\mathcal{H} \in Coh(\E)$.
It is clear that $\phi$ restricts to an isomorphism of vector spaces between $\mathbf{G}'$ and $\mathbf{H}_{\E}$, hence it is enough
to show that $\mathbf{G}=\mathbf{G}'$.

\vspace{.05in}

If $\mathcal{F}=\mathcal{H}_1 \oplus \cdots \oplus \mathcal{H}_r$ is a
decomposition of a sheaf $\mathcal{F}$ into a direct sum of semi-stable
objects  with
$\mu(\mathcal{H}_1) < \cdots < \mu(\mathcal{H}_r)$, we denote
$HN(\mathcal{F})=(\overline{\mathcal{H}_1}, \ldots,
\overline{\mathcal{H}_r})$ and call this vector the \textit{HN-type} of
$\mathcal{F}$. One can introduce an order on the set of HN types as follows :
$\bigl((r_1,d_1), \ldots,  (r_s,d_s)\bigr) \preceq \bigl((r'_1,d'_1), \ldots,
(r'_t,d'_t)\bigr)$ if there exists $l$ such that $(r_{s-i},d_{s-i})=(r'_{t-i},d'_{t-i})$ for $i <l$ while
$\frac{d_{s-l}}{r_{s-l}} > \frac{d'_{t-l}}{r'_{t-l}}$ or
$\frac{d_{s-l}}{r_{s-l}} = \frac{d'_{t-l}}{r'_{t-l}}$ and $d_{s-l} >
d'_{t-l}$.

Fix $\a \in K_0'\bigl(Coh(\E)\bigr)$. We shall  prove that any monomial $x_{\mathcal{F}_1} \cdots x_{\mathcal{F}_r}$ of weight $\a$ belongs to
$\mathbf{G}'$.
 For this, we argue successively by induction on the HN type $HN(\underline{\mathcal{F}})$ of the
sheaf $\mathcal{F}=\mathcal{F}_1\oplus \cdots \oplus \mathcal{F}_r$  and then on
the number $n_{\underline{\mathcal{F}}}$ of inversions in the sequence
$\bigl(\mu(\mathcal{F}_1), \ldots, \mu(\mathcal{F}_r)\bigr)$. Note that
if $HN(\underline{\mathcal{F}})$ is maximal, i.e.~if $\mu(\mathcal{F}_1)=\cdots=\mu(\mathcal{F}_r)=\nu$ then
$x_{\mathcal{F}_1} \cdots x_{\mathcal{F}_r} \in \bigoplus_{\mathcal{H} \in {\cC}_{\nu}}
\KK x_{\mathcal{H}} \subset \mathbf{G}'$; on the other hand, if $n_{\underline{\mathcal{F}}}=0$ then $\mu(\mathcal{F}_1)
\leq \cdots \leq \mu(\mathcal{F}_r)$ and $x_{\mathcal{F}_1} \cdots x_{\mathcal{F}_r} \in \mathbf{G}'$ by definition. So
let $x_{\mathcal{F}_1} \cdots x_{\mathcal{F}_r}$ be a monomial of weight $\a$ and assume that $x_{\mathcal{G}_1} \cdots x_{\mathcal{G}_s}$
belongs to $\mathbf{G}'$ whenever $HN(\underline{\mathcal{G}}) \succ HN(\underline{\mathcal{F}})$ or
$HN(\underline{\mathcal{G}}) = HN(\underline{\mathcal{F}})$ and $n_{\underline{\mathcal{G}}} < n_{\underline{\mathcal{F}}}$. If
$n_{\underline{\mathcal{F}}}=0$ then we are done, so we may assume that $\mu(\mathcal{F}_i) > \mu(\mathcal{F}_{i+1})$ for some $i$.
By Remark~2.7, we have
$$x_{\mathcal{F}_{i}} \cdot x_{\mathcal{F}_{i+1}} \in \KK x_{\mathcal{F}_{i+1}}\cdot x_{\mathcal{F}_i}
\oplus \bigoplus_{\substack{\mathcal{H} \in {\cC}[\mu(\mathcal{F}_{i+1}), \mu(\mathcal{F}_i)]\\ \mathcal{H} \neq \mathcal{F}_{i} \oplus \mathcal{F}_{i+1}}} \KK \underline{x}_{\mathcal{H}}.$$
Now observe that the number of inversions of $x_{\mathcal{F}_1} \cdots x_{\mathcal{F}_{i+1}} \cdot x_{\mathcal{F}_i}
\cdots x_{\mathcal{F}_r}$ is one less than $n_{\underline{\mathcal{F}}}$, while the HN-type the sheaf $\mathcal{F}_1 \oplus \cdots
\oplus \mathcal{F}_{i-1} \oplus \mathcal{H} \oplus \mathcal{F}_{i+2} \cdots \oplus \mathcal{F}_{r}$ is strictly greater
than $HN(\underline{\mathcal{F}})$ as soon as $\mathcal{H} \in {\cC}\bigl[\mu(\mathcal{F}_{i+1}),\mu(\mathcal{F}_i)\bigr]$ is
of class $\overline{\mathcal{F}_{i} \oplus \mathcal{F}_{i+1}}$ and $\mathcal{H} \neq \mathcal{F}_i \oplus \mathcal{F}_{i+1}$ .
We deduce using the induction hypothesis that $x_{\mathcal{F}_1} \cdots x_{\mathcal{F}_r}$ belongs
to $\mathbf{G}'$, as desired.\qed

\vspace{.2in}

\paragraph{} Let $\mathbf{D}{\mathbf{H}}_{\E}$ be the Drinfeld double of the topological bialgebra
${\mathbf{H}}_{\E}$ with respect to the Hopf
pairing $(\;\,,\;)$. Recall (see e.g.~\cite{XDrinfeld}) that this is an algebra generated by two copies of ${\mathbf{H}}_{\E}$,
which we denote by ${\mathbf{H}}^+_{\E}$ and ${\mathbf{H}}^-_{\E}$ to avoid confusion, subject to the following set of
relations  for any pair $g \in {\mathbf{H}}^+_\E$ and $ h \in {\mathbf{H}}^-_\E$~:
\begin{equation}\label{E:Rgh}
\sum_{i,j} h_i^{(1)\,-} g_j^{(2)\,+} \bigl(h_i^{(2)},g_j^{(1)}\bigr)=\sum_{i,j} g_j^{(1)\,+}h_i^{(2)\,-} \bigl(h_i^{(1)},g_j^{(2)}\bigr)
\tag{$R(g,h)$}
\end{equation}
(we use here the usual Sweedler notation ${\Delta}(x^\pm)=\sum_i x_i^{(1)\, \pm}
\otimes x_i^{(2)\, \pm}$).
Observe that although the coproduct takes value in a completion of ${\mathbf{H}}_{\E} \otimes {\mathbf{H}}_{\E}$,
the relation $\bigl(R(g,h)\bigr)$ contains only finitely many terms. Indeed it is enough to consider the case $g=[\mathcal{G}], h=[\mathcal{H}]$,
and then $h_i^{(2)}$ involves only sheaves which are subsheaves of $\mathcal{H}$, while $g_j^{(1)}$ involves only sheaves which are
quotients of $\mathcal{G}$. As $\mathrm{Hom}(\mathcal{G},\mathcal{H})$ is a finite set, there are only finitely many sheaves which
are both quotients of $\mathcal{G}$ and subsheaves of $\mathcal{H}$, hence the scalar product $(h_i^{(2)},g_j^{(1)})$ vanishes for
almost all values of $(i,j)$. The same holds for the right-hand side of $\bigl(R(g,h)\bigr)$.
If $h \in {\mathbf{H}}_{\E}$ then we write $h^+, h^-$ for the corresponding elements in ${\mathbf{H}}^+_{\E}$
and ${\mathbf{H}}^-_{\E}$ respectively.

\vspace{.1in}

\begin{prop}\label{P:relDH} The algebra $\mathbf{D}{\mathbf{H}}_{\E}$ is isomorphic to the $\KK$-algebra generated by two copies
${\mathbf{H}}^+_\E, {\mathbf{H}}^-_{\E}$ of the Hall algebra ${\mathbf{H}}_\E$
subject to the set of relations
\begin{equation}\label{E:rel22}
R\bigl([\mathcal{G}]^+,[\mathcal{H}]^-\bigr)\qquad \text{\;for\;any\;semi-stable\;} \mathcal{G}, \mathcal{H} \in Coh(\E).
\end{equation}
\end{prop}

\noindent
\textit{Proof.}
We have to show that the set of relations (\ref{E:rel22}) for semistable $a=[\mathcal{G}]^+, b=[\mathcal{H}]^-$ implies the set of
relations $R(a,b)$ for arbitrary $a,b$. By bilinearity of the relations $R(a,b)$ it is enough to prove this
in the case $a=[\mathcal{F}]^+, b=[\mathcal{K}]^-$ for some (arbitrary) sheaves $\mathcal{F},\mathcal{K}$.

\begin{lem}\label{L:relprod} Let $a,b \in {\mathbf{H}}^+_\E, c,d \in {\mathbf{H}}^-_{\E}$. The relation $R(ab,c)$ is
implied by the collection of all relations $R(a,c_k^{(1)})$ and $R(b,c_k^{(2)})$ for all $k \ge 1$. Similarly, $R(a,cd)$ follows from the collection of relations $R(a_k^{(1)},c)$ and
$R(a_k^{(2)},d)$.
\end{lem}

\noindent
We refer to Appendix B for a proof of this lemma.

\medskip

Now, let us consider the algebra $\mathbf{A}$ generated by ${\mathbf{H}}^+_\E$ and
${\mathbf{H}}^-_{\E}$ modulo relations (\ref{E:rel22}).
For any coherent sheaf $\mathcal{F}$ there exist  semi-stable sheaves
$\mathcal{G}_1, \ldots, \mathcal{G}_r$ such that $[\mathcal{F}] = v^{\sum_{i<j}\langle \mathcal{G}_i,
\mathcal{G}_j \rangle} [\mathcal{G}_1] \cdots [\mathcal{G}_r]$. Thus, in view of the above Lemma~\ref{L:relprod}, it is enough to prove
that $R\bigl([\mathcal{G}],[\mathcal{K}]\bigr)$ holds for semi-stable $\mathcal{G}$ and arbitrary $\mathcal{K}$. We shall  prove this
 by induction on the rank $r$ of $\mathcal{K}$. As any torsion sheaf is semi-stable, the statement is clear for $r=0$.
So let us assume that $R\bigl([\mathcal{G}],[\mathcal{K}']\bigr)$ holds for all semi-stable $\mathcal{G}$ and arbitrary
$\mathcal{K}'$ of rank less than $r$, and let $\mathcal{K}$ be a sheaf of rank $r$. If $\mathcal{K}$ is semi-stable then
there is nothing to prove, so we may assume that $\mathcal{K}$ splits into a non-trivial direct sum of semi-stable objects
$\mathcal{K}_1 \oplus \cdots \oplus \mathcal{K}_l$. Assume first that $\rank(\mathcal{K}_i) < r$ for all $i$.
Then $[\mathcal{K}] = v^{\sum_{i<j}\langle \mathcal{K}_i, \mathcal{K}_j \rangle}[\mathcal{K}_1] \cdots [\mathcal{K}_l]$, and by
Lemma~\ref{L:relprod} $R\bigl([\mathcal{G}],[\mathcal{K}]\bigr)$ is a consequence of the set of relations $R\bigl([\mathcal{G}]^{(i)}_j,[\mathcal{K}_i]\bigr)$.
 These hold in $\mathbf{A}$ by the induction hypothesis since $\rank(\mathcal{K}_i) <r$. The last case to consider is that of a sum
$\mathcal{K}=\mathcal{I} \oplus \mathcal{T}$ where $\mathcal{I}$ is a semi-stable vector bundle and $\mathcal{T}$ is a torsion sheaf.
As above, $R\bigl([\mathcal{G}],[\mathcal{K}]\bigr)$ is implied by the relations $R\bigl([\mathcal{G}]^{(1)}_i,[\mathcal{I}]\bigr)$ and
$R\bigl([\mathcal{G}]^{(2)}_j,[\mathcal{T}]\bigr)$. The second set of relations is satisfied by the induction hypothesis.
For the first set, let us again decompose $[\mathcal{G}]^{(1)}_i=v^{d_i} [\mathcal{V}_1]
\cdots [\mathcal{V}_t]$ for some semi-stable sheaves $\mathcal{V}_j$. As before, it is enough to see that
$R\bigl([\mathcal{V}_j],[\mathcal{I}]^{(j)}_k\bigr)$ holds for all $j,k$. But as $\mathcal{I}$ is a vector bundle, any sheaf appearing in
$[\mathcal{I}]^{(j)}_k$ is either semi-stable, or  splits as a direct sum of smaller rank sheaves. In both cases the induction
hypothesis applies. The Proposition is proved.\qed

\vspace{.1in}

\begin{prop}\label{P:DDoub} The multiplication map
$\mathbf{H}^+_{\E} \otimes \mathbf{H}^-_{\E} \stackrel{m}\lto
\mathbf{D}\mathbf{H}_{\E}$ is a vector space isomorphism.
\end{prop}

\noindent
\textit{Proof.} This statement is classical for Hopf algebras (see \cite{Joseph}, 3.2.4).
However, in our situation of topological bialgebras, extra
care needs to be taken because the coproduct $\Delta$ takes values in the  completion
$\H_{\E} \widehat\otimes \H_{\E}$. A proof of Proposition \ref{P:DDoub}
is given in Appendix B.
\qed

\vspace{.2in}

\paragraph{\textbf{3.2.}} It is useful to view $\mathbf{D}{\mathbf{H}}_{\E}$  as the (yet inexistent)
 Hall algebra of the root category $\mathcal{R}_{\E}$,
 where
${\mathbf{H}}^+_{\E}$ corresponds to the Hall algebra of $Coh(\E)$
and ${\mathbf{H}}^-_{\E}$ corresponds to
the Hall algebra of $Coh(\E)[1]$ (see, however \cite{Toen} or \cite{XX} for a recent approach to Hall algebras for derived 
categories). For  $\mathcal{F} \in Coh(\E)$ we put
$$\bigl[\mathcal{F}[\epsilon]\bigr]=\begin{cases} [\mathcal{F}]^+& \;\text{if}\; \epsilon =0\\
[\mathcal{F}]^-&\;\text{if}\;\epsilon =1. \end{cases}$$
We define the set of semi-stable objects of the root category
$\mathcal{R}_{\E}$ as
$\bigl\{\mathcal{F}[\epsilon]\bigr\}$, where  $\mathcal{F}$ is semi-stable and
$\epsilon \in \Z/2\Z$. Observe
that this set is invariant under auto-equivalences of $\mathcal{R}_{\E}$.

\begin{cor} The algebra $\mathbf{D}{\mathbf{H}}_{\E}$ is generated by the set of elements $[\mathcal{F}]$, where $\mathcal{F}$ runs among all semi-stable objects $\mathcal{F} \in \mathcal{R}_{\E}$.\end{cor}

\noindent
\textit{Proof.} This is a consequence of Lemma~\ref{L:relH} and Proposition~\ref{P:relDH}.\qed

\vspace{.2in}

Similarly to the case of the usual Hall numbers,
for any triple of objects $\kF, \kG, \kH$ of the derived category $D^b\bigl(Coh(\E)\bigr)$
we denote by $P^\kH_{\kF, \kG}$ the number of the distinguished triangles
$\bigl\{\kG \rightarrow \kH \rightarrow \kF \rightarrow \kG[1]
\bigr\}$ and $F^\kH_{\kF, \kG} = \frac{\displaystyle P^\kH_{\kF, \kG}}{\displaystyle a_\kF \cdot a_\kG}$. Next,
for any four objects $\kM, \kN, \kA$ and  $\kB$ of $Coh(\E)$ we denote by $C^{\kM, \kN}_{\kA, \kB}$
the number of the long exact sequences of the form
$$
\bigl\{0 \rightarrow \kN \rightarrow \kB \rightarrow \kA
\rightarrow \kM \rightarrow 0 \bigr\}.
$$
The following result of Kapranov \cite[Lemma 2.4.3]{Kap2}
plays a key role in our study of the Drinfeld double
$\mathbf{D}{\mathbf{H}}_{\E}$.

\begin{lem}\label{L:formulofKapr}
For for any four objects $\kM, \kN, \kA$ and  $\kB$ of $Coh(\E)$ we have:
$$
C^{\kM, \kN}_{\kA, \kB} =
\frac{\displaystyle P^{\kN[1] \oplus \kM}_{\kB[1], \kA}}{\displaystyle |\Ext(\kM, \kN)|}.
$$
\end{lem}

\noindent
Our next goal is to obtain an explicit form of the relations $R\bigl([\kF]^-, [\kG]^+\bigr)$, where
$\kF$ and $\kG$ are semi-stable sheaves on $\E$.

\begin{prop}\label{P:relationsDD}
Let $\kF$ and $\kG$ be a pair of semi-stable sheaves on an elliptic curve
 $\E$ with slopes  $\mu = \mu(\kF)$ and $\nu = \mu(\kG)$.
\begin{enumerate}
\item\label{it1ofProponDD} If $\mu < \nu$ then $R\bigl([\kF]^-, [\kG]^+\bigr)$ can be rewritten as
$$
[\kF]^- \cdot [\kG]^+ = v^{\langle \kF, \kG \rangle}\,
\sum\limits_{\kB, \kC} \, v^{- \langle \kC, \kB \rangle} \, F^{\kC[1] \oplus \kB}_{\kF[1], \kG} \,
[\kC]^+ \cdot [\kB]^-.
$$
\item\label{it2ofProponDD} If $\mu > \nu$ then $R\bigl([\kF]^-, [\kG]^+\bigr)$ reads  as
$$
[\kG]^+ \cdot [\kF]^- = v^{\langle \kG, \kF \rangle} \,
\sum\limits_{\kA, \kD} v^{- \langle \kD, \kA \rangle} \, F^{\kD \oplus \kA[-1]}_{\kG, \kF[-1]} \,
[\kD]^- \cdot [\kA]^+.
$$
\item\label{it3ofProponDD} Finally, if $\mu = \nu$ then we have
$$
\sum\limits_{\kA, \kD \in \cC_\mu} \, C^{\kA, \kD}_{\kF, \kG} \, [\kA]^+ \cdot [\kD]^- =
\sum\limits_{\kB, \kC \in \cC_\mu} \, C^{\kC, \kB}_{\kG, \kF} \, [\kC]^- \cdot [\kB]^-.
$$
\end{enumerate}
\end{prop}

\noindent\textit{Proof}. (\ref{it1ofProponDD}) Consider the first case when $\mu < \nu$.  Let
$$
\Delta\bigl([\kF]\bigr) =
\sum\limits_{\kA, \kB} v^{- \langle \kA, \kB\rangle} \frac{P^{\kF}_{\kA, \kB}}{a_\kF} [\kA] \otimes
[\kB] \quad \mbox{and} \quad
\Delta\bigl([\kG]\bigr) =
\sum\limits_{\kC, \kD} v^{- \langle \kC, \kD\rangle} \frac{P^{\kG}_{\kC, \kD}}{a_\kG} [\kC] \otimes
[\kD].
$$
Since $\kF$ and $\kG$ are semi-stable and $\mu(\kF) < \mu(\kG)$, we have
$\bigl([\kB], [\kC]\bigr) = 0$
for any proper subobject $\kB$ of $\kF$ and any proper quotient object $\kC$ of $\kG$. Hence, the relation $R\bigl([\kF]^-, [\kG]^+\bigr)$ has the following
shape:
$$
[\kF]^- \cdot [\kG]^+ =
\sum\limits_{\kA, \kB, \kC, \kD}
v^{- \langle \kA, \kB\rangle - \langle \kC, \kD \rangle} \frac{P^{\kF}_{\kA, \kB}
 P^{\kG}_{\kC, \kD}}{a_\kF a_\kG} \bigl([\kA], [\kD]\bigr) \; [\kC]^+ \cdot [\kB]^-
$$
$$
= \sum\limits_{\kB, \kC} v^{-\langle \kF -\kB, \kB \rangle - \langle \kC, \kG - \kC\rangle}
\frac{1}{a_\kF a_\kG} \sum\limits_{\kA} \frac{P^{\kF}_{\kA, \kB} P^{\kG}_{\kC, \kA}}{a_\kF} \;
[\kC]^+ \cdot [\kB]^-.
$$
Recall that the Euler form on $Coh(\E)$  is skew-symmetric, hence for any object
$\mathcal{I}$ of $Coh(\E)$  we have: $\langle \mathcal{I}, \mathcal{I}\rangle = 0$. Next, note  the following equality
of Hall coefficients:
$$
\sum\limits_{\kA} \frac{P^{\kF}_{\kA, \kB} P^{\kG}_{\kC, \kA}}{a_\kF} = C^{\kC, \kB}_{\kG, \kF}.
$$
Hence, the whole expression can be rewritten as
$$
[\kF]^- \cdot [\kG]^+ =  \sum\limits_{\kB, \kC} v^{-\langle \kF, \kB\rangle - \langle \kC, \kG\rangle} \,
\frac{C^{\kC, \kB}_{\kG, \kF}}{a_\kF a_\kG} \; [\kC]^+ \cdot [\kB]^-.
$$
Since $\Hom(\kC, \kB) = 0$ for any subobject $\kB$ of $\kF$ and any  quotient object
$\kC$ of $\kG$, Lemma \ref{L:formulofKapr}  implies that
$C^{\kC, \kB}_{\kG, \kF} = v^{-2 \langle \kC, \kB\rangle} P^{\kC[1] \oplus \kB}_{\kF[1], \kG}$. Hence, we
obtain:
$$
[\kF]^- \cdot [\kG]^+ = \sum\limits_{\kB, \kC}
v^{-\langle \kF, \kB \rangle - \langle \kC, \kG\rangle - 2\langle \kC, \kB\rangle}
\, F^{\kC[1] \oplus \kB}_{\kF[1], \kG} \, [\kC]^+ \cdot  [\kB]^-.
$$
To get the claim, it remains to note that
$$
\langle \kF, \kG\rangle + \langle \kF, \kB\rangle +
\langle \kC, \kG \rangle + \langle \kC, \kB \rangle =
\langle \kF + \kC, \kG + \kB\rangle = 0.
$$

\vspace{1mm}
\noindent
(\ref{it2ofProponDD}) In the case $\mu > \nu$, the derivation of the formula for
$R\bigl([\kF]^-, [\kG]^+\bigr)$  is similar to the case
(\ref{it1ofProponDD}) and is therefore left to the reader.

\vspace{1mm}
\noindent
(\ref{it3ofProponDD}) Finally, consider the case $\mu(\kF) = \mu = \mu(\kG)$.
Let sequences
$$
0 \longrightarrow \kB \longrightarrow \kF \longrightarrow \kA \longrightarrow 0
\quad \mbox{and} \quad
0 \longrightarrow \kD \longrightarrow \kG \longrightarrow \kC \longrightarrow 0
$$
be exact. Assume that $\bigl([\kB], [\kC]\bigr) \ne 0$, i.e.~$\kB \cong \kC$.
Then $\kB$ and $\kC$ are necessarily semi-stable of slope $\mu$. Hence, $\kA$ and
$\kD$ are semi-stable of slope $\mu$ as well. In other words, all four
objects $\kA, \kB, \kC$ and $\kD$ belong to the same abelian category
$\cC_\mu$. The relation  $R\bigl([\kF]^-, [\kG]^+ \bigr)$ can be rewritten as follows:
\begin{equation*}
\begin{split}
\sum\limits_{\substack{\kA, \kB, \kC, \kD \in \cC_\mu \\ \kB \cong \kC}}
v^{-\langle \kA, \kB \rangle - \langle \kC, \kD \rangle}
\, \frac{P^\kF_{\kA, \kB} P^\kG_{\kC, \kD}}{a_\kF a_\kG a_\kB} \; [\kA]^+ \cdot  [\kD]^-= \\
= \sum\limits_{\substack{\kA, \kB, \kC, \kD \in \cC_\mu \\ \kA \cong \kD}}
v^{-\langle \kA, \kB \rangle - \langle \kC, \kD \rangle}
\, \frac{P^\kF_{\kA, \kB} P^\kG_{\kC, \kD}}{a_\kF a_\kG a_\kA} \; [\kC]^-  \cdot [\kB]^+.
\end{split}
\end{equation*}
By Theorem \ref{T:Ati},  $\cC_\mu$ is equivalent to the category of coherent torsion sheaves
on $\E$. Hence,  the Euler form $\langle \;\, , \; \rangle$ vanishes on $\cC_\mu$ and 
 we obtain
the relation
\begin{equation*}
\sum\limits_{\kA, \kD \in \cC_\mu} \Bigl(\sum\limits_{\kB \in \cC_\mu}
\frac{P^\kF_{\kA, \kB} P^\kG_{\kB, \kD}}{a_\kB}\Bigr) \; [\kA]^+ \cdot [\kD]^- =
\sum\limits_{\kB, \kC \in \cC_\mu} \Bigl(\sum\limits_{\kA \in \cC_\mu}
\frac{P^\kF_{\kA, \kB} P^\kG_{\kC, \kA}}{a_\kA}\Bigr) \; [\kC]^- \cdot [\kB]^+,
\end{equation*}
which is obviously equivalent to the relation (\ref{it3ofProponDD}) of Proposition
\ref{P:relationsDD}.
\qed

\begin{theo}\label{P:braidact} Let $\Phi$ be an auto-equivalence of
$D^b\bigl(Coh(\E)\bigr)$. Then the  assignment $[\mathcal{F}] \mapsto [\Phi(\mathcal{F})]$, where
$\mathcal{F}$ is a semi-stable object of the root category  $\mathcal{R}_{\E}$, extends to a uniquely determined  algebra automorphism of $\mathbf{D}{\mathbf{H}}_{\E}$.
\end{theo}

\noindent
\textit{Proof}. Recall that  $\mathbf{D}{\mathbf{H}}_{\E}$ is a $\KK$-algebra
generated by the symbols $[\kF]^\pm$, where $\kF$ is a semi-stable coherent sheaf on $\E$
subject to the relations $P\bigl([\kF]^\pm, [\kG]^\pm\bigr)$
\begin{equation*}
[\kF]^\pm \cdot [\kG]^\pm = v^{- \langle \kF, \kG\rangle}
\sum\limits_{\kH \cong \kH_1 \oplus \kH_2 \oplus \dots \oplus \kH_t}
v^{\sum_{i < j} \langle \kH_i, \kH_j \rangle} \, F^{\kH}_{\kF, \kG} \, [\kH_1]^\pm \dots [\kH_t]^\pm,
\end{equation*}
where $\kF$ and $\kG$ are semi-stable, $\mu(\kF) < \mu(\kG)$ and $\kH = \kH_1 \oplus \kH_2 \oplus \dots \oplus \kH_t$ is a splitting into a direct sum of semi-stable objects such that
$\mu(\kH_1) < \mu(\kH_2) < \dots < \mu(\kH_t)$; together
with the relations $R\bigl([\kF]^\pm, [\kG]^\pm\bigr)$ of Proposition \ref{P:relationsDD}. In order to show that the group
$\Aut\bigl(D^b(Coh(\E))\bigr)$ acts on $\mathbf{D}{\mathbf{H}}_{\E}$ by algebra automorphisms, it is sufficient to check that all relations $P\bigl([\kF]^\pm, [\kG]^\pm\bigr)$ and
$R\bigl([\kF]^-, [\kG]^+\bigr)$ are preserved for
$\kF$ and $\kG$ semi-stable.

\vspace{2mm}
\noindent
Consider first the case of the relations $P\bigl([\kF], [\kG]\bigr)$. Let
$\mu = \mu(\kF), \nu = \mu(\kG)$ and $\Phi$ be an auto-equivalence
of $D^b\bigl(Coh(\E)\bigr)$.

\vspace{2mm}
\noindent
\underline{Case 1}. First assume that $\Phi(\kF) \cong \widehat\kF[i]$ and
$\Phi(\kG) \cong \widehat\kG[i]$ for some $i \in \mathbb{Z}$.
Let $\hat\mu = \mu(\widehat\kF)$ and $\hat\nu = \mu(\widehat\kG)$. If we assume
that $\mu > \nu$ then it automatically follows that  $\hat\mu > \hat\nu$. Moreover, $\Phi$ induces
an equivalence of exact categories $\cC[\nu, \mu] \rightarrow \cC[\hat\nu, \hat\mu]$. Hence, we
have an isomorphism of Hall algebras $H\bigl(\cC[\nu, \mu]\bigr)  \rightarrow
H\bigl(\cC[\hat\nu, \hat\mu]\bigr)$ preserving all Hall constants. In other words, the
relation  $P\bigl([\kF]^\pm, [\kG]^\pm\bigr)$ is  mapped to the relation
$P\bigl([\widehat\kF]^\pm, [\widehat\kG]^\pm\bigr)$.

\vspace{2mm}
\noindent
\underline{Case 2}. Assume $\Phi(\kF) \cong \widehat\kF[2i+1]$ and
$\Phi(\kG) \cong \widehat\kG[2i]$ for some $i \in \mathbb{Z}$.

\centerline{
\begin{picture}(14,16)
\thicklines
\put(7,7){\vector(0,1){7}}
\thinlines
\put(12.2,6){$\rank$}
\put(7,7){\line(0,-1){7}}
\put(7,7){\vector(1,0){7}}
\put(7,14.2){$\deg$}
\put(7,7){\line(-1,0){7}}
\put(7,7){\line(1,2){3}}
\put(10.2,13){$\mathcal{F}$}
\put(10,3){$Coh(\E)$}
\put(1,3){$Coh(\E)[1]$}
\put(10,13){\circle*{.2}}
\put(7,7){\line(-2,1){6}}
\put(1,10){\circle*{.2}}
\put(0,10.3){$\Phi(\mathcal{F})$}
\put(7,7){\line(4,-1){4}}
\put(11,6.2){$\mathcal{G}$}
\put(11,6){\circle*{.2}}
\put(7,7){\line(1,4){1}}
\put(8,11){\circle*{.2}}
\put(7.3,11.5){$\Phi(\mathcal{G})$}
\end{picture}}
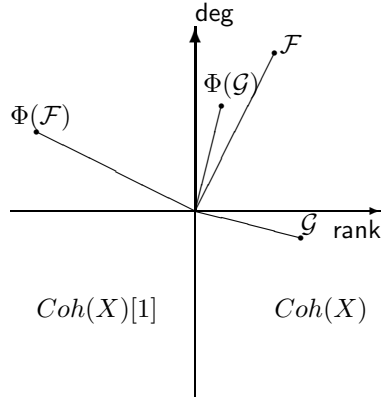
\centerline{Figure 2. Relation $P([\kF], [\kG])$, where $\mu(\kF) > \mu(\kG)$}

\vspace{2mm}
\noindent
First note that there exists
a slope $\kappa$, where $\nu \le \kappa < \mu$ such that $\Phi(\cC_\varphi) \in Coh(\E)[2i+1]$
for $\kappa < \varphi \le  \mu$ and  $\Phi(\cC_\varphi) \in Coh(\E)[2i]$ for $\nu \le
\varphi \le \kappa$. Next,
for any short exact sequence $0 \rightarrow \kG \rightarrow \kH \rightarrow \kF \rightarrow
0$ we can write 
$$\kH \cong \bigl(\kH_1 \oplus \dots \oplus \kH_t\bigr) \oplus
\bigl(\kH_{t+1} \oplus \dots \oplus \kH_n\bigr) \cong \kH' \oplus \kH'',
$$
where all objects $\kH_i$ are semi-stable, $\mu(\kH_1) < \mu(\kH_2) < \dots <
\mu(\kH_t) = \kappa < \kH_{t+1} < \dots < \mu(\kH_n)$,
$\kH' = \kH_1 \oplus \dots \oplus \kH_t$ and  $\kH'' = \kH_{t+1} \oplus \dots \oplus \kH_n$. In these
notations we have:
\begin{equation}\label{E:actofAut1}
[\kF]*[\kG] = v^{-\langle \kF, \kG\rangle}
\sum\limits_{\kH', \kH''} F^{\kH' \oplus \kH''}_{\kF, \kG} v^{\langle\kH', \kH''\rangle}
[\kH'] * [\kH''].
\end{equation}
Since $\Phi(\kF) \cong \widehat{\kF}[2i+1], \Phi(\kH'') \cong \widehat{\kH}''[2i+1]$, whereas
$\Phi(\kG) \cong \widehat{\kG}[2i]$ and $\Phi(\kH') \cong \widehat{\kH}'[2i]$, we have:
$\langle \kF, \kG\rangle = - \langle \widehat\kF, \widehat\kG\rangle$ and
$\langle \kH', \kH''\rangle = - \langle \widehat\kH', \widehat\kH''\rangle$. Hence, the
image of the relation (\ref{E:actofAut1}) is the following  equality in the Drinfeld double:
\begin{equation}\label{E:actofAut2}
[\widehat\kF]^- * [\widehat\kG]^+ =
v^{\langle \widehat\kF, \widehat\kG \rangle} \sum\limits_{\widehat{\kH}', \widehat{\kH}''}
v^{-\langle \widehat{\kH}', \widehat{\kH}''\rangle} \; F^{\widehat\kH' \oplus
\widehat\kH''[1]}_{\widehat{\kF}[1], \widehat\kG} \; [\kH']^+ * [\kH'']^-.
\end{equation}
It remains to note that $\mu(\widehat\kF) < \mu(\widehat\kG)$ and the equality
(\ref{E:actofAut2}) is nothing but the relation of the Drinfeld double
$R\bigl([\widehat\kF]^-, [\widehat\kG]^+\bigr)$.

\vspace{2mm}
\noindent
\underline{Case 3}. In a similar way, if  $\Phi(\kF) \cong \widehat\kF[2i]$ and
$\Phi(\kG) \cong \widehat\kG[2i-1]$ for some $i \in \mathbb{Z}$ then
the relation $P\bigl([\kF], [\kG]\bigr)$ is mapped to the relation
$R\bigl([\widehat\kF]^+, [\widehat\kG]^-\bigr)$.

\vspace{3mm}
\noindent
Now we check  the preservation of the relations of the Drinfeld double
$R\bigl([\widehat\kF]^\pm, [\widehat\kG]^\mp\bigr)$
for all semi-stable objects $\kF$ and $\kG$.

\vspace{2mm}
\noindent
\underline{Case 1}. First assume  $\mu(\kF) = \mu(\kG) = \mu$. Recall that any  auto-equivalence
$\Phi \in \Aut\bigl(D^b(Coh(\E))\bigr)$ induces an equivalence of abelian categories
$\cC_\mu \cong \cC_\nu$ for an appropriate slope $\nu$. In particular, we obtain:
$$\Phi\bigl(R\bigl([\kF]^\pm, [\kG]^\mp\bigr)\bigr)  = R\bigl([\widehat\kF]^\pm, [\widehat\kG]^\mp\bigr)
$$
where $\Phi(\kF) \cong  \widehat{\kF}[i]$ and $\Phi(\kG) \cong  \widehat{\kG}[i]$
for an appropriate  $i \in \mathbb{Z}$.

\vspace{2mm}
\noindent
\underline{Case 2}. Assume assume  $\mu(\kF) < \mu(\kG)$. Then there exists an auto-equivalence
$\Psi$ such that both complexes $\widehat\kF := \Psi\bigl(\kF[1]\bigr)$ and
$\widehat\kG := \Psi(\kG)$ belong to the heart of the standard t-structure $Coh(\E)$. Since we have already
shown that $R\bigl([\kF]^-, [\kG]^+\bigr) = \Psi\bigl(P([\widehat\kF]^+, [\widehat\kG]^+)\bigr)$, we have:
$$
\Phi\bigl(R\bigl([\kF]^-, [\kG]^+\bigr)\bigr)
 = \Phi\circ \Psi^{-1} \bigl(P([\widehat\kF]^+, [\widehat\kG]^+)\bigr).
 $$
Hence, this is again a relation either of the type
$P\bigl([\widehat\kF]^\pm, [\widehat\kG]^\pm\bigr)$ or of the type
$R\bigl([\widehat\kF]^\mp, [\widehat\kG]^\pm\bigr)$.

\vspace{2mm}
 \noindent
\underline{Case 3}. The remaining case   $\mu(\kF) > \mu(\kG)$ is similar to the former one and is left
to the reader. Theorem \ref{P:braidact} is proven.
\qed

\vspace{.1in}

\begin{cor} The group $\widehat{SL}(2,\Z)$ acts by algebra automorphisms on $\mathbf{D}{\mathbf{H}}_{\E}$.\end{cor}

\vspace{.1in}

\addtocounter{theo}{1}
\paragraph{\textbf{Remark \thetheo}} Theorem  \ref{P:braidact} is close in spirit to \cite{Kap2} (see also \cite{X} and \cite{P?}). Recently, it has been  generalized by Cramer, who proved that any
derived auto-equivalence between hereditary abelian categories gives rise to an isomorphism at the level of Hall algebras, see \cite{Cramer}.

\vspace{.1in}

\noindent
It turns out that the Drinfeld double $\mathbf{D}\mathbf{H}_{\E}$ carries one more symmetry~:

\begin{prop}\label{T:Duali}
The duality functor  $D= \mathbf{R}{\mathcal Hom}(-,\kO)$ induces an involutive
anti-isomorphism $[\mathcal{F}] \mapsto \bigl[D(\mathcal{F})\bigr]$ of the algebra $\mathbf{D}\mathbf{H}_{\E}$ satisfying the
relation
$$
D \circ \Phi =  i^* \circ [1] \circ \Phi \circ D,
$$
where $\Phi = T_\kO T_{\kO_{x_o}} T_\kO$ and
$i$ is an  involution of the curve $\E$ preserving $x_0$.
\end{prop}

\noindent
\textit{Proof}. The proof of the fact that $D$ is an anti-homomorphism of
the algebra $\mathbf{D}\mathbf{H}_{\E}$ is completely analogous to the
proof of the Theorem \ref{P:braidact} and is therefore skipped. The equality
relating the dualizing functor and the Fourier-Mukai transform is a
corollary of the Proposition~A.2. (see Appendix A).\qed

\vspace{.1in}

\addtocounter{theo}{1}
\paragraph{\textbf{Remark \thetheo}} \label{R:dual} Since the map $D$ sends vector bundles to vector bundles, it restricts to an antiinvolution of the subalgebra
$\H^{+,\mathsf{vec}}:= \vec{\bigotimes\limits_{-\infty < \mu < \infty}} \mathbf{H}^{+, (\mu)}$.

\vspace{.2in}

\section{The algebra ${\mathbf{U}}_\E$}

\vspace{.2in}

Our main object of study is a subalgebra ${\mathbf{U}}_\E$ of $\mathbf{D}{\mathbf{H}}_{\E}$, generated by certain ``averages'' of semi-stable sheaves. Before defining ${\mathbf{U}}_\E$ and giving some of its first properties,
we state some useful results on the classical Hall algebra, associated to the category  of torsion sheaves supported at
a point (or equivalently to the category of nilpotent representations of the Jordan quiver).

\vspace{.1in}

\paragraph{\textbf{4.1.}} We shall need the usual notions of $\nu$-integers~: if $\nu \neq \pm 1$ we set
$$[s]_{\nu}=\frac{\nu^{s}-\nu^{-s}}{\nu-\nu^{-1}}.$$
We shall  usually only use $[s]:=[s]_v$ where $v^2=\# \kk^{-1}$ is as in Section~2.1.
For a finite field $\ll$ fix $u \in \C$ such that $u^2=(\#\ll)^{-1}$.
Denote by  $\mathcal{N}_{\ll}$
the category of nilpotent representations over $\ll$ of the quiver consisting
of a single vertex and a single loop.
Then  there is exactly one indecomposable object
$I_{(r)}$ of length $r$ for any $r  \in \mathbb{N}$, and for a partition
 $\lambda=(\lambda_1, \ldots, \lambda_s)$ we write
$I_{\lambda}=I_{(\lambda_1)} \oplus \cdots \oplus I_{(\lambda_s)}$. The set
$\{I_{\lambda}\}$,  where $\lambda$ runs among all partitions is a complete collection of
non-isomorphic
objects in $\mathcal{N}_{\ll}$.
The structure of the Hall algebra $\mathbf{H}_{\mathcal{N}_{\ll}}$ of the category
$\mathcal{N}_{\ll}$ is completely described in \cite{Mac},
Chap.~II, (see also \cite{Mac} Chap. III. 3.4).
The following proposition summarizes those properties of $\mathbf{H}_{\mathcal{N}_{\ll}}$
that will be needed later on. Let us denote by
$\Lambda_t$ Macdonald's ring of symmetric functions, defined over the ring $\Q[t^{\pm 1}]$, and by $e_{\lambda}$ (resp. $p_{\lambda}$)
the elementary (resp. power-sum) symmetric functions.

\begin{prop}[\cite{Mac}]\label{P:Mac}  The  assignment $[I_{(1)^r}] \mapsto u^{r(r-1)}e_r$ extends to a bialgebra isomorphism $\Psi_{\ll}:
\mathbf{H}_{\mathcal{N}_{\ll}} \stackrel{\sim}{\to}
(\Lambda_t)_{|t=u^2}$. Set
$F_r= \Psi_{\ll}^{-1}(p_r)$. Then
\begin{enumerate}
\item[i)]$F_r= \sum_{|\lambda|=r} n_u\bigl(l(\lambda)-1\bigr)[I_{\lambda}]$, where
$n_{u}(l)=\prod_{i=1}^l (1-u^{-2i})$,
\item[ii)] $\Delta(F_r)=F_r \otimes 1 + 1 \otimes F_r$,
\item[iii)] $(F_r,F_s)= \delta_{r,s}
\frac{\displaystyle ru^r}{\displaystyle u^{-r}-u^r}$.
\end{enumerate}
\end{prop}
\begin{proof} Statements i), ii) and iii) may be found in \cite{Mac}, III. 7. Ex.2, I.5 Ex. 25 and III.4 (4.11) respectively
\end{proof}

\noindent
In particular, the scalar product $(\;\,,\;)$ on $\mathbf{H}_{\mathcal{N}_{\ll}}$
coincides, up to a renormalization, with the
Hall-Littlewood scalar product.

\vspace{.2in}

\paragraph{\textbf{4.2.}} Let $x$ be a closed point of $\E$. Since the residue field
$\kk_x$ at the point $x$ is of the same characteristic as $\kk$,
there is an equivalence of categories $ \mathcal{N}_{\kk_x} \stackrel{\sim}{\to}
\mathcal{T}or_x$ which provides us with an isomorphism
$\Psi_{\kk_x}: \mathbf{H}_{\mathcal{T}or_x} \stackrel{\sim}{\to} (\Lambda_t)_{|t=v^{2\deg(x)}}$,
where $v^2=\# \kk^{-1}$. For  $r \in \N$ we define
an
element $ T^{(\infty)}_{r,x} \in {\H}_{\E} $ by the equation
$$\frac{T^{(\infty)}_{r,x}}{[r]}=\begin{cases} 0 & \mathrm{if}\; r \not\equiv 0\;
(\mathrm{mod}\;{\deg(x)})\\
\frac{\deg(x)}{r}\Psi_{\kk_x}^{-1}\bigl(p_{\frac{r}{\deg(x)}}\bigr) & \mathrm{if} \;r
\equiv 0\;(\mathrm{mod}\;{\deg(x)})
\end{cases}$$
 and we put $T^{(\infty)}_r=\sum_x T_{r,x}^{(\infty)}$. Note that this sum is finite since there are only finitely many closed points
on $\E$ of a given degree.

\vspace{.05in}

Recall the subalgebras $\H_{\E}^{(\mu)}$ of $\H_{\E}$ defined in Section~2.4. In particular,
$\H_{\E}^{(\infty)}$ is the Hall algebra of the category of torsion sheaves on $\E$.
As the Hall algebra of $\mathcal{N}_{{\kk}_x}$ is commutative for any $x$, it follows that $\H_{\E}^{(\infty)}$ is commutative  and hence for any slope $\mu$ the algebra $\H_{\E}^{(\mu)}$ is commutative as well.

\vspace{.1in}

By definition, $T^{(\infty)}_r \in {\H}_{\E}^{(\infty)}$. For an arbitrary $\mu \in \Q$ we put $T_r^{(\mu)} =
\epsilon_{\mu,\infty}\bigl(T_r^{(\infty)}\bigr)$. As $\epsilon_{\mu_1,\mu_2} \circ \epsilon_{\mu_2,\mu_3} \simeq \epsilon_{\mu_1,\mu_3}$, we
have $\epsilon_{\mu_1,\mu_2}\bigl(T^{(\mu_2)}_{r}\bigr)=T^{(\mu_1)}_r$ for any $\mu_1, \mu_2$.

\vspace{.1in}

\addtocounter{theo}{1}
\paragraph{\textbf{Definition \thetheo}} Let ${\U}^+_{\E} \subset
{\H}_{\E}^+$ be the
$\KK$-subalgebra generated by all elements
$T_r^{(\mu)}$ for $r \geq 1$ and
$\mu \in \mathbb{Q} \cup \{\infty\}$, and let $\mathbf{U}^-_{\E} \subset
{\H}_{\E}^-$ be the (isomorphic) subalgebra
defined in a similar  way. We denote by
${\mathbf{U}}_\E$ the subalgebra of $\mathbf{D}{\mathbf{H}}_{\E}$ generated by
${\mathbf{U}}^+_{\E}$ and ${\U}^-_{\E}$.

\vspace{.1in}

It will be convenient for us to introduce one more type of notation~: if
$\mu=\frac{l}{n}$ with $n \geq 1$ and $l,n$ relatively prime, we put
$T_{(\pm rn,\pm rl)}=(T^{(\mu)}_r)^\pm \in \mathbf{U}^{\pm}_{\E}$.
 Similarly,
we put
$T_{(0,\pm r)}=(T^{(\infty)}_r)^\pm$ and $T_{(0,0)} = 1$.
We also set $\ZZ=\Z^2$, so that
$$\ZZ^{\pm}=\big\{ (q,p) \in \Z^2\;|\; \pm q >0\;\text{or}\;q=0, \pm p \geq 0\big\}, \qquad \ZZ^* =\Z^2 \backslash \{(0,0\}$$
and $\ZZ=\ZZ^+ \cup \ZZ^-$.
Thus, by definition, ${\U}^\pm_{\E}$ is the subalgebra of
$\mathbf{D}{\mathbf{H}}_{\E}$ generated by
$T_{(q,p)}$ for $(q,p) \in (\ZZ)^{\pm}$.
Note also that by construction the $\widehat{SL}(2,\Z)$-action on
$\mathbf{D}{\mathbf{H}}_{\E}$ preserves ${\mathbf{U}}_\E$. In fact, since $i^* \bigl(T_{(q,p)}\bigr) = T_{(q,p)}$ for any involution $i$ of $\E$, this action factors through $SL(2,\Z)$.

Finally, it will be necessary to consider a new system of generators for
$\mathbf{U}^{+}_{\E}$. Namely, for $\a \in \ZZ^+$ we put
$$\mathbf{1}_{\a}^{\mathsf{ss}}=\sum_{{\overline{\mathcal{H}}=\a};{\mathcal{H}
\in {\cC}_{\mu(\a)}}} [\mathcal{H}] \, \in \,  \mathbf{H}^+_{\E}[\a].
$$
This sum is finite. If $\a=(q,p)$ with $p,q$ relatively prime then
(see e.g.~\cite{S0}, Section~6.3.) we have
\begin{equation}\label{E:oopo}
1+\sum_{r \geq 1} \mathbf{1}^{\mathsf{ss}}_{r\a}s^r = \exp\bigg(\sum_{r \geq 1} \frac{T_{r\a}}{[r]}s^r\bigg).
\end{equation}
 In particular, $\mathbf{1}_{\a}^{\mathsf{ss}} \in \mathbf{U}^+_{\E}$ and the set $\bigl\{\mathbf{1}_{\a}^{\mathsf{ss}}\;|\; \a \in \ZZ^+\bigr\}$ indeed generates
$\mathbf{U}^+_{\E}$.

\vspace{.2in}

\paragraph{\textbf{4.3.}} Let us now introduce completions of $\mathbf{U}^+_{\E}$ and $\mathbf{U}^+_{\E} \otimes \mathbf{U}^+_{\E}$. Put
$\mathbf{U}^{\not\ge n}_{\E}[\alpha] := \mathbf{U}^+_\E[\alpha] \cap \H_\E^{\not\ge n}[\a]$ and
$\mathbf{U}^{\ge n}_{\E}[\alpha] := \mathbf{U}^{+}_{\E}/\mathbf{U}^{\not\ge n}_{\E}[\alpha]$.
We  can define
\begin{equation}
\widehat{\mathbf{U}}^+_{\E}[\alpha] := \lim_{\underset{n}{\longleftarrow}}
\mathbf{U}^{\ge n}_{\E}[\alpha],
\end{equation}
 then clearly
$\widehat{\mathbf{U}}^+_{\E}[\a] \subseteq
\widehat\H^+_\E[\alpha].
$
In the same way, we denote
\begin{equation}
\mathbf{U}^+_\E[\alpha] \widehat\otimes \mathbf{U}^+_\E[\beta] :=
\lim_{\underset{n,m}{\longleftarrow}}  \mathbf{U}^{\ge n}_{\E}[\alpha]
\otimes \mathbf{U}^{\ge m}_{\E}[\beta]
\subseteq \H^+_\E[\alpha] \widehat\otimes \H^+_\E[\beta].
\end{equation}
By definition, an element $a \in \widehat{\H}^+_{\E}[\a]$ belongs to $\widehat{\U}_{\E}^+[\a]$ if and only if $\jet_n(a) \in \U^{\geq n}_{\E}[\a]$ for all $n$. Similarly,
$a \in \H^+_\E[\a] \widehat{\otimes} \H^+_{\E}[\b]$ belongs to $\mathbf{U}^+_\E[\alpha] \widehat\otimes \mathbf{U}^+_\E[\beta]$ if and only if $\jet_{m,n}(a) \in \U^{\geq m}_{\E}[\a] \otimes \U^{\geq n}_{\E}[\b]$ for all $m,n$.

\vspace{.1in}
\noindent
Next, we set
$$
\widehat{\mathbf{U}}_\E^+ := \bigoplus\limits_{\alpha \in \ZZ^+}
\widehat{\mathbf{U}}_\E^+[\alpha]
\mbox{\textit  \, \, and \, \,}
\mathbf{U}_\E^+ \widehat\otimes \mathbf{U}_\E^+ =
\bigoplus\limits_{\alpha \in \ZZ^+}
\big(\prod\limits_{\beta + \gamma = \alpha}
\mathbf{U}_\E^+[\beta] \widehat\otimes \mathbf{U}_\E^+[\gamma]\big).
$$
The aim of this section is to prove the following result~:

\begin{prop}\label{P:bigebre} $\widehat{\U}_{\E}^+$ is a topological sub-bialgebra of $\widehat{\H}^+_{\E}$. That is,   $\widehat{\U}_{\E}^+$ is stable under the product, and we have 
${\Delta}_{\alpha, \beta}\bigl(\widehat{\mathbf{U}}_\E^+[\alpha + \beta]\bigr) \subset
\mathbf{U}_\E^+[\alpha] \widehat\otimes \mathbf{U}_\E^+[\beta]$.
\end{prop}
\noindent
\textit{Proof.} We first show that  $\widehat{\U}_{\E}^+$ is stable under multiplication. Let $a  \in \widehat{\mathbf{U}}_\E^+[\alpha]$ and
$b  \in \widehat{\mathbf{U}}_\E^+[\beta]$, and fix $u_n \in  {\mathbf{U}}_\E^+[\alpha]$, $v_n \in {\mathbf{U}}_\E^+[\b]$ so that
$\jet_n(a) = \jet_n(u_n)$,   and
$\jet_n(b) = \jet_n(v_n)$ for all $n$.
Then by the continuity of the product (Lemma \ref{L:cont})
for all $m$ we can find $n$ such that
$\jet_m(ab)= \jet_m(u_n v_n)$, hence $ab \in \widehat{\mathbf{U}}_\E^+[\alpha + \beta]$.  The same proof shows that $\U^+_{\E} \widehat{\otimes} \U^+_{\E}$ is a subalgebra of $\H^+_{\E} \widehat{\otimes}\H^+_{\E}$.

To prove the stability of $\widehat{\U}_{\E}^+$ under the coproduct, it is enough to show that $\Delta(\mathbf{1}_{\a}^{\mathsf{ss}}) \in \U^+_{\E} \widehat{\otimes} \U^+_{\E}$ for any $\a$.
For this, we introduce another set of generators, this time for $\widehat{\U}_{\E}^+$.
Namely, we denote
$$
\mathbf{1}_\alpha := \sum\limits_{\kF; \, \overline\kF = \alpha} [\kF] \in
\widehat{\H}^+_\E[\alpha]
$$
From the existence and splitting of the Harder-Narasimhan filtrations  we deduce
the following equality  in $\widehat{\H}_\E[\alpha]$
\begin{equation}\label{E:flip}
\mathbf{1}_{\a}=\mathbf{1}_{\a}^{\mathsf{ss}} + \sum_{t > 1}
\sum_{\substack{\a_1 + \cdots + \a_t=\a \\ \mu(\a_1)< \cdots < \mu(\a_t)}}
v^{\sum_{i<j}\langle \a_i,\a_j\rangle}\mathbf{1}_{\a_1}^{\mathsf{ss}} \cdots
\mathbf{1}_{\a_t}^{\mathsf{ss}},
\end{equation}
from which we conclude  that $\jet_n(\mathbf{1}_{\a}) \in \U^{\geq n}_{\E}[\a]$ for any $n$, and thus
$\mathbf{1}_{\a}$ belongs to $\widehat{\U}_{\E}^+$. Next, we use the following well-known property of Hall algebras (see e.g. \cite[Lemma 1.7]{SLectures})~:  for any $\a,\beta \in \ZZ^+$,
\begin{equation}\label{E:coprodun}
{\Delta}_{\a,\b}(\mathbf{1}_{\a+\b})= v^{\langle \a,\b \rangle}
\mathbf{1}_{\a}
\otimes \mathbf{1}_{\b}.
\end{equation}
It follows that for any polynomial $u=u(\mathbf{1}_{\a_1}, \ldots, \mathbf{1}_{\a_r})$ we have $\Delta(u) \subset \mathbf{U}^+_{\E}\widehat{\otimes}\mathbf{U}^+_{\E}$. The inclusion $\Delta(\mathbf{1}_{\a}^{\mathsf{ss}}) \in \U^+_{\E} \widehat{\otimes} \U^+_{\E}$ is thus a consequence of the continuity of the map $\Delta$ together with the next Lemma~:

\begin{lem} For any $\alpha \in \ZZ^+$ and any $n\in \mathbb{Z}$ there exists an integer $m(n)$,
a polynomial $u_n \in \KK\bigl[t_1,t_2,\dots,t_{m(n)}\bigr]$ and classes
$\a_1,\a_2,\dots,\a_{m(n)} \in \ZZ^+$ satisfying  $\sum\limits_{i=1}^{m(n)} \a_i = \a$ such that
$
\jet_n(\mathbf{1}^{\mathsf{ss}}_\a) = \jet_n\bigl(u_n(\mathbf{1}_{\a_1},\mathbf{1}_{\a_2},\dots, \mathbf{1}_{\a_{m(n)}})\bigr).
$
\end{lem}
\noindent
\textit{Proof}. We prove this lemma by induction on $\rank(\a)$. The case $\rank(\a) = 0$ is clear since
$\mathbf{1}^{\mathsf{ss}}_\a = \mathbf{1}_\a$. Assume that $\a = (r,d)$ and that the assertion is proven for all classes $\b$ such that
$\rank(\b) < r$.
From the formula (\ref{E:flip})  we get the following expression in $\widehat{\H}_\E[\a]$:
\begin{equation*}
\mathbf{1}^{\mathsf{ss}}_{\a}=\mathbf{1}_{\a} -
 \sum_{\substack{\rank(\a) = \rank(\gamma)\\ n \le  \mu(\gamma)  < \mu(\a)}}\hspace{-.15in}
v^{\langle \gamma,\a-\gamma\rangle}\mathbf{1}^{\mathsf{ss}}_{\gamma}\mathbf{1}^{\mathsf{ss}}_{\a-\gamma} -  \sum_{t>1} \hspace{-.1in}
\sum_{\substack{\b_1+\cdots +
\b_t=\a\\ \rank(\b_i) < \rank(\a) \\ n \le  \mu(\b_1)< \cdots < \mu(\b_t)}}\hspace{-.15in}v^{\sum_{i<j}\langle \b_i,\b_j\rangle}\mathbf{1}^{\mathsf{ss}}_{\b_1} \cdots \mathbf{1}^{\mathsf{ss}}_{\b_t} + r_n,
\end{equation*}
where $\jet_n(r_n) = 0$. Note that  both  sums in the right hand side of the equality
are finite. In particular, there
are  finitely many classes $\gamma = (r,d')$ such that $d>d'$ and $\mu(\gamma) \ge  n$. Applying
the above formula  an
appropriate number of times to the element $\mathbf{1}^{\mathsf{ss}}_\gamma$
we obtain
\begin{equation}\label{E:recurs2}
\mathbf{1}^{\mathsf{ss}}_{\a}=\mathbf{1}_{\a} + \sum_{\substack{i=1 \\ \rank(\gamma_i) = \rank(\a)}}^k
\mathbf{1}_{\gamma_i} p_i +
\sum_{\substack{j=1 \\ \beta_1 + \cdots + \beta_{n(j)} = \a \\ \rank(\beta_i) < \rank(\a)}}^l
q_j \, \mathbf{1}^{\mathsf{ss}}_{\beta_1} \mathbf{1}^{\mathsf{ss}}_{\beta_2} \dots \mathbf{1}^{\mathsf{ss}}_{\beta_{n(j)}} +  r'_n,
\end{equation}
where  $\jet_n(r'_n) = 0$,
$p_i$ are polynomials in elements of type  $\mathbf{1}_{(0,l)}$ and  $q_j$ are scalars.
Now by the continuity of the product (Lemma \ref{L:cont}) there exists an integer $N$ such that for all classes
$\beta_1, \beta_2,\dots,\beta_{n(j)}$ occurring in the decomposition (\ref{E:recurs2})
of $\mathbf{1}^{\mathsf{ss}}_\a$  and for any $x_1 \in \widehat{\H}_\E[\beta_1], x_2 \in \widehat{\H}_\E[\beta_2], \dots, x_{n(j)} \in
\widehat{\H}_\E[\beta_{n_j}]$ we have
$$
\jet_n(x_1 x_2 \dots x_{n(j)}) = \jet_n\bigl(\jet_N(x_1)\jet_N(x_2) \dots \jet_N(x_{n(j)})\bigr).
$$
Approximating the elements $\mathbf{1}^{\mathsf{ss}}_{\beta_1}, \mathbf{1}^{\mathsf{ss}}_{\beta_2},
\dots,
\mathbf{1}^{\mathsf{ss}}_{\beta_{n(j)}}$
up to the order $N$
by polynomials $u_{\beta_1}, u_{\beta_2},\dots, u_{\beta_{n(j)}}$ in classes
$\mathbf{1}_\gamma$, we obtain the desired polynomial $u_n$, approximating
$\mathbf{1}^{\mathsf{ss}}_\alpha$ up to the order $n$.
This concludes the proof of the Lemma and of Proposition~\ref{P:bigebre}. \qed

\vspace{.2in}

\paragraph{\textbf{4.4.}}  We now come to the main result of this Section. In Section 2.4. we described a PBW-type basis of ${\mathbf{H}}^+_{\E}$, which by Proposition \ref{P:DDoub} extends to a PBW basis
of ${\mathbf{H}}_{\E}$.

We give a similar
construction for $\mathbf{U}^+_{\E}$.
For $\mu \in \mathbb{Q} \cup \{\infty\}$ let us denote by
$\mathbf{U}_{\E}^{\pm, (\mu)} \subset \mathbf{U}_{\E}^{\pm}$ the subalgebra generated by $\bigl\{(T^{(\mu)}_r)^{\pm}\;|\;r \geq 1\bigr\}$. We also
let $\vec{\bigotimes\limits_{\mu}}$ stand for the restricted ordered tensor product
(see Section~2.4.).

\begin{theo}\label{L:bial} The multiplication map induces isomorphisms of
$\KK$-vector spaces 
\begin{equation}\label{E:pbwU}
 \vec{\bigotimes\limits_{\mu}} \mathbf{U}^{\pm,(\mu)}_{\E}
\stackrel{\sim}{\to}
\mathbf{U}^{\pm}_{\E}, \qquad
\vec{\bigotimes\limits_{\mu}} \mathbf{U}^{+,(\mu)}_{\E} \otimes
 \vec{\bigotimes\limits_{\mu}}
\mathbf{U}^{-,(\mu)}_{\E}\stackrel{\sim}{\to} {\mathbf{U}}_{\E}.
\end{equation}
Moreover, $\mathbf{U}^{\pm}_{\E}$ is a
topological  bialgebra:
${\Delta}_{\a,\b}(\mathbf{U}^{\pm}_{\E}[\a+\b]) \subset
\mathbf{U}^{\pm}_{\E}[\a]  \otimes {\mathbf{U}}^{\pm}_{\E}[\b]$.
\end{theo}

\noindent
\textit{Proof.}  In order to prove the above result, it will be convenient to consider  $\overline{\mathbf{U}}_\E^+[\alpha] :=
\widehat{\mathbf{U}}^+_\E[\alpha] \cap \H^+_\E[\alpha]$, where the intersection is
taken in $\widehat{\H}^+_\E[\alpha]$. Of course, $\mathbf{U}^+_\E[\alpha] \subset \overline{\mathbf{U}}^+_\E[\alpha]$ and as it will turn out in the end  $\overline{\mathbf{U}}^+_\E[\alpha] =  \mathbf{U}^+_\E[\alpha]$,
but \textit{a priori}
$\overline{\mathbf{U}}_\E^+[\alpha]$
might be bigger. Observe however that $\overline{\mathbf{U}}^+_\E[\alpha] =  \mathbf{U}^+_\E[\alpha]$ for any  class $\a = (0,d)$, $d \in \mathbb{Z}_{>0}$.  It is easy to see that
$\overline{\U}^+_{\E}=\bigoplus\limits_{\a
\in \ZZ^+}\overline{\mathbf{U}}_\E^+[\alpha]$
is a subalgebra of $\H^+_{\E}$. In addition, it is also a sub-bialgebra~:

\begin{lem}\label{L:huh6} For any $\a,\b$ we have 
$\Delta_{\a,\b}\bigl(\overline{\U}_{\E}^+[\a+\b]\bigr) \subset
\overline{\U}_{\E}^+[\a]\otimes \overline{\U}_{\E}^+[\b]$.
\end{lem}
\noindent
\textit{Proof.} By Proposition~\ref{P:bigebre} It is enough to show that
\begin{equation}\label{E:huh}
\bigl(\mathbf{U}^+_\E[\alpha] \widehat\otimes \mathbf{U}^+_\E[\beta]\bigr) \cap
\bigl(\H^+_\E[\alpha] \otimes \H^+_\E[\beta]\bigr) = \overline{\mathbf U}^+_\E[\alpha] \otimes
\overline{\mathbf U}^+_\E[\beta],
\end{equation}
where the intersection is
taken in $\H^+_\E[\alpha] \widehat\otimes \H^+_\E[\beta]$.  Let $V_{\a,\b}$ stand for the left hand side
 of (\ref{E:huh}). The inclusion $\overline{\mathbf U}^+_\E[\alpha] \otimes
\overline{\mathbf U}^+_\E[\beta] \subset V_{\a,\b} $ is obvious
since $\overline{\mathbf U}^+_\E[\alpha] \otimes
\overline{\mathbf U}^+_\E[\beta] \subset \mathbf{U}^+_\E[\alpha] \widehat\otimes \mathbf{U}^+_\E[\beta]$ and
$\overline{\mathbf{U}}^+_\E[\alpha] \otimes \overline{\mathbf{U}}^+_\E[\beta]
\subset \H^+_\E[\alpha] \otimes \H^+_\E[\beta]$.  For any sheaf $\mathcal{F}$ of class $\gamma$ let
$\pr_{\mathcal{F}}: \widehat{\H}^+_{\E}[\gamma] \to \C$ be the linear form picking the coefficient of $[\mathcal{F}]$. To prove the reverse inclusion, it is enough to show that for any $\mathcal{F}$ of class $\a$ and any $\mathcal{G}$ of class $\b$ we have
\begin{equation}\label{E:huh2}
(\pr_{\mathcal{F}} \otimes 1)(V_{\a,\b}) \subset \overline{\U}^+_{\E}[\b] \qquad
\mbox{and} \qquad
(1\otimes \pr_{\mathcal{G}})(V_{\a,\b}) \subset \overline{\U}^+_{\E}[\a].
\end{equation}
Indeed, if (\ref{E:huh2}) holds then
$$V_{\a,\b} \subset \bigl(\overline{\U}^+_{\E}[\a] \otimes \H^+_{\E}[\b]\bigr) \cap
\bigl(\H^+_{\E}[\a] \otimes \overline{\U}^+_{\E}[\b]\bigr) =\overline{\U}^+_{\E}[\a] \otimes \overline{\U}^+_{\E}[\b].$$
Finally, we prove (\ref{E:huh2}). Let $v \in V_{\a,\b}$, and let $\mathcal{F},\mathcal{G}$ be sheaves of class $\a$ and $\b$ respectively. Choose $m \in \Z$ such that $\mathcal{F}, \mathcal{G} \in
Coh_{\geq m}$. As $v \in \mathbf{U}^+_\E[\alpha] \widehat\otimes \mathbf{U}^+_\E[\beta]$ for any $m'<m$ we have $v \in \U^+_{\E}[\a] \otimes \U^+_{\E}[\b] + 
\bigl(\H^{\not\geq m'}_{\E}[\a] \widehat{\otimes} \H^+_{\E}[\b] + \H^+_{\E}[\a]\widehat{\otimes} \H^{\not\geq m'}_{\E}[\b]\bigr)$ from which we deduce that $(\pr_{\mathcal{F}} \otimes 1)(v) \in \U^+_{\E}[\b] + \widehat{\H}^{\not\geq m'}_{\E}[\b]$ and $(1 \otimes \pr_{\mathcal{G}})(v) \in \U^+_{\E}[\a] + \widehat{\H}^{\not\geq m'}_{\E}[\a]$. Equation (\ref{E:huh2}) follows and the Lemma is proved. \qed

\vspace{.1in}
\noindent
Let $\overline{\U}_{\E} \subset \mathbf{D}\H_{\E}$ be the subalgebra generated by two copies
$\overline{\U}^{\pm}_{\E} \subset \H^{\pm}_{\E}$ of $\overline{\U}^+_{\E}$.

\begin{cor}\label{C:huh3} The algebra $\overline{\U}_{\E}$ is isomorphic to the Drinfeld double
of $\overline{\U}^+_{\E}$, and the multiplication map induces an isomorphism $\overline{\U}^+_{\E}
\otimes \overline{\U}^-_{\E}   \simeq \overline{\U}_{\E}$.
\end{cor}
\noindent
\textit{Proof.} Since $\overline{\mathbf{U}}^+_\E$ is a topological bialgebra,
by the same proof as for Proposition~\ref{P:DDoub} we see that
$\mathbf{D}\overline{\mathbf{U}}^+_\E \cong \overline{\mathbf{U}}_\E^+
\otimes \overline{\mathbf{U}}_\E^-$ and
$\mathbf{D}\overline{\mathbf{U}}^+_\E$ is isomorphic to the
subalgebra of $\mathbf{D}\H_\E$ generated by $\overline{\mathbf{U}}^+_\E$ and
$\overline{\mathbf{U}}^-_\E$.\qed

\vspace{.1in}

Recall that we defined for any $\nu \in \Q \cup \{\infty\}$ a subalgebra $\H_{\E}^{+,(\nu)}$ and that we set $\U_{\E}^{+,(\nu)}=\U^+_{\E} \cap \H_{\E}^{+,(\nu)}$. By Lemma~\ref{L:PBW}, we have a tensor product decomposition $m:  \vec{\bigotimes\limits_{\nu}} \H^{+,(\nu)}_{\E} \stackrel{\sim}{\to} \H^+_{\E}$. Let $\H_{\E}^{+,\mathsf{vec}}=m\big(\vec{\bigotimes\limits_{\nu<\infty}} \H^{+,(\nu)}_{\E} \big)$ be the subspace spanned by the classes of the vector bundles, so that the multiplication map
$m: \H_{\E}^{+,\mathsf{vec}} \otimes \H^{+,(\infty)}_{\E} \to \H^+_{\E}$ is an isomorphism.
The following property of $\overline{\U}^+_{\E}$ will be crucial for our purposes.

\begin{lem}\label{L:huh4} We have $\overline{\U}^+_{\E} \subset \H^{+, \mathsf{vec}}_{\E}
\otimes \U^{+,(\infty)}_\E$.
\end{lem}

\noindent
\textit{Proof.} Consider an element 
 $u \in \overline{\U}^+_{\E}[\a]$. Viewing it as  an element of
 $m\bigl(\H^{+,\mathsf{vec}}_{\E} \otimes \H_{\E}^{+,(\infty)}\bigr)$ we can expand it 
 as a finite sum $u=\sum_l u_l$ with $u_l=\sum_i u'_{l,i} \cdot u''_{l,i}$, where
$u'_{l,i} \in \H^{+, \mathsf{vec}}_{\E}$ and $u''_{l,i} \in \H^{+,(\infty)}_{\E}\bigl[(0,l)\bigr]$
for all  $i,l$. Let
$\pi: \H^+_{\E} \to \H^{+,\mathsf{vec}}_{\E}$ denote the  projection of the Hall algebra on its subspace.
Observe that as any coherent sheaf $\mathcal{F}$ has a unique maximal torsion subsheaf, we get
\begin{equation}\label{E:huh5}
(\pi \otimes 1)\Delta_{\a-(0,l),(0,l)}(u)=v^{\bigl\langle \a-(0,l),(0,l)\bigr\rangle} \sum_i u'_{l,i} \otimes u''_{l,i}.
\end{equation}
On the other hand, by Lemma~\ref{L:huh6} we have
$$\Delta_{\a-(0,l),(0,l)}(u) \in \overline{\U}^+_{\E}\bigl[\a-(0,l)\bigr] \otimes \overline{\U}^+_{\E}\bigl[(0,l)\bigr]=\overline{\U}^+_{\E}\bigl[\a-(0,l)\bigr] \otimes {\U}^+_{\E}\bigl[(0,l)\bigr].$$
 But then from (\ref{E:huh5}) we obtain
 $\sum_i u'_{l,i} \otimes u''_{l,i} \in \H^{+, \mathsf{vec}}_{\E} \otimes \U_{\E}^+\bigl[(0,l)\bigr]$for all $l$, and hence $u \in \H^{+, \mathsf{vec}}_{\E} \otimes \U^{+, (\infty)}_{\E}$ as wanted.  \qed

\vspace{.1in}

After these preliminaries we are now
ready to prove that the multiplication map induces an isomorphism
$\vec{\bigotimes\limits_{\nu}} \U^{+,(\nu)}_{\E} \simeq \U^+_{\E}$. For any $\nu \in \Q \cup \{\infty\}$ there exists $\gamma \in SL(2,\Z)$ such that $\gamma(\nu)=\infty$. Recall that
the group $\widehat{SL}(2,\Z)$ acts on $\mathbf{D}\H_{\E}$ and preserves $\U_{\E}$. Moreover, this action is compatible with the decomposition $\mathbf{D}\H_{\E} \simeq
\vec{\bigotimes\limits_{\mu}} \H^{+,(\mu)}_{\E}
\otimes \vec{\bigotimes\limits_{\mu}} \H^{-,(\mu)}_{\E}$ (i.e. it permutes the subalgebras $\H^{\pm,(\mu)}_{\E}$). Hence, using Corollary~\ref{C:huh3} and Lemma~\ref{L:huh4} we obtain the chain of inclusions
\begin{equation}\label{E:huh7}
\gamma(\U^+_{\E}) \subset
\overline{\U}^+_{\E} \otimes \overline{\U}^-_{\E}  \subset   \big(\vec{\bigotimes\limits_{\mu <\infty}} \H^{+,(\mu)}_{\E} \otimes
\U_{\E}^{+,(\infty)}\big) \otimes  \big(\vec{\bigotimes\limits_{\mu <\infty}} \H^{-,(\mu)}_{\E} \otimes
\U_{\E}^{-,(\infty)}\big).
\end{equation}
But then, applying $\gamma^{-1}$ to (\ref{E:huh7}) and using the equality
$\gamma^{-1}\bigl(\U^{+,(\infty)}_{\E}\bigr) =\U^{+,(\nu)}_{\E}$ we see that
$\U^+_{\E} \subset
\vec{\bigotimes\limits_{\mu <\nu}} \H^{+,(\mu)}_{\E} \otimes \U^{+,(\nu)}_{\E}
\otimes \vec{\bigotimes\limits_{\mu>\nu}} \H^{+,(\mu)}_{\E}$. As this is true for all $\nu$, we get
$\U^+_{\E} \subset \vec{\bigotimes\limits_{\mu}}  \U^{+,(\mu)}_{\E}$ and finally
$\U^+_{\E} =\vec{\bigotimes\limits_{\mu}} \U^{+,(\mu)}_{\E}$. Of course, this also proves the equality $\U^-_{\E} =\vec{\bigotimes\limits_{\mu}} \U^{-,(\mu)}_{\E}$. The second statement in
Theorem \ref{L:bial} is now a consequence of Corollary~\ref{C:huh3} and the next result~:

\begin{lem}\label{L:Uisclos}
The two algebras $\mathbf{U}^+_\E$ and $\overline{\mathbf{U}}_\E^+$ coincide.
\end{lem}

\noindent
\textit{Proof}. Recall that the condition $u \in \overline{\mathbf{U}}_\E^+[\alpha]$ means that
$u \in \mathbf{H}_\E^+[\alpha]$ and for all $n$ there exists $u_n \in
\mathbf{U}^+_\E[\alpha]$ such that $\jet_n(u) = \jet_n(u_n)$. But note that for
$n \gg 0$  we have  $u = \jet_n(u)$ and as $\U_\E^+ \simeq \vec{\bigotimes\limits_{\nu}} \U^{+,(\nu)}_{\E}$, we get $\jet_n(u_n) \in \mathbf{U}^+_\E[\alpha]$. Therefore,  $\mathbf{U}_\E^+
= \overline{\mathbf{U}}_\E^+$ and as a corollary,
$$\overline{\mathbf{U}}_\E = \mathbf{D}\overline{\mathbf{U}}^+_\E =
\overline{\mathbf{U}}_\E^+ \otimes \overline{\mathbf{U}}_\E^- =
\mathbf{U}^+_\E \otimes \mathbf{U}^-_\E =
\mathbf{D}\mathbf{U}^+_\E = \mathbf{U}_\E.
$$
This concludes the proof of Theorem~\ref{L:bial}.\qed

\vspace{.2in}

\paragraph{\textbf{4.5.}} We finish this section with several important computations regarding $\U^+_{\E}$. They  will be used in a crucial way in the next section. Let us set, for $i \geq 1$
$$c_i(\E)= \#\E(\mathbb{F}_{q^i}) \, v^i [i]/i.$$

\begin{lem}\label{L:scalar} For any $\x=(q,p) \in \ZZ^+$ we have $$(T_{\x},T_{\x})=\frac{c_r(\E)}{(v^{-1}-v)},$$
where $r = \mathsf{gcd}(q,p) \in \N$.
\end{lem}
\noindent
\textit{Proof.} Using the $SL(2,\Z)$ action, we may restrict ourselves to the case of $\x=(0,r)$ with $r>0$. We have
$$\bigl(T_{(0,r)},T_{(0,r)}\bigr)=\sum_{d | r} \sum_{x: \; \deg(x) = d}
\bigl(T^{(\infty)}_{r,x}, T^{(\infty)}_{r,x}\bigr).$$
If $x$ is of degree $d$, it follows from Proposition~\ref{P:Mac}, iii) that
$$\bigl(T^{(\infty)}_{r,x},
T^{(\infty)}_{r,x}\bigr)=\frac{ \displaystyle d[r]^2}{\displaystyle r(q^r-1)}=\frac{\displaystyle v^r[r]d}{\displaystyle r(v^{-1}-v)}.$$
The statement of the
Lemma follows from the equation
$$\sum_{d | r} \sum_{x: \; \deg(x)=d} d=\#\E(\mathbb{F}_{q^r}).
$$\qed

We now turn to the coproduct. Define elements $\Theta_{\x} \in \U_{\E}$ by equating the
 coefficients of the following generating series:
\begin{equation}
\sum_i \Theta_{i\x_0}s^i = \exp\bigg((v^{-1}-v)\sum_{r \geq 1}T_{r\x_0}s^r\bigg),
\end{equation}
for any $\x_0 \in \ZZ^*$ such that $\deg(\x_0)=1$.

\begin{lem}\label{L:coprodun} For any $p \in \Z$ we have:
$$
\Delta(T_{(1,p)})=T_{(1,p)} \otimes 1 + \sum_{l \geq 0} \Theta_{(0,l)}  \otimes T_{(1,p-l)}.
$$
\end{lem}
\noindent
\textit{Proof.}  Up to a twist  by a line bundle, it is enough to consider $T_{(1,0)}= \mathbf{1}^{\mathsf{ss}}_{(1,0)}$.
By the proof of Proposition \ref{P:bigebre}, we have
$$
\begin{array}{lcl}
\Delta_{(1,-n),(0,n)}
(\mathbf{1}_{(1,0)}) & = & v^n \mathbf{1}_{(1,-n)} \otimes \mathbf{1}_{(0,n)},\\
\Delta_{(0,n),(1,-n)}
(\mathbf{1}_{(1,0)}) & = & v^{-n} \mathbf{1}_{(0,n)} \otimes\mathbf{1}_{(1,-n)},\\
\Delta_{(0,n),(0,m)}(\mathbf{1}_{(0,n+m)}) & = &\mathbf{1}_{(0,n)} \otimes \mathbf{1}_{(0,m)}
\end{array}
$$
and moreover
$$\mathbf{1}^{\mathsf{ss}}_{(1,0)}=\sum_{n \geq 0} v^n \mathbf{1}_{(1,-n)}\chi_n, $$
where
$\chi_n=\sum_{r >0}(-1)^r \sum_{l_1+\cdots + l_r=n} \mathbf{1}_{(0,l_1)} \cdots \mathbf{1}_{(0,l_r)}$.
Denote by
$$
\mathbf{1}(s)=\sum_{l \geq 0} \mathbf{1}_{(0,l)}s^l, \qquad \chi(s)=\sum_{l \geq 0} \chi_l s^l
$$
the generating functions of $\mathbf{1}_{(0,n)}$ and $\{\chi_n\}$.
It is easy to see  that the elements  $\{\chi_n\}$ are completely determined by the relations
$\sum_{i+j=l} \mathbf{1}_{(0,i)} \chi_j=\delta_{l,0}$,  which can  be rewritten in the form  $\mathbf{1}(s) \chi(s)=1$.
 In particular, from the formula  for the coproduct
  we have $\Delta(\mathbf{1}(s))=\mathbf{1}(s) \otimes \mathbf{1}(s)$ from which  we deduce that
$\Delta\bigl(\chi(s)\bigr)=\chi(s) \otimes \chi(s)$, i.e
$\Delta_{(0,n),(0,m)}(\chi_{n+m})=\chi_n \otimes \chi_m$. This implies that
\begin{equation*}
\Delta_{(1,-l),(0,l)}\bigl(\mathbf{1}^{ss}_{(1,0)}\bigr)=\sum_{k \geq 0}
v^{l+k}\mathbf{1}_{(1,-l-k)}\chi_{k} \otimes \bigl(\mathbf{1}_{(0,l)} +
\mathbf{1}_{(0,l-1)}\chi_1 + \cdots + \chi_l\bigr).
\end{equation*}
Using $\sum_{i+j=l} \mathbf{1}_{(0,i)} \chi_j=\delta_{l,0}$, we get
$\Delta_{(1,-l),(0,l)}\bigl(\mathbf{1}^{\mathsf{ss}}_{(1,0)}\bigr)=
\delta_{l,0}\mathbf{1}^{\mathsf{ss}}_{(1,0)} \otimes 1$.
A similar computation shows that
$$\Delta_{(0,l), (1,-l)}(\mathbf{1}^{\mathsf{ss}}_{(1,0)})=\sum_{k \geq 0}\bigl(v^{k-l} \mathbf{1}_{(0,l)} +
v^{2+k-l} \mathbf{1}_{(0,l-1)}\chi_1 + \cdots + v^{k+l} \chi_l\bigr) \otimes \mathbf{1}_{(1,-l-k)}
\chi_{k}.
$$
Hence, setting $\theta_l=\sum_{k=0}^l v^{2k-l} \mathbf{1}_{(0,l-k)}\chi_k$ we obtain
$$\Delta(\mathbf{1}^{\mathsf{ss}}_{(1,0)}) =\mathbf{1}^{\mathsf{ss}}_{(1,0)} \otimes 1 + \sum_{l \geq 0} \theta_l \otimes \mathbf{1}^{\mathsf{ss}}_{(1,-l)}.$$
Finally, we claim that the elements $\theta_l$ can  be characterized through
the relation $\sum\limits_{l \ge 0} \theta_l s^l=
\exp\big((v^{-1}-v)\sum\limits_{r \geq 1} T_{(0,r)}s^r\big)$. To see this, note that by definition
$\sum\limits_{l \geq 0} \mathbf{1}_{(0,l)}s^l =
\exp\bigg(\sum\limits_{r \geq 1} \frac{\displaystyle T_{(0,r)}}{\displaystyle [r]} s^r\bigg)$, hence
$\sum\limits_{l \geq 0} \chi_l s^l = \exp\bigg(-\sum\limits_{r \geq 1} \frac{\displaystyle T_{(0,r)}}{\displaystyle [r]} s^r\bigg)$. But then
\begin{equation*}\begin{split}
\sum_{l \geq 0}  \theta_l s^l=&\mathbf{1}(v^{-1}s)\chi(vs) = \exp\bigg(\sum_{r \geq 1}
v^{-r} \frac{T_{(0,r)}}{[r]} s^r -\sum_{r \geq 1} v^r \frac{T_{(0,r)}}{[r]} s^r\bigg)\\
= & \exp\bigg((v^{-1}-v)\sum_{r \geq 1} T_{(0,r)}s^r\bigg)
\end{split}
\end{equation*}
as desired.
\qed

\vspace{.1in}
\noindent
The final computation which we shall need is the following. Set
$$\mathbf{1}^{\mathsf{vec}}_{\a}=
\sum_{\substack{\mathcal{F}\;{vec.\;bdle} \\ \overline{\mathcal{F}}=\a}} [\mathcal{F}] \in \widehat{\U}^+_{\E}.$$

\begin{lem}\label{L:comp} For any $n \geq 0$ and for any $\a=(r,d) \in \ZZ^+$ we have
\begin{equation}\label{E:compeq1}
\bigl[T_{(0,n)},\mathbf{1}_{\a}\bigr]=c_n(\E) \frac{v^{rn}-v^{-rn}}{v^n-v^{-n}} \mathbf{1}_{\a+(0,n)},
\end{equation}
\begin{equation}\label{E:compeq2}
\bigl[T_{(0,n)},\mathbf{1}^{\mathsf{vec}}_{\a}\bigr]=c_n(\E) \frac{v^{rn}-v^{-rn}}{v^n-v^{-n}} \mathbf{1}^{\mathsf{vec}}_{\a+(0,n)}.
\end{equation}
\end{lem}
\begin{proof} Since $\mathbf{1}_{\a}=\sum_{d \geq 0} v^{\bigl\langle \a, (0,d)\bigr\rangle} \mathbf{1}^{\mathsf{vec}}_{\a-(0,d)} \mathbf{1}_{(0,d)}$ and  $\bigl[T_{(0,n)},
 \mathbf{1}_{(0,d)}\bigr]=0$ for all $n$ and $d$, the equation (\ref{E:compeq1}) is a consequence of (\ref{E:compeq2}). We shall  thus only deal with (\ref{E:compeq2}).

Assume first that $\rank(\a)=1$. Up to twisting by a line bundle, we may assume that $\a=(1,0)$. Note that $\mathbf{1}^{\mathsf{vec}}_{(1,0)}=T_{(1,0)}$.
There exist elements $S_0, \ldots, S_n$ with $S_i$ belonging to the algebra generated by $T_{(0,1)}, \ldots, T_{(0,i)}$ such that $T_{(0,n)}T_{(1,0)}$ is equal to a linear combination
\begin{equation}\label{E:ree1}
T_{(0,n)}T_{(1,0)}=\sum_{i=0}^{n}T_{(1,n-i)}S_i,
\end{equation}
We first compute $S_n$. Let us write $T_{(0,n)}=\sum_{\mathcal{T}} w_{\mathcal{T}} [\mathcal{T}]$ and $S_n=\sum_{\mathcal{T}} u_{\mathcal{T}} [\mathcal{T}]$, for some scalars $w_{\mathcal{T}},u_{\mathcal{T}} \in \KK$. Observe that
a term of the form $[\mathcal{L} \oplus \mathcal{T}]$, for a line bundle $\mathcal{L}$ of degree zero and a torsion sheaf $\mathcal{T}$ of degree $n$, only appears on the right
hand side of (\ref{E:ree1}) in $T_{(1,0)}S_n$, and with a coefficient
equal to $u_{\mathcal{T}}v^{-n}$. On the other hand, the coefficient of $[\mathcal{L} \oplus \mathcal{T}]$ in the left hand side is equal
to $v^nw_{\mathcal{T}} {F}_{\mathcal{T},\mathcal{L}}^{\mathcal{L} \oplus \mathcal{T}}=v^{-n}w_{\mathcal{T}}$. Hence $u_{\mathcal{T}}=w_{\mathcal{T}}$ for all $\mathcal{T}$ and $S_n=T_{(0,n)}$.

Now we show that $S_i=0$ for $i \neq 0,n$. By Proposition~\ref{P:topbial},
$\Delta\bigl([T_{(0,n)},T_{(1,0)}]\bigr) = \bigl[\Delta(T_{(0,n)}),\Delta(T_{(1,0)})\bigr]$. By Proposition~\ref{P:Mac}, ii),
$\Delta\bigl(T_{(0,n)}\bigr)=T_{(0,n)}
\otimes 1 + 1 \otimes T_{(0,n)}$. Let $C = \Delta\bigl([T_{(0,n)},T_{(1,0)}]\bigr)$.  From  Lemma \ref{L:coprodun}  we deduce the formula
\begin{equation}\label{E:ree2}
 C = \bigl[T_{(0,n)}, T_{(1,0)}\bigr] \otimes 1 + 1 \otimes \bigl[T_{(0,n)},T_{(1,0)}\bigr]+\sum_{l \geq 1}
\Theta_{(0,l)} \otimes \bigl[T_{(0,n)},T_{(1,-l)}\bigr].
\end{equation}
 Let $i_0$ be the maximal value of
$i$ distinct from $n$ for which
$S_i\neq 0$.  Note that $\Delta_{(1,n-i_0),(0,i_0)}\bigl(T_{(1,n-i)}S_i\bigr)=0$ if $i<i_0$, while we have
$\Delta_{(1,n-i_0),(0,i_0)}\bigl(T_{(1,n-i_0)}S_{i_0}\bigr)
=v^{n}T_{(1,n-i_0)} \otimes S_{i_0}$. But on the other hand, for any $j >0$, (\ref{E:ree2}) implies that
$\Delta_{(1,n-j),(0,j)}\bigl([T_{(0,n)},T_{(1,0)}]\bigr)=0$.
Hence $i_0=0$, $\bigl[T_{(0,n)},T_{(1,0)}\bigr]=z_0 T_{(1,n)}$ for some $z_0 \in \KK$. In order to determine the value of $z_0$ we compute the scalar product 
$\bigl(T_{(0,n)}T_{(1,0)},\mathbf{1}_{(1,n)}\bigr)$ in two different ways.
By Proposition~\ref{P:Mac} ii),  $T_{(0,n)}$ is orthogonal to the subalgebra generated by $T_{(0,i)}$ for $i <n$.
Hence, using (\ref{E:oopo}) and (\ref{E:flip}), we obtain
$\mathbf{1}_{(1,n)}=T_{(1,n)} + \frac{\displaystyle v^n}{\displaystyle [n]}T_{(1,0)}T_{(0,n)} + u$ where $u \in
\bigl(\KK T_{(1,n)} \oplus
\KK T_{(1,0)}T_{(0,n)}\bigr)^{\perp}$, and using Lemma~\ref{L:scalar} we get
\begin{equation}\label{E:z01}
\begin{split}
\bigl(T_{(0,n)}T_{(1,0)},\mathbf{1}_{(1,n)}\bigr)=&\frac{v^n}{[n]} \bigl(T_{(1,0)}T_{(0,n)},
T_{(1,0)}T_{(0,n)}\bigr) + z_0 \bigl(T_{(1,n)},T_{(1,n)}\bigr)\\
=&\frac{v^n}{[n]}\frac{c_n(\E)c_1(\E)}{(v^{-1}-v)^2} + z_0 \frac{c_1(\E)}{v^{-1}-v}.
\end{split}
\end{equation}
On the other hand, we have $\bigl(T_{(0,n)}T_{(1,0)},\mathbf{1}_{(1,n)}\bigr)
= \bigl(T_{(0,n)} \otimes T_{(1,0)},
\Delta(\mathbf{1}_{(1,n)})\bigr)$ and by
(\ref{E:coprodun}) we have $\Delta(\mathbf{1}_{(1,n)})=
\frac{\displaystyle v^{-n}}{\displaystyle [n]}T_{(0,n)}\otimes T_{(1,0)} + u'$ where $u' \in 
\bigl(\KK T_{(0,n)}
\otimes T_{(1,0)}\bigr)^{\perp}$. It follows that
\begin{equation}\label{E:z02}
\bigl(T_{(0,n)}T_{(1,0)},\mathbf{1}_{(1,n)}\bigr) = \frac{v^{-n}}{[n]} \bigl(T_{(0,n)},T_{(0,n)}\bigr)\bigl(T_{(1,0)},T_{(1,0)}\bigr)=
\frac{v^{-n}}{[n]}\frac{c_1(\E)c_n(\E)}{(v^{-1}-v)^2}.
\end{equation}
Combining (\ref{E:z01}) and (\ref{E:z02}) we finally obtain $z_0=c_n(\E)$ as wanted.

Now let $r=\rank(\a)$ be arbitrary. Repeating the argument above, we have $\bigl[T_{(0,n)},\mathbf{1}^{\mathsf{vec}}_{(r,d)}\bigr] \in \widehat{\U}_{\E
}^{\mathsf{vec}}$. Let $\mathcal{V}$ be a vector bundle of class $\a+(0,n)$ and let $\mathcal{T}$ be a torsion sheaf of degree $n$. The coefficient
of $[\mathcal{V}]$ in $\bigl[[\mathcal{T}],\mathbf{1}^{\mathsf{vec}}_{\a}\bigr]$ is easily seen to be equal to $v^{rn} \big| \text{Hom}^{\mathsf{surj}}(\mathcal{V},\mathcal{T})\big| / a_{\mathcal{T}}$, where $\text{Hom}^{\mathsf{surj}}(\mathcal{V},\mathcal{T})$ stands for the set of surjective maps $\mathcal{V} \tto \mathcal{T}$. By the inclusion-exclusion principle, we have
\begin{equation*}
\begin{split}
\big| \text{Hom}^{\mathsf{surj}}(\mathcal{V},\mathcal{T})\big|&=\big|\text{Hom}(\mathcal{V},
\mathcal{T})\big| -\sum_{\mathcal{T}' \subsetneq \mathcal{T}} \big|\text{Hom}(\mathcal{V},
\mathcal{T}')\big| + \sum_{\mathcal{T}'' \subsetneq \mathcal{T}' \subsetneq \mathcal{T}} \big|\text{Hom}(\mathcal{V},\mathcal{T}'')\big| -\cdots\\
 &=v^{-2rn} -\sum_{\mathcal{T}' \subsetneq \mathcal{T}} v^{-2r \deg(\mathcal{T}')} + \sum_{\mathcal{T}'' \subsetneq \mathcal{T}' \subsetneq \mathcal{T}}
 v^{-2r \deg(\mathcal{T}'')} - \cdots
\end{split}
\end{equation*}
(the sum is finite since there are only finitely many subsheaves of $\mathcal{T}$). The above expression only depends on $\mathcal{T}$ and the rank $r$. Hence 
$\bigl[T_{(0,n)}, \mathbf{1}^{\mathsf{vec}}_{\a}\bigr]=u_r \mathbf{1}^{\mathsf{vec}}_{\a +(0,n)}$ for some $u_r \in \KK$, which remains to be determined. For this, we use the iterated coproduct map $\Delta_{1, \ldots , 1}$. We have, by (\ref{E:coprodun})
\begin{equation}\label{E:xuy1}
\Delta_{1, \ldots , 1}\bigl(\mathbf{1}_{(r,l)}\bigr) = \sum_{l_1 + \cdots + l_n=l} v^{\sum_{i<j} (l_j-l_i)} \mathbf{1}_{(1,l_1)} \otimes \cdots \otimes \mathbf{1}_{(1,l_r)},
\end{equation}
while by Proposition~\ref{P:Mac}, ii),
\begin{equation}\label{E:xuy2}
\begin{split}
\Delta_{1, \ldots , 1}& \bigl([T_{(0,n)},\mathbf{1}_{(r,d)}]\bigr)\\
&=\bigg[\sum_{j=1}^r    1 \otimes \cdots \otimes T_{(0,n)} \otimes \cdots \otimes 1, \Delta_{1, \ldots, 1}(\mathbf{1}_{(r,d)})\bigg]\\
&=u_1\sum_{j=1}^r \sum_{d_1 + \cdots + d_n=d} v^{\sum_{i<j} (d_j-d_i)} \mathbf{1}_{(1,d_1)} \otimes \cdots \otimes \mathbf{1}_{(1, d_j+n)} \otimes \cdots \otimes \mathbf{1}_{(1,d_r)}.
\end{split}
\end{equation}
Comparing (\ref{E:xuy1}) with (\ref{E:xuy2}) and using the case $r=1$ treated above we get
$$u_r=u_1v^{(r+1)n} \sum_{j=1}^r v^{-2jn}=c_n(\E) \frac{v^{rn}-v^{-rn}}{v^n-v^{-n}}$$
as wanted. We are done.
\end{proof}

\vspace{.2in}

\section{The algebra $\UU_{\s,\bs}$}

\vspace{.2in}

\paragraph{\textbf{5.1.}} The aim of this section is to give a presentation for ${\mathbf{U}}_{\E}$ by
generators and relations.
Since it is  convenient to depict elements of ${\mathbf{U}}_{\E}$
graphically, a few notational preparations are in order.
Let $\mathbf{o}$ stand for the origin in $\ZZ$. By a path in $\ZZ$ we shall  mean a finite sequence $\p=(\x_1, \x_2, \ldots, \x_r)$
of non-zero elements of $\ZZ$, which we represent as the piecewise-linear curve in $\ZZ$ joining
$\mathbf{o},\x_1, \x_1+\x_2, \ldots,
\x_1+\cdots + \x_r$. Let $\widehat{\x\y} \in [0, 2\pi[$ denote the angle between the segments
$\mathbf{o}\x$ and $\mathbf{o}\y$. A path $\mathbf{p}=(\x_1, \ldots, \x_r)$ will be called \textit{convex}
if $0 \leq \widehat{\x_1\x_2} \leq \widehat{\x_1 \x_3} \leq \cdots \leq
\widehat{\x_{1}\x_r}<2\pi$.
 Put $L_0=\N (0,-1)$ and let $\mathbf{Conv}'$ be the
collection of all convex paths $\mathbf{p}=(\x_1, \ldots, \x_r)$ satisfying
$\widehat{\x_1 L_0} \geq \widehat{\x_2  L_0} \geq \cdots \geq \widehat{\x_r L_0} \geq 0$.
Two convex paths $\p=(\x_1, \ldots, \x_r)$ and $\mathbf{q}=(\y_1, \ldots, \y_s)$ in $\mathbf{Conv}'$ will be called equivalent if
$\{\x_1, \ldots, \x_r\}=\{\y_1, \ldots, \y_s\}$, i.e.~if $\p$ is the result of permuting
together several segments of $\mathbf{q}$ of the same slope. We denote by
$\mathbf{Conv}$ the set of equivalence classes of convex paths in
$\mathbf{Conv}'$. We shall  only consider convex paths up to equivalence, and we shall simply refer to elements of $\mathbf{Conv}$
as ``paths''.
 We also introduce $\mathbf{Conv}^+$ (resp. $\mathbf{Conv}^-$) as the set
of convex paths $(\x_1, \ldots, \x_s)$ satisfying $\widehat{\x_1 L_0} \geq \cdots \geq \widehat{\x_s L_0}
\geq \pi$ (resp. $\pi>\widehat{\x_1 L_0} \geq \cdots \geq \widehat{\x_s L_0}
\geq 0$). Concatenation of paths then yields an identification $\mathbf{Conv} \simeq \mathbf{Conv}^+ \times \mathbf{Conv}^-$.
\centerline{
\begin{picture}(36,15)(0,-2.5)
\multiput(0,0.5)(1,0){11}{\circle*{.2}}
\multiput(0,1.5)(1,0){11}{\circle*{.2}}
\multiput(0,2.5)(1,0){11}{\circle*{.2}}
\multiput(0,3.5)(1,0){11}{\circle*{.2}}
\multiput(0,4.5)(1,0){11}{\circle*{.2}}
\multiput(0,5.5)(1,0){11}{\circle*{.2}}
\multiput(0,6.5)(1,0){11}{\circle*{.2}}
\multiput(0,7.5)(1,0){11}{\circle*{.2}}
\multiput(0,8.5)(1,0){11}{\circle*{.2}}
\multiput(0,9.5)(1,0){11}{\circle*{.2}}
\multiput(0,10.5)(1,0){11}{\circle*{.2}}
\multiput(12.5,0.5)(1,0){11}{\circle*{.2}}
\multiput(12.5,1.5)(1,0){11}{\circle*{.2}}
\multiput(12.5,2.5)(1,0){11}{\circle*{.2}}
\multiput(12.5,3.5)(1,0){11}{\circle*{.2}}
\multiput(12.5,4.5)(1,0){11}{\circle*{.2}}
\multiput(12.5,5.5)(1,0){11}{\circle*{.2}}
\multiput(12.5,6.5)(1,0){11}{\circle*{.2}}
\multiput(12.5,7.5)(1,0){11}{\circle*{.2}}
\multiput(12.5,8.5)(1,0){11}{\circle*{.2}}
\multiput(12.5,9.5)(1,0){11}{\circle*{.2}}
\multiput(12.5,10.5)(1,0){11}{\circle*{.2}}
\multiput(25,0.5)(1,0){11}{\circle*{.2}}
\multiput(25,1.5)(1,0){11}{\circle*{.2}}
\multiput(25,2.5)(1,0){11}{\circle*{.2}}
\multiput(25,3.5)(1,0){11}{\circle*{.2}}
\multiput(25,4.5)(1,0){11}{\circle*{.2}}
\multiput(25,5.5)(1,0){11}{\circle*{.2}}
\multiput(25,6.5)(1,0){11}{\circle*{.2}}
\multiput(25,7.5)(1,0){11}{\circle*{.2}}
\multiput(25,8.5)(1,0){11}{\circle*{.2}}
\multiput(25,9.5)(1,0){11}{\circle*{.2}}
\multiput(25,10.5)(1,0){11}{\circle*{.2}}
\put(5,5.5){\circle*{.3}}
\put(4.75,4.75){\Small{$\mathbf{o}$}}
\put(17.5,5.5){\circle*{.3}}
\put(17.25,4.75){\Small{$\mathbf{o}$}}
\put(30,5.5){\circle*{.3}}
\put(29.75,4.75){\Small{$\mathbf{o}$}}
\thicklines
\put(5,5.5){\line(1,1){1}}
\put(6,6.5){\circle*{.3}}
\put(6,6.5){\line(2,1){2}}
\put(8,7.5){\circle*{.3}}
\put(8,7.5){\line(-1,-1){1}}
\put(7,6.5){\circle*{.3}}
\put(7,6.5){\line(1,0){2}}
\put(9,6.5){\circle*{.3}}
\put(9.2,6.8){\Small{$\mathbf{p}_1$}}
\put(5,5.5){\line(1,-2){1}}
\put(6,3.5){\circle*{.3}}
\put(6,3.5){\line(1,-1){2}}
\put(8,1.5){\circle*{.3}}
\put(8,1.5){\line(1,3){1}}
\put(9,4.5){\circle*{.3}}
\put(9,4.5){\line(0,1){1}}
\put(9,5.5){\circle*{.3}}
\put(9.2,4.75){\Small{$\mathbf{p}_2$}}
\put(5,5.5){\line(-1,4){1}}
\put(4,9.5){\circle*{.3}}
\put(4,9.5){\line(-1,0){3}}
\put(1,9.5){\circle*{.3}}
\put(1,9.5){\line(0,-1){4}}
\put(1,5.5){\circle*{.3}}
\put(1.2,5.75){\Small{$\mathbf{p}_3$}}
\put(17.5,5.5){\line(-1,1){1}}
\put(16.5,6.5){\circle*{.3}}
\put(16.5,6.5){\line(-1,-1){2}}
\put(14.5,4.5){\circle*{.3}}
\put(14.5,4.5){\line(1,-1){4}}
\put(18.5,.5){\circle*{.3}}
\put(18.7,.8){\Small{$\mathbf{p}_4$}}
\put(17.5,5.5){\line(0,1){3}}
\put(17.5,8.5){\circle*{.3}}
\put(17.5,8.5){\line(-1,0){4}}
\put(16.5,8.5){\circle*{.3}}
\put(13.5,8.5){\circle*{.3}}
\put(13.5,8.5){\line(0,-1){6}}
\put(13.5,2.5){\circle*{.3}}
\put(13.25,1.8){\Small{$\mathbf{p}_5$}}
\put(30,5.5){\line(1,-1){2}}
\put(32,3.5){\circle*{.3}}
\put(32,3.5){\line(1,2){1}}
\put(33,5.5){\circle*{.3}}
\put(33,5.5){\line(0,1){2}}
\put(33,7.5){\circle*{.3}}
\put(33,7.5){\line(-1,0){4}}
\put(29,7.5){\circle*{.3}}
\put(29,7.5){\line(0,-1){5}}
\put(29,2.5){\circle*{.3}}
\put(29,2.5){\line(1,-1){1}}
\put(30,1.5){\circle*{.3}}
\put(30.2,1.05){\Small{$\mathbf{p}_6$}}
\put(.5,-1){\Small{$\mathbf{p}_1$ is not convex}}
\put(.5,-2){\Small{$\mathbf{p}_2 \in \mathbf{Conv}^+$,$\mathbf{p}_3 \in \mathbf{Conv}^-$}}
\put(12,-1){\Small{$\mathbf{p}_4$ is convex but $\mathbf{p}_4 \not\in \mathbf{Conv}$}}
\put(12,-2){\Small{$\mathbf{p}_5 \in \mathbf{Conv}$, but $\mathbf{p}_5 \not\in \mathbf{Conv}^{\pm}$}}
\put(25.5,-1){\Small{$\mathbf{p}_6$ is not convex}}
\end{picture}}
\centerline{Figure 3. Examples of paths}

\vspace{.05in}

Observe that distinct paths could give rise to the same polygonal line in $\ZZ$: for instance $\mathbf{p}=\bigl((0,1),(0,1)\bigr)$ and
$\mathbf{p}'=\bigl((0,2)\bigr)$.
To a path $\p=(\x_1, \ldots, \x_r)$ we associate the element
$T_{\p}:=T_{\x_1} \cdots T_{\x_r} \in \mathbf{U}_{\E}$.
This expression is well-defined since  $\U^{\pm, (\mu)}_{\E}$ is
commutative for all slopes $\mu$.
Moreover, it follows from Theorem \ref{L:bial} that the set of elements 
$\bigl\{T_\p\;|\p \in \mathbf{Conv}^{\pm}\bigr\}$ is a
$\KK$-basis of $\U^{\pm}_{\E}$.

\vspace{.1in}

\addtocounter{theo}{1}
\paragraph{\textbf{Remark \thetheo}}  The group $SL(2,\Z)$ naturally acts on the set of paths.
For any ray  $L$ in $\ZZ$ starting at the origin we can  define the set $\mathbf{Conv}^L$ by replacing $L_0$ by $L$, and any
$\sigma \in SL(2,\Z)$ maps bijectively $\mathbf{Conv}^L$ to $\mathbf{Conv}^{\sigma(L)}$.
In particular, $\{T_{\p}\;|\p \in \mathbf{Conv}^L\}$ is a $\KK$-basis of $\U_{\E}$ for any $L$. Such a
choice of $L$ corresponds to a choice of a $t$-structure in the derived category $D^b\bigl(Coh({\E})\bigr)$.

\vspace{.1in}

For $\x,\y \in \ZZ^*$ we let $\boldsymbol{\Delta}_{\x,\y}$ stand for
the triangle with corners  $\mathbf{o}, \mathbf{x}, \mathbf{x+y}$. If
$\widehat{\x\y} < \pi$ then $T_{\y}T_{\x}$ (corresponding to the path
$(\y,\x)$) can  be written as a linear combination of elements
$T_{\p}$ where  $\p$ runs through
 the set of \textit{convex} paths lying in
$\boldsymbol{\Delta}_{\x,\y}$. Indeed, this is a reformulation of Remark~2.7 when $\x,\y \in \ZZ^+$, and follows for an
arbitrary pair $(\x,\y)$ by ${SL}(2,\Z)$-invariance of ${\U}_{\E}$.

\centerline{
\begin{picture}(28,12)
\multiput(0,0.5)(1,0){11}{\circle*{.2}}
\multiput(0,1.5)(1,0){11}{\circle*{.2}}
\multiput(0,2.5)(1,0){11}{\circle*{.2}}
\multiput(0,3.5)(1,0){11}{\circle*{.2}}
\multiput(0,4.5)(1,0){11}{\circle*{.2}}
\multiput(0,5.5)(1,0){11}{\circle*{.2}}
\multiput(0,6.5)(1,0){11}{\circle*{.2}}
\multiput(0,7.5)(1,0){11}{\circle*{.2}}
\multiput(0,8.5)(1,0){11}{\circle*{.2}}
\multiput(0,9.5)(1,0){11}{\circle*{.2}}
\multiput(0,10.5)(1,0){11}{\circle*{.2}}
\multiput(17,0.5)(1,0){11}{\circle*{.2}}
\multiput(17,1.5)(1,0){11}{\circle*{.2}}
\multiput(17,2.5)(1,0){11}{\circle*{.2}}
\multiput(17,3.5)(1,0){11}{\circle*{.2}}
\multiput(17,4.5)(1,0){11}{\circle*{.2}}
\multiput(17,5.5)(1,0){11}{\circle*{.2}}
\multiput(17,6.5)(1,0){11}{\circle*{.2}}
\multiput(17,7.5)(1,0){11}{\circle*{.2}}
\multiput(17,8.5)(1,0){11}{\circle*{.2}}
\multiput(17,9.5)(1,0){11}{\circle*{.2}}
\multiput(17,10.5)(1,0){11}{\circle*{.2}}
\put(5,5.5){\circle*{.3}}
\put(22,5.5){\circle*{.3}}
\put(4.75,4.75){\Small{$\mathbf{o}$}}
\put(21.75,4.75){\Small{$\mathbf{o}$}}
\put(5,5.5){\line(1,4){1}}
\put(6,9.5){\circle*{.3}}
\put(5.75,9.85){\Small{$\mathbf{y}$}}
\put(6,9.5){\line(4,-3){4}}
\put(10,6.5){\circle*{.3}}
\put(9.5,6.75){\Small{$\mathbf{x}+\mathbf{y}$}}
\put(5,5.5){\line(4,-3){4}}
\put(9,2.5){\circle*{.3}}
\put(8.75,1.75){\Small{$\mathbf{x}$}}
\put(9,2.5){\line(1,4){1}}
\put(5,5.5){\line(5,1){5}}
\put(7,4.75){\Small{$\boldsymbol{\Delta}_{\mathbf{x},\mathbf{y}}$}}
\put(22,5.5){\line(1,4){1}}
\put(23,9.5){\circle*{.3}}
\put(22.75,9.85){\Small{$\mathbf{y}$}}
\put(23,9.5){\line(4,-3){4}}
\put(27,6.5){\circle*{.3}}
\put(26.5,6.75){\Small{$\mathbf{x}+\mathbf{y}$}}
\put(26,2.5){\circle*{.3}}
\put(25.75,1.75){\Small{$\mathbf{x}$}}
\put(22,5.5){\line(5,1){5}}
\thicklines
\put(22,5.5){\line(1,0){3}}
\put(25,5.5){\circle*{.3}}
\put(25,5.5){\line(2,1){2}}
\put(22,5.5){\line(3,-1){3}}
\put(25,4.5){\circle*{.3}}
\put(25,4.5){\line(1,0){1}}
\put(26,4.5){\circle*{.3}}
\put(26,4.5){\line(1,2){1}}
\put(26,2.5){\circle*{.3}}
\put(26,3.5){\line(1,3){1}}
\put(26,3.5){\circle*{.3}}
\put(22,5.5){\line(4,-3){4}}
\put(26,2.5){\line(1,4){1}}
\put(22,5.5){\line(2,-1){4}}
\end{picture}}
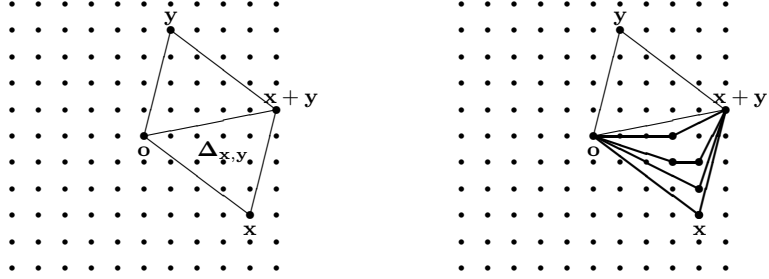
\centerline{Figure 4. The triangle $\boldsymbol{\Delta}_{\x,\y}$ and some convex paths in it}

\vspace{.1in}

Several arguments in this Section are based on Pick's formula, which we recall~: for any pair of non-colinear points $\x,\y \in \ZZ$ :
\begin{equation}\label{E:Pick}
|\det(\x,\y)| = \deg(\x) + \deg(\y) + \deg(\x+\y) -2 + 2 \#(\Delta_{\x,\y} \cap \ZZ).
\end{equation}

\vspace{.2in}

\paragraph{\textbf{5.2.}} We shall  describe ${\U}_{\E}$ as an abstract algebra
of paths modulo a minimal set of ``straightening'' relations given below.
If $\x=(q,p) \in \ZZ^*$ we write $\deg(\x)= \mathsf{gcd}(q,p) \in \N$. For non-collinear
$\x,\y \in \ZZ^*$ we set $\epsilon_{\x,\y}=\mathsf{sign}\bigl(\det(\x,\y)\bigr) \in \{\pm 1\}$.

\vspace{.15in}

\addtocounter{theo}{1}
\paragraph{\textbf{Definition \thetheo}} Fix $\s,\bs \in \C^*$ with $\sigma, \bs \not\in \{\pm 1\}$ and set
$\nu=(\s \bs)^{-1/2}$ and
 $$c_i(\s,\bs)=
(\s^{i/2}-\s^{-i/2})(\bs^{i/2}-\bs^{-i/2})[i]_{\nu}/i.$$
Let
$\UU_{\s,\bs}$ be the $\C$-algebra generated by $\{t_{\x}\;|\x \in \ZZ^*\} $ modulo the following set of relations
\begin{enumerate}
\item[i)] If $\x,\x'$ belong to the same line in $\ZZ$ then
$$[t_\x,t_{\x'}]=0,$$
\item[ii)] Assume that $\x,\y \in \ZZ^*$ are such that $\deg(\x)=1$ and that
$\boldsymbol{\Delta}_{\x,\y}$ has no interior lattice point. Then
$$[t_\y,t_{\x}]=\epsilon_{\x,\y}c_{\deg(\y)}(\s,\bs)
\frac{\theta_{\x+\y}}{\nu^{-1}-\nu}$$
where  the elements $\theta_{\z}$, $\z \in \ZZ^*$ are defined by the following generating series
\begin{equation}\label{E:formulatheta}
\sum_i \theta_{i\x_0}s^i = \exp\bigg((\nu^{-1}-\nu)\sum_{r \geq 1}t_{r\x_0}s^r\bigg),
\end{equation}
for any $\x_0 \in \ZZ^*$ such that $\deg(\x_0)=1$.
\end{enumerate}

\vspace{.15in}
\noindent
Observe that $\theta_{\z}=(\nu^{-1}-\nu)t_{\z}$ whenever $\deg(\z)=1$.
We also denote by $\UU^{\pm}_{\s,\bs}$ the subalgebra of $\UU_{\s,\bs}$ generated by $t_{\x}$ for
$\x \in \ZZ^{\pm}$.

\vspace{.1in}

\begin{lem} For any $\gamma \in  {SL}(2,\Z)$ we have an algebra automorphism 
$\Phi_\gamma: \UU_{\s,\bs} \rightarrow  \UU_{\s,\bs}$ given by the formula 
$\Phi_\gamma(t_{\x})= t_{\gamma(\x)}$  for any $\x \in \ZZ^*$. 
\end{lem}
\noindent
\textit{Proof.} Obvious.\qed

\vspace{.2in}

\paragraph{\textbf{5.3.}} Now let $\#\E(\mathbb{F}_{q^r})$ stand for the number of rational points of $\E$ over
$\mathbb{F}_{q^r}$ and recall that $v = q^{-1/2}$.  By a theorem of Hasse (see e.g.~
\cite{Ha}, App. C) there exist conjugate algebraic numbers $\sigma, \overline{\sigma}$, satisfying $\sigma \overline{\sigma}=q$,
such that $$\#\E(\mathbb{F}_{q^r})=q^r+1-(\sigma^r + \overline{\sigma}^r)$$ for any $r \geq 1$.
These  numbers $\sigma, \overline{\sigma}$ are the eigenvalues of the
Frobenius automorphism  acting on $H^1\bigl(\E_{\overline{\mathbb{F}_q}}, \qlb\bigr)$.
Note that $c_i(\s,\bs) = v^i[i]\#\E(\mathbb{F}_{q^i})/i=c_i(\E)$.

\begin{theo}\label{T:main2} The assignment $\Omega: t_{\x} \mapsto T_\x $ for $\x \in \ZZ^*$ extends to an isomorphism
$\Omega: \UU_{\s,\bs} \simeq {\U}_{\E} \otimes_{\KK} \C$.
\end{theo}

\vspace{.1in}

\addtocounter{theo}{1}
\paragraph{\textbf{Example \thetheo}} Before giving the proof of this theorem, let us illustrate the use of the straightening
relation ii). We shall compute $t_{(1,2)}t_{(1,-1)}$, which corresponds to the path 
$\bigl((1,2),(1,-1)\bigr)$ not
belonging to $\mathbf{Conv}$. By ii) we have $\bigl[t_{(0,1)},t_{(1,1)}\bigr]=c_1  t_{(1,2)}$, hence
\begin{equation*}
\begin{split}
\bigl[t_{(1,2)},t_{(1,-1)}\bigr]&=\frac{1}{c_1}\Bigl\{ \bigl[[t_{(0,1)},t_{(1,-1)}],t_{(1,1)}\bigr] +
\bigl[t_{(0,1)},[t_{(1,1)},t_{(1,-1)}]\bigr]\Bigr\}\\
&=[t_{(1,0)},t_{(1,1)}] + \bigl[t_{(0,1)},t_{(2,0)}+\frac{1}{2}(v^{-1}-v)t_{(1,0)}^2\bigr]
\end{split}
\end{equation*}
where we have used ii) in each term, and the relation $$\frac{\displaystyle
\theta_{(2,0)}}{\displaystyle v^{-1}-v}=t_{(2,0)}+
\frac{1}{2}(v^{-1}-v)t_{(1,0)}^2.$$
Now, by ii) again we have
$\bigl[t_{(1,0)},t_{(1,1)}\bigr]=-c_1t_{(2,1)}$, 
$\bigl[t_{(0,1)},t_{(2,0)}\bigr]=c_2t_{(2,1)}$
and $\bigl[t_{(0,1)},t_{(1,0)}\bigr]=t_{(1,1)}$. Hence, we obtain
\begin{equation*}
\begin{split}
\bigl[t_{(1,2)},t_{(1,-1)}\bigr]&=(c_2-c_1)t_{(2,1)}+\frac{1}{2}(v^{-1}-v)c_1(t_{(1,1)}t_{(1,0)}+
t_{(1,0)}t_{(1,1)})\\
&=(c_2-c_1)t_{(2,1)} +
\frac{1}{2}(v^{-1}-v)c_1\big(c_1 t_{(2,1)}+2t_{(1,0)}t_{(1,1)}\big)
\end{split}
\end{equation*}
Gathering terms, we get
$$t_{(1,2)}t_{(1,-1)}=
t_{(1,-1)}t_{(1,2)} +\frac{1}{2}(v^{-1}-v)c_1t_{(1,0)}t_{(1,1)}+
c_1\bigl([3]-vc_1\bigr)t_{2,1}.
$$
Observe that all three paths $\bigl((1,-1),(1,2)\bigr), \bigl((1,0),(1,1)\bigr)$,
 $\bigl((2,1)\bigr)$ belong to $\mathbf{Conv}$.

\vspace{.2in}

We begin the proof of Theorem~\ref{T:main2}. Let us first show that the map $\Omega$ is well-defined, i.e.~that
relations i) and ii) hold in $\U_\E$. By the ${SL}(2,\Z)$-invariance of
$\UU_{\s,\bs}$ and $\U_{\E}$ it is enough to prove relation i) for
$\x=(0,r), \x'=(0,r')$. The subalgebra $\mathbf{H}^{(\infty)}_{\E}$ of $\mathbf{H}_{\E}$
is stable under the coproduct (as any
subsheaf or quotient of a torsion sheaf is again a torsion sheaf) and can be described as the product over all points $ x \in \E$
of the Hall bialgebras of the categories $\mathcal{N}_{\kk_x}$. By Proposition~\ref{P:Mac},~ii),
$\Delta\bigl(T^{(\infty)}_{r,x}\bigr)=T^{(\infty)}_{r,x} \otimes 1 + 1 \otimes T^{(\infty)}_{r,x} $. Hence, from the definition of the Drinfeld double
we get $\bigl[T_{(0,r)},T_{(0,r')}\bigr]=0$ as desired.

\vspace{.1in}

Let us prove the relation ii). Assume that $\x,\y$ are as in ii). Since $\deg(\x)=1$ we cannot have
$\deg(\y) = \deg(\x+\y)=2$. On the other hand, it is easy to see that if $\deg(\y) \geq 2$ and $\deg(\x+\y) \geq 3$, or if $\deg(\x+\y) \geq 2$
and $\deg(\y) \geq 3$ then
$\boldsymbol{\Delta}_{\x,\y}$ contains interior lattice points. In conclusion, we either have $\deg(\y)=1$ or $\deg(\x+\y)=1$. We split
our argument according to this dichotomy.

\vspace{.05in}

\noindent
\textit{Case a.1}. We have $\deg(\x+\y)=1$ and $\epsilon_{\x,\y} >0$. Up to the ${SL}(2,\Z)$-action, we
may fix $\x=(1,0)$ and if $\det(\x,\y)=r$ then we may furthermore assume that $\y=(s,r)$ for some $0 \leq s <r$. Using Pick's formula (\ref{E:Pick}), we deduce that
there are no points inside $\boldsymbol{\Delta}_{\x,\y}$ if and only if $\deg(\y)=r$, which implies
 $\y=(0,r)$. Then relation ii) follows from Lemma~\ref{L:comp}.

\vspace{.1in}

\noindent
\textit{Case a.2} We have $\deg(\x+\y)=1$ and $\epsilon_{\x,\y}<0$. Without loss of generality, we may  assume that $\x=(r_1,d_1),
\y=(r_2,d_2)$ with $r_1 >0$ and $r_2>0$. Now let us use the antiautomorphism $D$ of Proposition~\ref{T:Duali}. Note that
$D(T_{(r,d)})=T_{(r,-d)}$ and
$\epsilon_{D(\x),D(\y)}>0$ hence the desired relation follows from case a.1 above.

\vspace{.1in}

\noindent
\textit{Case b.} We have $\deg(\y)=1$. In that situation, simple application of Pick's formula (\ref{E:Pick}) shows that
$\deg(\x+\y)=\big|\det(\x,\y)\big|$, and, after exchanging the role of $\x$ and $\y$ if necessary and using the
${SL}(2,\mathbb{Z})$ invariance we may assume that $\x=(1,n), \y=(-1,l)$.
The expression for the commutator $\bigl[T_{\y},T_{\x}\bigr]$ can be now
derived from the definition of the Drinfeld double together with Lemma~\ref{L:coprodun}~:
if $n+l>0$ then $\epsilon_{(1,n),(-1,l)}=1$ and the relation $R\bigl(T_{(1,n)},T_{(-1,l)}\bigr)$ is
\begin{equation*}
T_{(-1,l)}T_{(1,n)}=T_{(1,n)}T_{(-1,l)}+\Theta_{(0,n+l)}(T_{(-1,l)},T_{(-1,l)})=
T_{(1,n)}T_{(-1,l)}+c_1\frac{\Theta_{(0,n+l)}}{v^{-1}-v},
\end{equation*}
and if $n+l<0$ then $\epsilon_{(1,n),(-1,l)}=-1$ and the relation $R\bigl(T_{(1,n)},T_{(-1,l)}\bigr)$ is
\begin{equation*}
T_{(1,n)}T_{(-1,l)}=T_{(-1,l)}T_{(1,n)}+\Theta_{(0,n+l)}\bigl(T_{(-1,l)},T_{(-1,l)}\bigr)=
T_{(-1,l)}T_{(1,n)}+c_1\frac{\Theta_{(0,n+l)}}{v^{-1}-v}.
\end{equation*}
This concludes the proof of relation ii).
\vspace{.1in}

By the above, $\Omega$ is well-defined and extends to a surjective algebra morphism $\Omega: \UU_{\s,\bs} \tto \U_{\E}\otimes \C$.
Moreover, this morphism is ${SL}(2,\Z)$-equivariant.
 In the rest of the proof, we construct an inverse of $\Omega$.  We first concentrate on the ``positive'' subalgebra $\UU^+_{\s,\bs}$
of $\UU_{\s,\bs}$. For any path $\mathbf{p}=(\x_1, \ldots, \x_r)$ we set $t_{\mathbf{p}}
=t_{\x_1}\cdots t_{\x_r}$. Note that from the surjectivity of $\Omega$ and
Proposition~\ref{L:bial} it follows that the elements
$\bigl\{t_{\mathbf{p}}\;|\; \mathbf{p}\in \mathbf{Conv}^+\bigr\}$ are linearly independent.

\vspace{.1in}

\begin{lem}\label{L:Aplus} The subalgebra $\UU^+_{\s,\bs}$ is equal to $\bigoplus_{\mathbf{p} \in \mathbf{Conv}^+} \C t_{\mathbf{p}}$.
\end{lem}
\noindent
\textit{Proof.}
The inclusion is obvious in one direction. For the other inclusion, we have
to show that  any path $\mathbf{p}$ in $\ZZ^+$  can  be ``straightened''
using the relations  ii). By an argument, which is at all steps  similar to the
proof of Lemma~\ref{L:relH},
it is sufficient to show that for any $\x,\y \in \ZZ^+$ with $\mu(\y) > \mu(\x)$, we have
\begin{equation}\label{E:hypderec}
t_{\y}t_{\x} \in \bigoplus_{\mathbf{p} \in I_{\x,\y}} \C t_{\mathbf{p}},
\end{equation}
where by definition $I_{\x,\y}$ is the set of convex paths in $\boldsymbol{\Delta}_{\x,\y}$ joining $\mathbf{o}$ to $\x+\y$.
We shall achieve this by induction on $\det(\y,\x)$.

If $\det(\y,\x)=1$ then (e.g.~by Pick's formula) $\deg(\x)=\deg(\y)=\deg(\x+\y)=1$ thus $t_{\y}t_{\x}=t_{\x}t_{\y} + c_1 t_{\x+\y}$ by
relation ii).
So let us fix $d>1$ and let us assume that (\ref{E:hypderec}) holds
for any $\x',\y'$ satisfying $\det(\y',\x')<d$.

If $\mathbf{p}=(\x_1, \ldots, \x_r)$ is any path in $\ZZ^+$ we put $\mathbf{p}^{\#}=(\x_{\sigma(1)}, \ldots, \x_{\sigma(r)})$ where
$\sigma $ is the least length permutation satisfying $\mu(\x_{\sigma(1)}) \leq
\mu(\x_{\sigma(2)}) \leq \cdots \leq \mu(\x_{\sigma(r)})$, and
we denote by $a(\mathbf{p})$ the area of the polygon bounded by $\mathbf{p}$ and $\mathbf{p}^{\#}$. Observe that if $\mathbf{p}'$ is a
subpath of $\mathbf{p}$ then $a(\mathbf{p}') \leq a(\mathbf{p})$. Also, if $\z, \mathbf{w} \in \ZZ^+$ are such that $\mu(\z) >
\mu(\mathbf{w})$ then $a\bigl((\z,\mathbf{w})\bigr) = \det(\mathbf{z},\mathbf{w})$.

\centerline{
\begin{picture}(10,10)
\multiput(0,0.5)(1,0){11}{\circle*{.2}}
\multiput(0,1.5)(1,0){11}{\circle*{.2}}
\multiput(0,2.5)(1,0){11}{\circle*{.2}}
\multiput(0,3.5)(1,0){11}{\circle*{.2}}
\multiput(0,4.5)(1,0){11}{\circle*{.2}}
\multiput(0,5.5)(1,0){11}{\circle*{.2}}
\multiput(0,6.5)(1,0){11}{\circle*{.2}}
\multiput(0,7.5)(1,0){11}{\circle*{.2}}
\multiput(0,8.5)(1,0){11}{\circle*{.2}}
\put(1,5.5){\circle*{.3}}
\put(.75,5.75){\Small{$\mathbf{o}$}}
\thicklines
\put(1,5.5){\line(1,0){1}}
\put(2,5.5){\circle*{.3}}
\put(2,5.5){\line(1,2){1}}
\put(3,7.5){\circle*{.3}}
\put(3,7.5){\line(1,-2){1}}
\put(4,5.5){\circle*{.3}}
\put(4,5.5){\line(1,-1){2}}
\put(6,3.5){\circle*{.3}}
\put(6,3.5){\line(2,1){2}}
\put(8,4.5){\circle*{.3}}
\put(8,4.5){\line(1,0){2}}
\put(10,4.5){\circle*{.3}}
\put(10,4.5){\line(0,1){1}}
\put(10,5.5){\circle*{.3}}
\put(1,5.5){\line(1,-2){1}}
\put(2,3.5){\circle*{.3}}
\put(2,3.5){\line(1,-1){2}}
\put(4,1.5){\circle*{.3}}
\put(4,1.5){\line(1,0){3}}
\put(5,1.5){\circle*{.3}}
\put(7,1.5){\circle*{.3}}
\put(7,1.5){\line(2,1){2}}
\put(9,2.5){\circle*{.3}}
\put(9,2.5){\line(1,2){1}}
\put(5.25,.7){\Small{$\mathbf{p}^{\#}$}}
\put(5.25,4.75){\Small{$\mathbf{p}$}}
\put(4.5,2.7){\Small{$a(\mathbf{p})$}}
\end{picture}}
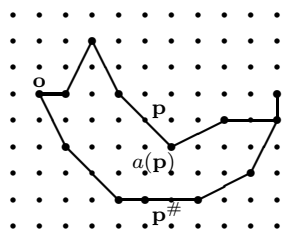
\centerline{Figure 5. The area $a(\mathbf{p})$ of a path in $\ZZ^+$.}

\vspace{.1in}

\noindent
\textbf{Claim.} For any path $\mathbf{p}$ in $\ZZ^+$ satisfying $a(\mathbf{p}) <d$ we have $t_{\mathbf{p}}
\in \bigoplus_{\mathbf{p} \in \mathbf{Conv}^+} \C t_{\mathbf{p}}$.\\
\noindent
\textit{Proof of Claim}. The assertion is true by definition if $a(\p)=0$. If $a(\p) >0$ then $\p=(\x_1, \ldots, \x_r)$ with
$\mu(\x_1) \leq \cdots \leq \mu(\x_s) > \mu(\x_{s+1})$ for some $s$. We have $\det(\x_{s}, \x_{s+1}) \leq a(\p) <d$ hence
$t_{\x_s} t_{\x_{s+1}}=\sum_i u_i t_{\mathbf{q}_i}$ for some $\mathbf{q}_i \in I_{\x_{s+1},\x_{s}}$ and, setting
$\mathbf{p}_i=(\x_1, \ldots, \x_{s-1}, \mathbf{q}_i, \x_{s+2}, \ldots, \x_{r})$ we get  $t_{\p}=\sum_{i} u_i t_{\p_i}$. It is clear
that for all $i$ both $\mathbf{p}_i$ and $\mathbf{p}_i^{\#}$ strictly lie inside the polygon bounded by $\p$ and $\p^{\#}$, so that
$a(\p_i) < a(\p)$.

\centerline{
\begin{picture}(15,10)
\multiput(0,0.5)(1,0){11}{\circle*{.2}}
\multiput(0,1.5)(1,0){11}{\circle*{.2}}
\multiput(0,2.5)(1,0){11}{\circle*{.2}}
\multiput(0,3.5)(1,0){11}{\circle*{.2}}
\multiput(0,4.5)(1,0){11}{\circle*{.2}}
\multiput(0,5.5)(1,0){11}{\circle*{.2}}
\multiput(0,6.5)(1,0){11}{\circle*{.2}}
\multiput(0,7.5)(1,0){11}{\circle*{.2}}
\multiput(0,8.5)(1,0){11}{\circle*{.2}}
\put(1,5.5){\circle*{.3}}
\put(.75,5.75){\Small{$\mathbf{o}$}}
\thicklines
\put(1,5.5){\line(1,0){1}}
\put(2,5.5){\circle*{.3}}
\put(2,5.5){\line(1,2){1}}
\put(3,7.5){\circle*{.3}}
\put(3,7.5){\line(1,-2){1}}
\put(4,5.5){\circle*{.3}}
\put(4,5.5){\line(1,-1){2}}
\put(6,3.5){\circle*{.3}}
\put(6,3.5){\line(2,1){2}}
\put(8,4.5){\circle*{.3}}
\put(8,4.5){\line(1,0){2}}
\put(10,4.5){\circle*{.3}}
\put(10,4.5){\line(0,1){1}}
\put(10,5.5){\circle*{.3}}
\put(1,5.5){\line(1,-2){1}}
\put(2,3.5){\circle*{.3}}
\put(2,3.5){\line(1,-1){2}}
\put(4,1.5){\circle*{.3}}
\put(4,1.5){\line(1,0){3}}
\put(5,1.5){\circle*{.3}}
\put(7,1.5){\circle*{.3}}
\put(7,1.5){\line(2,1){2}}
\put(9,2.5){\circle*{.3}}
\put(9,2.5){\line(1,2){1}}
\thinlines
\multiput(2,5.5)(.5,0){4}{\line(1,0){.5}}
\put(1,5.5){\line(1,-1){2}}
\put(3,3.5){\line(1,0){5}}
\put(8,3.5){\line(2,1){2}}
\put(3,3.5){\circle*{.3}}
\put(4,3.5){\circle*{.3}}
\put(6,3.5){\circle*{.3}}
\thicklines
\put(12,5){\line(1,0){1}}
\put(13.25,5){\Small{:\;$\mathbf{p}$,$\mathbf{p}^{\#}$}}
\thinlines
\put(12,4){\line(1,0){1}}
\put(13.25,4){\Small{:\;$\mathbf{p}_i$,$\mathbf{p}_i^{\#}$}}
\end{picture}}
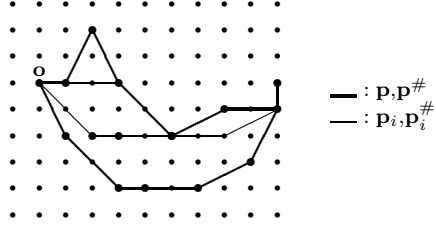
\centerline{Figure 6. The area of a path $\mathbf{p}$ before and after one straightening.}

\vspace{.1in}
\noindent
We may iterate this process until we are only left with paths $\mathbf{q}$ satisfying $a(\mathbf{q})=~0$. Hence $t_{\mathbf{p}}
\in \bigoplus_{\mathbf{p} \in \mathbf{Conv}^+} \C t_{\mathbf{p}}$ and the claim is proven.\qed

\vspace{.1in}

Now let us fix $\x,\y$ such that $\mu(\y)>\mu(\x)$ and $\det(\y,\x)=d$.
If $\boldsymbol{\Delta}_{\x,\y}$ has no interior lattice point then
(see the proof of relation ii) above) either $\deg(\x)=\deg(\y)=\deg(\x+\y)=2$, or $\deg(\x)=1$ or
$\deg(\y)=1$. In the first case, we can  assume up to the
${SL}(2,\Z)$-action that $\y=(2,0)$ and $\x=(0,2)$. We
leave  to the reader to check that repeated applications
of ii) as in Example~5.5 lead to the equality
$$
t_{(0,2)}t_{(2,0)}
=t_{(2,0)}t_{(0,2)} +  c  t_{(1,1)}^2 +
c_2\big(\frac{c_2}{c_1}-2\big)t_{(2,2)},
$$
where $c =
\frac{(v^{-1}-v)}{2}
\left(\frac{c_2}{c_1}(c_2 + \frac{c_2}{c_1}-1)+  \frac{(v^{-1}-v)}{2}c_2(1-c_1) \right).$
In the last two cases, relation ii) directly yields (\ref{E:hypderec}). So we may assume that $\boldsymbol{\Delta}_{\x,\y}$
contains interior lattice points.

Let us choose $\z \in \boldsymbol{\Delta}_{\x,\y}$ so that the triangle $\mathbf{o}\z\x$ has no interior points and $\deg(\z)= \deg(\x-\z)=1$.
Note that (\ref{E:hypderec}) is stable under the action of ${SL}(2, \Z)$, hence without loss of generality we can assume that $\x-\z \in \ZZ^+$.
By construction, $\z$ and  $\x-\z$ satisfy both conditions of the  relation ii),  hence
$[t_{\z},t_{\x-\z}]=c_{1}\frac{\theta_{\x}}{v^{-1}-v}=c_{1} {t_\x} +u$ for some $u$ belonging to the subalgebra
$\bigl\langle t_{\x_0}, \ldots, t_{(\deg(\x)-1)\x_0} \bigr\rangle$ generated by $t_{\x_0}, \ldots, t_{(\deg(\x)-1)\x_0}$, where
$\x_0=\frac{\x}{\deg(\x)}$.  Therefore,
\begin{equation}\label{E:xyz}
c_1[t_{\y},t_{\x}]=\bigl[[t_{\y},t_{\z}\bigr],t_{\x-\z}] + 
\bigl[t_{\z},[t_{\y},t_{\x-\z}]\bigr] -[t_{\y},u].
\end{equation}
\centerline{
\begin{picture}(10,9)(0,1)
\multiput(0,1.5)(1,0){11}{\circle*{.2}}
\multiput(0,2.5)(1,0){11}{\circle*{.2}}
\multiput(0,3.5)(1,0){11}{\circle*{.2}}
\multiput(0,4.5)(1,0){11}{\circle*{.2}}
\multiput(0,5.5)(1,0){11}{\circle*{.2}}
\multiput(0,6.5)(1,0){11}{\circle*{.2}}
\multiput(0,7.5)(1,0){11}{\circle*{.2}}
\multiput(0,8.5)(1,0){11}{\circle*{.2}}
\put(1,4.5){\circle*{.3}}
\put(.75,4.75){\Small{$\mathbf{o}$}}
\thicklines
\put(1,4.5){\line(2,3){2}}
\put(3,7.5){\circle*{.3}}
\put(3,7.5){\line(3,-2){3}}
\put(6,5.5){\circle*{.3}}
\put(1,4.5){\line(5,1){5}}
\put(1,4.5){\line(3,-2){3}}
\put(4,2.5){\circle*{.3}}
\put(4,2.5){\line(2,3){2}}
\thinlines
\put(1,4.5){\line(2,-1){2}}
\put(1,4.5){\line(1,-1){1}}
\put(3,3.5){\circle*{.3}}
\put(2,3.5){\circle*{.3}}
\put(3,3.5){\line(1,-1){1}}
\put(2,3.5){\line(2,-1){2}}
\put(2.75,7.85){\Small{$\mathbf{y}$}}
\put(3.75,1.95){\Small{$\mathbf{x}$}}
\put(2.75,3.85){\Small{$\mathbf{z}$}}
\put(.6,2.85){\Small{$\mathbf{x}-\mathbf{z}$}}
\put(6.1,5.7){\Small{$\mathbf{x}+\mathbf{y}$}}
\end{picture}}
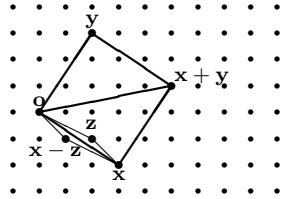
\centerline{Figure 7. The decomposition $\x=\mathbf{z} + (\x-\mathbf{z})$.}

\vspace{.1in}
\noindent
Note that $c_1 \neq 0$ since $|\s|=|\bs|=\sqrt{q}$.
As $\z$ is an interior point of $\boldsymbol{\Delta}_{\x,\y}$ we have $\z=\alpha \x + \beta \y$ for some
$\a,\b \in ]0,1[$ satisfying $\a>\b$. It follows that $\det(\y,\z)<d$ and $\det(\y,\x-\z)<d$, and by the induction hypothesis
$$[t_{\y},t_{\z}] \in \bigoplus_{\p \in I_{\z,\y}} \C t_{\p},\qquad [t_\y,t_{\x-\z}] \in
\bigoplus_{\mathbf{q} \in I_{\x-\z,\y}} \C t_{\mathbf{q}}.$$
Next, as $\mu(\x-\z) < \mu(\z) < \mu(\y)$ we have, for any $\p \in I_{\z,\y}$, 
$\bigl((\x-\z),\p\bigr) \in \mathbf{Conv}^+$ and
$(\p,(\x-\z))^{\#}=\bigl((\x-\z),\p\bigr)$. Thus $a\bigl((\x-\z,\p)\bigr) = 0$ and
$a\bigl((\p,\x-\z)\bigr) = \det(\y+\z,\x-\z)=(1+\b-\a) \det(\y,\x) < d$. It follows by the Claim that $\bigl[[t_{\y},t_{\z}],t_{\x-\z}\bigr] \in
\bigoplus_{\mathbf{Conv}^+} \C t_{\mathbf{p}}$.
In a similar manner, for any $\mathbf{q} \in I_{\x-\z,\y}$ we have  
$a\bigl((\z,\mathbf{q})\bigr) < a\bigl((\z,\y,\x-\z)\bigr)=\det(\y+\z,\x-\z)=
\det(\y,\x)-\det(\x+\y,\z) <d$ since $\det(\x+\y,\z)>0$; and $
a\bigl((\mathbf{q},\z))<a((\y,\x-\z,\z)\bigr)= \det(\y,\x)=d$. Thus
$\bigl[t_{\z},[t_{\y}, t_{\x-\z}]\bigr] \in
\bigoplus_{\mathbf{Conv}^+} \C t_{\mathbf{p}}$.
Finally, let us write $u=\sum_{j=1}^{\deg(\x)-1} a_j t_{j\x_0}$ with $a_j \in
\bigl\langle t_{\x_0}, \ldots, t_{(\deg(\x)-1)\x_0}\bigr\rangle$ of
weight $(\deg(\x)-j)\x_0$. By the induction hypothesis, $t_{\y}a_j \in
\bigoplus_{\p} \C t_{\p}$ where $\p$ ranges in $I_{(\deg(\x)-j)\x_0,\y}$. But as for any such $j$ and $\p$ we have $a\bigl((\p,j\x_0)\bigr) =
\frac{\deg(\x)-j}{\deg(\x)}d <d$ the Claim implies that $t_{\y}u \in \bigoplus_{\mathbf{Conv}^+} \C t_{\mathbf{p}}$.
Hence all together, by (\ref{E:xyz}), $t_{\y}t_{\x} \in \bigoplus_{\mathbf{Conv}^+} \C t_{\mathbf{p}}$.
Finally, let us write $t_{\y}t_{\x}=\sum_{\p \in \mathbf{Conv}^+} c_{\p}t_{\p}$. Applying $\Omega$, we get
$T_{\y}T_{\x}=\sum_{\p} c_{\p} T_{\p}$. By Remark~2.7, we have $T_{\y}T_{\x} \in \bigoplus_{\p \in \mathbf{I}_{\x,\y}}
\C T_{\p}$ so that $c_{\p}=0$ for $\p \not\in \mathbf{I}_{\x,\y}$. Therefore $t_{\y}t_{\x} \in \bigoplus_{\p \in \mathbf{I}_{\x,\y}}
\C t_{\p}$ as desired.
This closes the induction step and proves Lemma~\ref{L:Aplus}.\qed

\vspace{.1in}

Now we are ready to finish the proof of Theorem~\ref{T:main2}. Define $\UU^-_{\s,\bs}$ in the same way as $\UU^+_{\s,\bs}$ by
replacing $\ZZ^+$ by $\ZZ^-$. By Lemma~\ref{L:Aplus}, $\UU^-_{\s,\bs}$ is equal to $\bigoplus_{\mathbf{p} \in
\mathbf{Conv}^-}\C t_{\x}$.
The map $\Omega$ restricts to isomorphisms $\UU^{\pm}_{\s,\bs} \simeq {\U}^{\pm}_{\E}\otimes \C$.
By Theorem~\ref{L:bial} and Corollary~\ref{C:huh3},
${\U}_{\E}$ is generated by
$\mathbf{U}^{\pm}_{\E}$ modulo the collection of relations $R(g,h)$ for sums of classes of
semi-stable sheaves $g \in \U^+_{\E}$ and $h \in \U^-_{\E}$. Now, if
$g$ and $h$ are as above and
$\mu(g)=\mu(h)$ then $R(g,h)$ expresses the fact that $\U^{+,(\mu)}_{\E}$ and $\U^{-,(\mu)}_{\E}$ commute. By
relation~ii), $R\bigl(\Omega^{-1}(g), \Omega^{-1}(h)\bigr)$ holds in $\UU_{\s,\bs}$. If on the other hand $\mu(g) \neq \mu(h)$ then there
exists $\gamma \in {SL}(2,\Z)$ such that ${\gamma}(g), {\gamma}(h) \in \U^+_{\E}$. In that situation, applying ${\gamma}$ to $R(g,h)$ yields a relation $R^{\gamma}\bigl({\gamma} (g),{\gamma}(h)\bigr)$  {in} $\U^+_{\E}$. We deduce that
$R^{\gamma}(\Omega^{-1}\circ{\gamma}(g)), \Omega^{-1}\circ {\gamma}(h)))$ holds in $\UU^+_{\s,\bs}$.  As $\UU_{\s,\bs}$ carries
an action of ${SL}(2,\Z)$ compatible with $\Omega$, it follows that
$R\bigl(\Omega^{-1}(g),\Omega^{-1}(h)\bigr)$ holds in
$\UU_{\s,\bs}$. Therefore, $\Omega^{-1}$ extends to a morphism
$\U_{\E}\otimes \C \to \UU_{\s,\bs}$, which is the desired
inverse to $\Omega$. The theorem is proved.\qed

\vspace{.2in}

\paragraph{\textbf{5.4.}} We still assume that $(\s, \bs)$ is associated to an elliptic curve $\E$. The proof of Theorem~\ref{T:main2} in fact gives the following. Let ${}'\UU^+_{\s,\bs}$ be the $\C$-algebra generated by elements $t_{\x}$ for $\x \in \ZZ^+$ subject to relations i) and ii) of Section~5.2.

\begin{cor} The assignment $\Omega~:t_{\x} \mapsto T_{\x}$ for $\x \in \ZZ^+$ extends to an algebra isomorphism ${}'\UU^+_{\s,\bs} \stackrel{\sim}{\to}\U^+_{\E}$.
In other words, the natural morphism ${}'\UU^+_{\s,\bs} \to \UU^+_{\s,\bs}$ is an isomorphism.
\end{cor}
\vspace{.2in}

\section{Further results~: integral form and  central extension}

\noindent
In this section, we gather some  useful properties of the algebras $\UU_{\s,\bs}$ and $\U_{\E}$.

\vspace{.1in}

\paragraph{\textbf{6.1.}}  For any smooth projective curve $\E$ Kapranov \cite{Kap} considered\footnote{at least implicitly} a natural subalgebra $\H^{sph}_{\E}$ of $\H_{\E}$ which we call the \textit{spherical Hall algebra} of $\E$. By definition, $\H^{sph}_{\E}$ is generated by the elements
$\bigl\{\mathbf{1}_{(0,d)}\;|\; d \in \N\bigr\} \cup
\bigl\{\mathbf{1}^{\mathsf{ss}}_{(1,l)}\;|\; l \in \Z\bigr\}$. In the language of automorphic forms used in \cite{Kap}, these generators are the simplest and most natural cuspidal elements of $\H_{\E}$. In the 
case of an elliptic curve $\E$ it turns out that our algebra $\U^+_{\E}$ coincides with $\H^{sph}_{\E}$. This is an easy consequence of the following corollary of Theorem~\ref{T:main2}~:

\vspace{.1in}

\begin{cor}\label{C:Uisfing} The algebra ${\U}_{\E}^+$ is generated by
$\bigl\{T_{\a}\;| \rank(\a) \leq 1\bigr\}$. Similarly, the algebra
${\U}_{\E}$ is generated by either of the following two sets~:
$$\bigl\{T_{(\pm 1,0)}, T_{(0,\pm 1)}\bigr\}, \qquad \bigl\{T_{(1,0)}, T_{(0,1)}, T_{(-1,-1)}
\bigr\}.$$
\end{cor}
\noindent
\textit{Proof.} We prove the first statement by induction. Denote by $\mathfrak{W}$ the subalgebra generated by $\{T_{\a}\;|
\rank(\a) \leq 1\}$ and assume that $T_{r,s} \in \mathfrak{W}$ for any $(r,s) \in \ZZ^+$ with $r < n \geq 2$. Fix $\z=(n,p) \in
\ZZ^+$, and let $\x$ be the point of $\bigl\{(r,s)\:|\;r <n\bigr\}$ closest to the segment $\mathbf{o}\z$. By construction there are no
interior lattice points in $\Delta_{\x,\z-\x}$ and thus $\bigl[T_{\x},T_{\z-\x}\bigr] =u \theta_{\z}$ for some $u \neq 0$. By the induction
hypothesis, we have $T_{\x}, T_{\z-\x} \in \mathfrak{W}$, and $\theta_{\z} \in (v^{-1}-v) T_{\z} \oplus \mathfrak{W}$. We deduce
that $T_{\z} \in \mathfrak{W}$ as wanted.

\vspace{.05in}

Let us  deal with the second assertion. As before,  denote by $\mathfrak{W}$ the subalgebra generated
by $\bigl\{T_{(\pm 1,0)}, T_{(0,\pm 1)}\bigr\}$.
We have, for any $l \in \Z$, $\bigl[T_{(0,\pm 1)}, T_{(1,l)}\bigr]=\pm c_1
T_{(1,l\pm 1)}$ and it follows that $T_{(1,l)} \in \mathfrak{W}$ for any $l \in \Z$. Similarly,
$T_{(-1,l)} \in \mathfrak{W}$ for any
$l \in \Z$. But then, considering commutators $\bigl[T_{(-1,l)},T_{(1,l')}\bigr]$, we have
$\Theta_{(0,n)} \in \mathfrak{W}$ for any $n$ as well.
The subalgebra generated by $\{\Theta_{(0,n)}\}$ and the one generated by $\{T_{(0,n)}\}$ being equal, we see that $\mathfrak{W}$
contains all $T_{(r,n)}$ with $|r| \leq 1$.
Applying the first statement of the corollary, we get
${\U}_{\E}^{\pm}  \subset \mathfrak{W}$, from which we deduce that
$\mathfrak{W}={U}_{\E}$.

\vspace{.05in}

The last statement follows the second statement together with the relations
$\bigl[T_{(1,0)},T_{(-1,-1)}\bigr]=c_1  T_{(0,-1)}$ and
$\bigl[T_{(-1,-1)},T_{(0,1)}\bigr]= c_1  T_{(-1,0)}$.\qed

 \vspace{.15in}

Kapranov exhibited certain relations satisfied by the generators $\bigl\{\mathbf{1}_{(0,d)},\mathbf{1}^{\mathsf{ss}}_{(1,l)}\;|\;d >0, l \in \Z\bigr\}$, for any curve $\E$. These are the so-called \textit{functional equations} for Eisenstein series. When $\E$ is an elliptic curve, they take the following form. Put
$${E}^{+}(t)=\sum_{p \in \Z} \mathbf{1}^{\mathsf{ss}}_{(1,p)}t^p, \qquad \psi^{+}(s)=\sum_{d \geq 0} \mathbf{1}_{(0,d)}s^d.$$
Then (see \cite{Kap}, Thm.~3.3.)
\begin{equation}\label{E:funct1}
{E}^{+}(t_1){E}^{+}(t_2)=\frac{\zeta_{\E}(t_1/t_2)}{\zeta_{\E}(t_2/t_1)}
{E}^{+}(t_2){E}^{+}(t_1)
\end{equation}
\begin{equation}\label{E:funct2}
\psi^{+}(t_1)E^+(t_2)=\zeta(\s^{-1/2}\bs^{-1/2}t_1/t_2)E^+(t_2){E}^{+}(t_1),
\end{equation}
\begin{equation}\label{E:funct3}
\psi^+(t_1)\psi^+(t_2)=\psi^+(t_2)\psi^+(t_1),
\end{equation}
where $\zeta_{\E}(t)=\frac{\displaystyle (1-\sigma t)(1-\bs t)}{\displaystyle (1-t)(1-qt)}$ is the
\emph{zeta function} of $\E$.
It is known however that relations (\ref{E:funct1}\,--\,\ref{E:funct3}) do not exhaust the complete
list of relations of
$\U^+_{\E}$. In other words, if $\mathbb{U}^+_{\E}$ denotes the algebra generated by some elements
$T_{(1,l)}, T_{(0,d)}$ subject to relations (\ref{E:funct1}\,--\,\ref{E:funct3}) above then there is a nontrivial surjective algebra homomorphism $\mathbb{U}^+_{\E} \tto \mathbf{U}^+_{\E}$.
One may hope to use the description of $\mathbf{U}^+_{\E}$ given in this paper to explicitly describe the kernel of this map $\mathbb{U}^+_{\E} \tto \mathbf{U}^+_{\E}$. This appears to us to be a very interesting problem: using Kapranov's interpretation of the Hall algebra in terms of automorphic forms for $GL(n)$ over a function field, elements of this kernel correspond to some new, higher rank relations satisfied by residues of Eisenstein series.

\vspace{.1in}

The reader will find yet another presentation of $\U^+_{\E}$ in \cite{SV2}, Section~9, this time in terms of shuffle (or \textit{Feigin-Odesskii}) algebras.

\vspace{.2in}

\paragraph{\textbf{6.2.}} We have only defined so far an algebra
$\UU_{\s,\bs}$ for complex values of $\s$ and $\bs$. However, it is also natural to consider a
version of $\UU_{\s,\bs}$, where $\s$ and $\bs$ are formal parameters.
Put $\Rb=\C\bigl[\s^{\pm 1/2},\bs^{\pm 1/2}\bigr]$, 
$\Kb= \mathrm{Frac}(\Rb)=\C\bigl(\s^{1/2},\bs^{1/2}\bigr)$ where $\s, \bs$ are now formal variables, and consider the $\Kb$-algebra $\UU_{\Kb}$ generated by
elements $\bigl\{t_{\x}\;|\; \x \in \ZZ^*\bigr\}$ modulo
the relations i) and ii) of Section~5.2.We also set $\vv=(\s\bs)^{-1/2}$, $\tilde{t}_{\x}=t_{\x}/[\deg(\x)]_{\vv}$, and for an
arbitrary path $\mathbf{p}=(\x_1, \ldots, \x_r)$ we put  $\tilde{t}_{\mathbf{p}}=t_{\p}/[\p]_{\vv}$
with $[\p]_{\vv}=[\deg(\x_1)]_{\vv} \cdots [\deg(\x_r)]_{\vv}$. Finally, we let $\UU_{\Rb}$
stand for the $\Rb$-subalgebra of $\UU_{\Kb}$ generated by 
$\bigl\{\tilde{t}_{\x}\;|\; \x \in \ZZ^*\bigr\}$. Subalgebras $\UU^{\pm}_{\Rb}$ are defined in a similar fashion. There is an obvious action of $SL(2,\Z)$ on $\UU_{\Kb}$ and $ \UU_{\Rb}$.

\begin{prop}\label{P:gen1} The following hold~:
\begin{enumerate}
\item[i)] $\UU^{\pm}_{\Kb}=\bigoplus_{\mathbf{p} \in \mathbf{Conv}^{\pm}} \Kb \tilde{t}_{\mathbf{p}},$
\item[ii)] there is a triangular decomposition
$\UU_{\Kb}=\UU^{+}_{\Kb}
\otimes \UU^-_{\Kb}$;
in particular, we have  $\UU_{\Kb} = \bigoplus_{\mathbf{p} \in \mathbf{Conv}} \Kb \tilde{t}_{\mathbf{p}}.$
\end{enumerate}
\end{prop}
\begin{proof} We begin with i). Let us first show that the elements
$\bigl\{\tilde{t}_{\p}\;|\; \p \in \textbf{Conv}^+\bigr\}$ are linearly independent. For this we shall
 use a specialization argument. Let ${}'\UU^+_{\Rb}$ be the $\Rb$-algebra generated by some elements $\bigl\{{}'\tilde{t}_{\x}\;|\; \x \in \ZZ^+\bigl\}$ modulo the relations i) and ii) in Section~5.2 (these relations have coefficients in $\Rb$ when written in terms of the generators $\tilde{t}_{\x}$). By construction there is a canonical map $ {}'\UU^+_{\Rb} \to {}'\UU^+_{\Rb} \otimes_{\Rb} \Kb = \UU_{\Kb}, u \mapsto u \otimes 1$ whose image is $\UU^+_{\Rb}$. Moreover, for any elliptic curve $\E$ with Frobenius eigenvalues $\{\a, \overline{\a}\}$ there is a specialization morphism
$$\ev_{\E}~: {}'\UU^+_{\Rb} \to ({}'\UU^+_{\Rb})_{|\substack{\s=\a \\\bs=\overline{\a}}} = \UU^+_{\a,\overline{\a}} \simeq \U^+_{\E}.$$
Now assume that $\sum_{\p \in \textbf{Conv}^+} z_{\p} \tilde{t}_{\p}$ is a nontrivial (finite) linear relation in $\UU_{\Rb}$, with $z_{\p} \in \Rb$. Then $c:=\sum_{\p \in \textbf{Conv}^+} z_{\p}{}'\tilde{t}_{\p}$ is a torsion element of ${}'\UU^+_{\Rb}$. Let $Z $ denote its support, which
a strict subvariety of $\mathsf{Spec}(\Rb) \simeq \C^* \times \C^*$. We have
\begin{equation}\label{E:torsi}
\ev_{\E}(c) = \ev_{\E}\bigg(\sum_{\p} z_{\p} {}'\tilde{t}_{\p}\bigg)=\sum_{\p} z_{\p} \widetilde{T}_{\p} \in \U^+_{\E}
\end{equation}
where $\widetilde{T}_{\p}=T_{\p}/[\p]$. If $(\a, \overline{\a}) \not\in Z$ then $\ev_{\E}(c)=0$ and (\ref{E:torsi}) yields a nontrivial linear dependence relation between the elements 
$\bigl\{ \widetilde{T}_{\p}\;|\; \p \in \textbf{Conv}^+\bigl\}$, in contradiction with Theorem~\ref{L:bial}. It remains to find an elliptic curve with $(\a,\overline{\a}) \not\in Z$. For all prime powers $q$ let $N(q)$ be the number of possible Frobenius eigenvalues $\{\a, \overline{\a}\}$ for an elliptic curve over $\mathbb{F}_q$ (i.e. the number of isogeny classes of elliptic curves over $\mathbb{F}_q$). Then $\lim_{q \to \infty} N(q)=\infty$ (this is, for instance, a consequence of the main theorem in \cite{Honda}). But by Bezout's theorem the number of intersection points between $Z$ and $Y_q= \bigl\{(y,y')\;|\; yy'=q\bigr\}$ is bounded as $q \to \infty$. This provides the existence of the required elliptic curve, and concludes the proof of the linear independence of the elements $\bigl\{\tilde{t}_{\p}\;|\; \p \in \textbf{Conv}^+\bigr\}$ in $\UU^+_{\Rb}$ and hence in $\UU^+_{\Kb}$. The same arguments as in Lemma~\ref{L:Aplus} now show that $\UU^+_{\Kb}=\bigoplus_{\p \in \textbf{Conv}^+} \Kb \widetilde{t}_{\p}$.

\vspace{.1in}

We turn our attention to ii). We shall first show that the multiplication map $\UU^+_{\Kb} \otimes \UU^-_{\Kb} \to \UU_{\Kb}$ is surjective.
For this, using i), it is enough to see that
\begin{equation}\label{E:zoubi}
\tilde{t}_{\y}\tilde{t}_{\mathbf{p}} \in \UU^{+}_{\Kb}  \UU^-_{\Kb}
\end{equation}
for any $\y \in \ZZ^-$ and $\mathbf{p} \in \mathbf{Conv}^+$.
We say that a path $\mathbf{p}=(\x_1, \ldots, \x_r)$ is concave if $(\x_r,\ldots, \x_1)$ is convex. Let $\mathbf{Conc}^{\pm}$
denote the set of concave paths in $\ZZ^{\pm}$, and put $\mathbf{Conc}\simeq \mathbf{Conc}^+ \times \mathbf{Conc}^-$.
A symmetric version of Lemma~\ref{L:Aplus} and i) above shows that $\UU^{\pm}_{\Kb}=\bigoplus_{\mathbf{p}
\in \mathbf{Conc}^{\pm}} \Kb \tilde{t}_{\mathbf{p}}$. In particular, for any $\x,\y$ in $\ZZ$ with $\widehat{\x\y}>\pi$
we can  write $\tilde{t}_{\y}\tilde{t}_{\x}$ as a linear combination of elements $\tilde{t}_{\mathbf{p}}$
for \textit{concave} paths $\mathbf{p}$ lying in the triangle $\Delta_{\x,\y}$ (compare with Section~5.1.).
Now choose $\y \in \ZZ^-$ and $\x \in \ZZ^+$. If $\widehat{\x\y}>\pi$ then by the above remark $\tilde{t}_{\y}\tilde{t}_{\x}
\in  \bigoplus_{\mathbf{p} \in \mathbf{Conc}} \Kb\tilde{t}_{\mathbf{p}}$; if $\widehat{\x\y}=
\pi$ then $[\tilde{t}_{\y},\tilde{t}_{\x}] =0$; and if $\widehat{\x\y}<\pi$ then
$\tilde{t}_{\y}\tilde{t}_{\x} \in \bigoplus_{\mathbf{p} \in \mathbf{Conv}}
\Kb\tilde{t}_{\mathbf{p}}$. In all cases,
\begin{equation}\label{E:zouzou}
\tilde{t}_{\y}\tilde{t}_{\x} \in \UU^{+}_{\Kb} \UU^-_{\Kb}
\end{equation}
We shall  prove (\ref{E:zoubi}) by induction on the rank of $-\y$. If $\rank(-\y)=0$ (i.e.~if $\y=(l,0)$ 
for some $l \in \mathbb{N}^-$)
then $\bigl[\tilde{t}_{\y},\tilde{t}_{\x}\bigr] \in \UU^+_{\Kb} $ for any $\x \in \ZZ^+$.
Thus $\bigl[\tilde{t}_{\y}, \tilde{t}_{(\x_1, \ldots, \x_r)}\bigr]=
\sum_{i=1}^r \tilde{t}_{(\x_1, \ldots, \x_{i-1})}\bigl[\tilde{t}_{\y},\tilde{t}_{\x_i}\bigr] \tilde{t}_{(\x_{i+1},\ldots,\x_r)}
\in \UU^+_{\Kb} $. Now fix $\y \in \ZZ^-$ such that $-\y$ is of positive rank and
assume that (\ref{E:zoubi}) holds for all $\y'$ of smaller rank. Observe that if $\rank(\x)>0$
then from (\ref{E:zouzou}) we have
$[\tilde{t}_{\y},\tilde{t}_{\x}] =\sum_i u_i \tilde{t}_{\mathbf{p}_i^+}\tilde{t}_{\mathbf{p}^-_i}$ with
$\mathbf{p}_i^{\pm} \in \mathbf{Conv}^{\pm}$ and $\mathbf{p}^-_i=(\z^{(i)}_1, \ldots ,\z^{(i)}_{l_i})$
satisfying $\rank(-\z^{(i)}_j) < \rank(-\y)$. As a consequence, by the induction hypothesis we have
$\tilde{t}_{\mathbf{p}^-_i} \UU^+_{\Kb}
\subset  \UU^{+}_{\Kb}\otimes
\UU^-_{\Kb}$.
Next, if $\rank(\x)=0$ then $[\tilde{t}_{\y},\tilde{t}_{\x}] \in  \UU_{\Kb}^-$.
From these two facts we deduce that if $(\x_1, \ldots, \x_r) \in \mathbf{Conv}^+$ then
$[\tilde{t}_{\y},\tilde{t}_{\p}]=\sum_{i=1}^r \tilde{t}_{(\x_1, \ldots, \x_{i-1})}[\tilde{t}_{\y},
\tilde{t}_{\x_i}] \tilde{t}_{(\x_{i+1},\ldots,\x_r)} \in \UU^+_{\Kb}  \UU^-_{\Kb}$,
as wanted. This closes the induction and proves the surjectivity of the map
$\UU^+_{\Kb} \otimes \UU^-_{\Kb} \to \UU_{\Kb}$. It only remains to see that the elements
$\tilde{t}_{\mathbf{p}^+}\tilde{t}_{\mathbf{p}^-}$ for $\mathbf{p}^{\pm} \in \mathbf{Conv}^{\pm}$
 and $\a \in \ZZ$ are linearly independent over $\Kb$. For this, we may argue in the same fashion as in i) above using a specialization
 argument.
\end{proof}

\vspace{.1in}

We view $\UU_{\Rb}$ and $\UU_{\Kb}$ as generic versions of the Hall algebra $\U_{\E}$. 
 Moreover, one can lift various notions from $\U_{\E}$ to these generic forms. For instance we set
$$\widehat{\UU}^+_{\Kb}=\bigoplus_{\a \in \ZZ} \widehat{\UU}_{\Kb}^+[\a], \qquad \widehat{\UU}^+_{\Kb}[\a]=\prod_{\substack{\p \in \textbf{Conv}^+\\ wt(\p)=\a}} \Kb \tilde{t}_{\p}$$
and  define elements $\mathbf{1}^{\mathsf{ss}}_{\a} \in \UU^{+}_{\Kb}, \mathbf{1}_{\a} \in \widehat{\UU}^+_{\Kb}$ for any $\a \in \ZZ^+$ by the formulas
$$1+\sum_{l \geq 1} \mathbf{1}^{\mathsf{ss}}_{r\a_0} s^l = \exp\bigg(\sum_{l \geq 1} \tilde{t}_{l\a_0} s^l \bigg)$$
for any $\a_0$ such that $\deg(\a_0)=1$ and
$$\mathbf{1}_{\a}=\mathbf{1}_{\a}^{\mathsf{ss}} + \sum_{t > 1}
\sum_{\substack{\a_1 + \cdots + \a_t=\a \\ \mu(\a_1)< \cdots < \mu(\a_t)}} \hspace{-.1in}
\vv^{\sum_{i<j}\langle \a_i,\a_j\rangle}\mathbf{1}_{\a_1}^{\mathsf{ss}} \cdots
\mathbf{1}_{\a_t}^{\mathsf{ss}}.
$$
The elements $\bigl\{\mathbf{1}^{\mathsf{ss}}_{\a}\;|\; \a \in \ZZ^+\bigl\} $ belong to and actually generate over $\Rb$ the subalgebra $\UU^+_{\Rb}$, while the elements  
$\bigl\{\mathbf{1}_{\a}\;|\; \a \in \ZZ^+\bigl\}$ belong to and topologically generate over $\Rb$ the subalgebra $\widehat{\UU}^+_{\Rb}$. It is clear that the elements $\mathbf{1}^{\mathsf{ss}}_{\a}$
and $\mathbf{1}_{\a}$ specialize, for each given elliptic curve $\E$, to the corresponding elements of the Hall algebras $\U^+_{\E}$ and $ \widehat{\U}^+_{\E}$. Using the generators 
$\bigl\{\mathbf{1}_{\a}\;|\; \a \in \ZZ^+\bigr\}$ we may define a comutiplication  $\Delta$ on $\UU^+_{\Kb}$ by means of the formula (\ref{E:coprodun}). This comultiplication preserves $\UU^+_{\Rb}$.

\vspace{.1in}

Let us now give a more precise description of the integral form $\UU_{\Rb}$~:

\begin{prop}\label{P:gen2} The following proposition hold~:
\begin{enumerate}
\item[i)] $\UU^{\pm}_{\Rb}=\bigoplus_{\mathbf{p} \in \mathbf{Conv}^{\pm}} \Rb \tilde{t}_{\mathbf{p}},$
\item[ii)] there is a triangular decomposition
$\UU_{\Rb}=\UU^{+}_{\Rb}
\otimes \UU^-_{\Rb}$;
in particular, we have  $\UU_{\Rb} = \bigoplus_{\mathbf{p} \in \mathbf{Conv}} \Rb \tilde{t}_{\mathbf{p}},$
\item[iii)] for any $\alpha, \overline{\alpha} \in \C \backslash \{\pm 1\}$ we have $(\UU_{\Rb})_{\big|\substack{\s=\alpha \\ \bs =\overline{\alpha}}} = \UU_{\alpha, \overline{\alpha}}$,
\item[iv)] we have:  $(\UU_{\Rb})_{\big|\substack{\s=1 \\ \bs =1}} = \C[\tilde{t}_{\x}]_{\x \in \ZZ^*}$ is  a (commutative) polynomial algebra.
\end{enumerate}
\end{prop}
\noindent
\textit{Proof.} To prove statement i), we have to check that $\UU^+_{\Rb}$ is linearly spanned over $\Rb$ by $\bigl\{\tilde{t}_{\p}\;|\; \p \in \textbf{Conv}^+\bigr\}$. Let us temporarily denote by $\mathbf{V}$ the space $\bigoplus_{\p \in \text{Conv}^+} \Rb \tilde{t}_{\p}$. It is enough to show that $\mathbf{V}$ is a subalgebra of $\UU^+_{\Kb}$, i.e. that $\tilde{t}_{\p} \in \mathbf{V}$ for an \textit{arbitrary} path $\p$ in $\ZZ^+$. For this we proceed along the lines of Lemma~\ref{L:Aplus}, whose notations we shall  freely use. It is sufficient to show that
\begin{equation}\label{E:loop}
\tilde{t}_{\x} \tilde{t}_{\y} \in \mathbf{V}
\end{equation}
for any $\x, \y \in \ZZ^+$. We argue by induction on $\big|\det(\x,\y)\big| \in \N$. The claim (\ref{E:loop}) is clear if $\det(\x,\y)=0$. Let us fix an integer $l >0$ and assume that (\ref{E:loop}) holds for any pair $\x', \y'$ with $\big|\det(\x',\y')\big| <l$. Then, as in Lemma~\ref{L:Aplus} we have $\tilde{t}_{\q} \in \mathbf{V}$ for any path $\q$ with $a(\q) <l$. Let us fix a pair $\x,\y $ such that $\det(\x,\y)=l$. Up to $SL(2,\Z)$-action, we may assume that $\x=(0,n)$ and $\y=(r,d)$. Because the change of basis matrix between $\{\tilde{t}_{\z}\}$ and $\{\mathbf{1}^{\mathsf{ss}}_{\z}\}$ is invertible over $\Rb$, it is equivalent to prove that $[\tilde{t}_{\x},\mathbf{1}^{\mathsf{ss}}_{\y}] \in \mathbf{V}$. By Proposition~\ref{P:gen1} i) we have $[\tilde{t}_{\x},\mathbf{1}^{\mathsf{ss}}_{\y}] \in \bigoplus_{\p \in I_{\x,\y}} \Kb \tilde{t}_{\p}$. We have to show that all the coefficients belong to $\Rb$.
For this we write
\begin{equation}\label{E:klok1}
\begin{split}
\mathbf{1}_{\y}=\mathbf{1}_{(r,d)}=& \mathbf{1}^{\mathsf{ss}}_{(r,d)} + \sum_{k \geq 1} \vv^{rk} \mathbf{1}^{\mathsf{ss}}_{(r,d-k)} \mathbf{1}^{\mathsf{ss}}_{(0,k)} +\\
&\quad +  \sum_{\substack{(r_1,d_1) + \cdots +(r_l,d_l)=(r,d) \\
\frac{d_1}{r_1} < \cdots < \frac{d_l}{r_l};\;
r_1<r}}\hspace{-.1in} \vv^{\sum_{i<j}(r_id_j-r_jd_i)} \mathbf{1}^{\mathsf{ss}}_{(r_1,d_1)} \cdots \mathbf{1}^{\mathsf{ss}}_{(r_l,d_l)}.
\end{split}
\end{equation}
Thus
\begin{equation}\label{E:klok2}
\begin{split}
\tilde{t}_{(0,n)} \mathbf{1}^{\mathsf{ss}}_{\y}&=[\tilde{t}_{(0,n)},\mathbf{1}_{\y}]-\mathbf{1}_{\y}\tilde{t}_{(0,n)} -\sum_{k \geq 1} \vv^{rk} \tilde{t}_{(0,n)}\mathbf{1}^{\mathsf{ss}}_{(r,d-k)} \mathbf{1}^{\mathsf{\mathsf{ss}}}_{(0,k)} -\\
&- \sum_{\substack{(r_1,d_1) + \cdots +(r_l,d_l)=(r,d) \\
\frac{d_1}{r_1} < \cdots < \frac{d_l}{r_l};\;
r_1<r}} \hspace{-.1in} \vv^{\sum_{i<j}(r_id_j-r_jd_i)} \tilde{t}_{(0,n)}\mathbf{1}^{\mathsf{ss}}_{(r_1,d_1)} \cdots \mathbf{1}^{\mathsf{ss}}_{(r_l,d_l)}.
\end{split}
\end{equation}
Observe that the infinite sums in (\ref{E:klok2}) become finite after projection to $\bigoplus_{\p \in I_{\x,\y}} \tilde{t}_{\p}$. In the second sum
in (\ref{E:klok2}) we have $r_1<r$ therefore $\det\bigl((0,n),(r_1,d_1)\bigr) < rn=l$ and by our induction hypothesis we may straighten
$\tilde{t}_{(0,n)}\mathbf{1}^{\mathsf{ss}}_{(r_1,d_1)}=\sum_{\q \in I_{\x,(r_1,d_1)}} u_{\q} \tilde{t}_{\q}$. For any convex path $\q \in I_{\x,(r_1,d_1)}$ we have
$\mathbf{a}\bigl(\q \cup ((r_2,d_2), \ldots, (r_l,d_l)\bigr)<l$ and hence $\tilde{t}_{\q} \mathbf{1}^{\mathsf{ss}}_{(r_2,d_2)} \cdots \mathbf{1}^{\mathsf{ss}}_{(r_l,d_l)} \in \mathbf{V}$. By Lemma~\ref{L:comp} $\bigl[\tilde{t}_{(0,n)},\mathbf{1}_{(r,d)}\bigr] \in \mathbf{V}$. Finally, after projection to $\bigoplus_{\p \in I_{\x,\y}} \Kb \tilde{t}_{\p}$ we have $\mathbf{1}_{\y} \tilde{t}_{(0,n)} \in \mathbf{V}$. All together, working modulo $\mathbf{V}$ and projecting to $\bigoplus_{\p \in I_{\x,\y}} \Kb \tilde{t}_{\p}$ we get:
\begin{equation}\label{E:klok3}
\tilde{t}_{(0,n)} \mathbf{1}^{\mathsf{ss}}_{(r,d)} \equiv - \sum_{k \geq 1} \vv^{rn} \tilde{t}_{(0,n)} \mathbf{1}^{\mathsf{ss}}_{(r,d-k)} \mathbf{1}^{\mathsf{ss}}_{(0,k)}.
\end{equation}
Substituting (\ref{E:klok3}) into itself (i.e. developing each $\tilde{t}_{(0,n)}
\mathbf{1}^{\mathsf{ss}}_{(r,d-k)}$ according to (\ref{E:klok3})) sufficiently many times
 yields an expression
 \begin{equation}\label{E:klok4}
\tilde{t}_{(0,n)} \mathbf{1}^{\mathsf{ss}}_{(r,d)} \equiv  \sum_{k \geq N}  \tilde{t}_{(0,n)} \mathbf{1}^{\mathsf{ss}}_{(r,d-k)} w_k,
\end{equation}
where $w_k \in \Rb\bigl[\tilde{t}_{(0,1)}, \tilde{t}_{(0,2)}, \ldots\bigr]$. For $N \gg 0$ the right-hand side
 of (\ref{E:klok4}) vanishes after projection to $\bigoplus_{\p \in I_{\x,\y}} \Kb \tilde{t}_{\p}$. It follows that $\tilde{t}_{(0,n)} \mathbf{1}^{\mathsf{ss}}_{(r,d)} \in \mathbf{V}$ as desired. Statement i) is proven.

\vspace{.1in}

The proof of ii) is completely parallel to that of Proposition~\ref{P:gen1} ii). For iii), notice that by Theorem~\ref{T:main2} there exists an algebra map
$$\ev^*_{\E}~: \U_{\E} \simeq \UU_{\a,\overline{\a}} \to (\UU_{\Rb})_{\big|\substack{\s = \a\\ \bs=\overline{\a}}}, \qquad  \tilde{t}_{\x} \mapsto \tilde{t}_{\x}.$$
This map clearly sends the $\C$-basis $\bigl\{\tilde{t}_{\p}\;|\; \p \in \textbf{Conv}^+\bigr\}$ of 
$\U^+_{\E} \otimes_\KK \C$ to the $\C$-basis $\bigl\{\tilde{t}_{\p}\;|\; \p \in 
\textbf{Conv}^+\bigr\}$ of $(\UU_{\Rb})_{\big|\substack{\s = \a\\ \bs=\overline{\a}}}$. It remains to prove iv). For this we shall show that $[\UU_{\Rb}, \UU_{\Rb}] \in c_1(\s,\bs) \UU_{\Rb}$. Obviously, it is enough to prove that $[\tilde{t}_{\x},\tilde{t}_{\y}]\in c_1(\s,\bs)\UU_{\Rb}$ for any $\x,\y \in \ZZ^*$. Using the $SL(2,\Z)$-action we may assume that $\x =(0,n)$ for $n \geq 0$ and that $\y \in \ZZ^+$. By Lemma~\ref{L:comp} we have
$$\bigl[\tilde{t}_{(0,l)},\mathbf{1}_{\a}\bigr] = \frac{c_n(\s,\bs)[\rank(\a)]_{\vv^n}}{[n]_{\vv}} \mathbf{1}_{\a+(0,n)}$$
and it is easy to check that $c_n(\s,\bs)/[n]_{\vv} \in c_1(\s,\bs) \Rb$. Since the elements $\mathbf{1}_{\a}$ topologically generate $\widehat{\UU}^+_{\Rb}$, we may approximate $\tilde{t}_{\y}$ up to any degree of precision by a polynomial with $\Rb$-coefficients in the $\mathbf{1}_{\a}$. We conclude using the continuity of the multiplication (see Lemma~\ref{L:cont}).\qed

\vspace{.1in}
\noindent
As a consequence of iv) above and Weyl's theorem (see \cite{Weyl}) there is a natural isomorphism
\begin{equation*}
%\begin{split}
(\UU_\Rb)_{|\substack{\s=1\\ \bs=1}} \stackrel{\sim}{\lto} \C[x_1^{\pm 1}, \ldots, y_1^{\pm 1},\ldots]^{\mathfrak{S}_{\infty}}=\LL,  \quad 
\tilde{t}_{(r,d)} \mapsto \sum_i x_i^ry_i^d.
%\end{split}
\end{equation*}
Hence $\UU_{\Rb}$ may be thought of as a flat deformation of the ring of invariants $\LL$.

\vspace{.2in}

\paragraph{\textbf{6.3.}} There is an obvious $\mathfrak{S}_2$-symmetry in $\UU_{\Kb}$~: numbers 
$\s, \bs$ corresponding  to the two Frobenius eigenvalues in $H^1(\E_{\overline{\kk}}, \qlb)$, are interchangeable. Less obvious is the fact that this $\mathfrak{S}_2$-symmetry may be upgraded to an $\mathfrak{S}_3$-symmetry. To see this, we simply renormalize the generators. Set
$$u_{\x}=\frac{t_{\x}}{c_{\deg(\x)}(\s,\bs)}, \qquad (\x \in \ZZ^*)$$
and for any $i \geq 1$ put
\begin{equation}\label{E:clevercoeff}
\alpha_i=\alpha_i(\s,\bs)=(1-\s^i)(1-\bs^i)(1-(\s\bs)^{-i}) /i.
\end{equation}
The defining relations in Section~5.2 may now be rewritten as
\begin{enumerate}
\item[i)] For a pair of collinear  $\x,\x'$  we have
$$[u_\x,u_{\x'}]=0.$$
\item[ii)] Assume that $\x,\y \in \ZZ^*$ are such that $\deg(\x)=1$ and that
$\boldsymbol{\Delta}_{\x,\y}$ has no interior lattice point. Then
$$[u_\y,u_{\x}]=\epsilon_{\x,\y}
\frac{\theta_{\x+\y}}{\alpha_1}$$
where  the elements $\theta_{\z}$, $\z \in \ZZ^*$ are obtained by equating the Fourier
 coefficients of the collection of relations
\begin{equation}\label{E:formulatheta2}
\sum_i \theta_{i\x_0}s^i = \exp\bigg(\sum_{r \geq 1}\alpha_r u_{r\x_0}s^r\bigg),
\end{equation}
for any $\x_0 \in \ZZ^*$ such that $\deg(\x_0)=1$.
\end{enumerate}

\vspace{.1in}

In this presentation it is obvious that $\UU_{\Kb}$ is equipped with an $\mathfrak{S}_3$ family of $\C$-automorphisms $\Theta_{\gamma}$ for $\gamma \in \mathrm{Perm}\{\s,\bs, (\s\bs)^{-1}\}$ simply defined by $\Theta_{\gamma}(u_{\x})=u_{\x}$, $\Theta_{\gamma}(\bullet)=\gamma(\bullet)$ for $\bullet \in \{\s,\bs,(\s\bs)^{-1}\}$. This symmetry may seem puzzling at first glance~: for any fixed elliptic curve $\E$ over a finite field $\mathbb{F}_q$ we have $|\s|=|\bs|=q^{1/2}$ while 
$\big|(\s\bs)^{-1}\big|=q^{-1}$.

\vspace{.2in}

\paragraph{\textbf{6.4.}} In order to define  the coproduct $\Delta$ of a Hall algebra $\mathbf{H}$ 
or to construct the  Drinfeld double of $\mathbf{H}$, it is usually necessary to add an extra commutative `Cartan' subalgebra $\boldsymbol{\mathcal{K}}$ to $\mathbf{H}$ (see e.g. \cite{SLectures}). In the present case of the category of coherent sheaves over an elliptic curve we could avoid doing so because the symmetrized Euler form vanishes. However adding the corresponding ``Cartan'' subalgebra $\boldsymbol{\mathcal{K}}$ provides a natural central extension $\widetilde{\H}_{\E}$ of $\mathbf{H}$ (and similarly for $\U_{\E}$ and $ \UU_{\Kb}$). This central extension is also 
important in applications (see e.g. \cite{SV2})

For the sake of brevity, we only write down the relations in $\widetilde{\UU}_{\Kb}$, using the rescaled presentation of Section~6.3.

\vspace{.15in}

\addtocounter{theo}{1}
\paragraph{\textbf{Definition \thetheo}} Let $\widetilde{\UU}_{\Kb}$ be the $\Kb$-algebra defined by generators $\{\boldsymbol{\kappa}_{\a}\;|\; \a \in \ZZ\}$ and $\{u_{\x}\;|\; \x \in \ZZ^*\}$ modulo the following set of relations~:
\begin{enumerate}
\item[i)] the subalgebra $\boldsymbol{\mathcal{K}}$ generated by $\{\boldsymbol{\kappa}_{\a}\;|\;\a \in \ZZ\}$ is central and we have
$$\boldsymbol{\kappa}_0=1, \qquad \boldsymbol{\kappa}_{\a}\boldsymbol{\kappa}_{\beta}=\boldsymbol{\kappa}_{\a+\beta},$$
\item[ii)] if $\x,\y$ belong to the same line in $\ZZ$ then
\begin{equation*}
[u_{\y},u_{\x}]=\delta_{\x,-\y}
\frac{\boldsymbol{\kappa}_{\x}-\boldsymbol{\kappa}^{-1}_{\x}}{\alpha_{\deg(\x)}}
\end{equation*}
\item[iii)] if $\x,\y \in \ZZ^*$ are such that $\deg(\x)=1$ and that
$\boldsymbol{\Delta}_{\x,\y}$ has no interior lattice point then
\begin{equation*}
[u_\y,u_{\x}]=\epsilon_{\x,\y}\boldsymbol{\kappa}_{\a(\x,\y)}
\frac{\theta_{\x+\y}}{\alpha_{1}},
\end{equation*}
where
\begin{equation*}
\a(\x,\y)=\begin{cases}
\epsilon_\x(\epsilon_{\x}\x+\epsilon_{\y}\y-\epsilon_{\x+\y}(\x+\y))/2 & \text{if}\; \epsilon_{\x,\y}=1,\\
\epsilon_\y(\epsilon_{\x}\x+\epsilon_{\y}\y-\epsilon_{\x+\y}(\x+\y))/2
& \text{if}\; \epsilon_{\x,\y}=-1,
\end{cases}
\end{equation*}
and where the elements $\theta_{\z}$, $\z \in \ZZ^*$, are given
by
\begin{equation*}
\sum_i \theta_{i\x_0}s^i=\exp\big(\sum_{r \geq
1}\alpha_r u_{r\x_0}s^r\big),
\end{equation*}
for any $\x_0 \in \ZZ^*$ such that $\deg(\x_0)=1$, where the coefficients $\alpha_r$ are given by (\ref{E:clevercoeff}).
\end{enumerate}

\vspace{.1in}

Note that by relation ii), the algebra $\widetilde{\UU}_{\Kb}$ contains many copies of the Heisenberg algebra (one for each line in $\ZZ$). Hence $\widetilde{\UU}_{\Kb}$ can be thought of as a flat deformation of a Heisenberg algebra over $\ZZ$.

\vspace{.1in}
\noindent
The triangular decomposition of $\widetilde{\UU}_{\Kb}$ now takes the form
\begin{equation}\label{E:trianng}
\widetilde{\UU}_{\Kb} \simeq \UU^+_{\Kb} \otimes \boldsymbol{\mathcal{K}} \otimes \UU^-_{\Kb}.
\end{equation}
One consequence of the central extension is that the group $SL(2,\Z)$ no longer acts on $\widetilde{\UU}_{\Kb}$~~: only its universal cover $\widetilde{SL}(2,\Z)$ does. There is a short exact sequence
$$
1 \longrightarrow   \Z \longrightarrow  \widetilde{SL}(2,\Z) \longrightarrow
  SL(2,\Z) \longrightarrow   1.
$$
For any slope $\frac{q}{p} \in \mathbb{Q} \cup \{\infty\}$ and any $\gamma \in \widetilde{SL}(2,\Z)$ we define a winding number $n(\gamma, \frac{q}{p})$ as follows. There is a natural action of $SL(2,\Z)$ on the circle $S^1=(\mathbb{R}^2\setminus \{0\})/\mathbb{R}^{+\star}$.
Using the identification $S^1= \mathbb{R}/2\Z$, we can  uniquely lift this action to an $\widetilde{SL}(2,\Z)$-action on $\mathbb{R}$.
Any $(q,p) \in \ZZ^*$ gives rise to an element $(q:p) \in S^1$ and if $\overline{(q:p)} \in \mathbb{R}$ is any
lift of $(q:p)$ then
\begin{equation}\label{E:Defwind}
n\bigg(\tilde{\gamma},\frac{q}{p}\bigg)=
\begin{cases}
\#\bigl(\Z \cap [\overline{(q,p)}, \tilde{\gamma}(\overline{(q:p)})]\bigr) & \text{if}\; \tilde{\gamma}(\overline{(q:p)}) \geq \overline{(q:p)}\\
-\#\bigl(\Z \cap [\tilde{\gamma}(\overline{(q,p)}), \overline{(q:p)}]\bigr) & \text{otherwise.}
\end{cases}
\end{equation}

One checks that the following rule gives rise to an $\widetilde{SL}(2,\Z)$-action on $\widetilde{\UU}_{\Kb}$ by automorphisms~:
\begin{equation}\label{E:aaction}
\Phi(\boldsymbol{\kappa}_{\x})= \boldsymbol{\kappa}_{\Phi(\x)}, \qquad \Phi ( u_{\x})= u_{\Phi(\x)} \boldsymbol{\kappa}_{\Phi(\x)}^{n(\Phi,\mu(\x))}.
\end{equation}

\vspace{.1in}

We leave it to the reader to define the integral forms, specializations, etc of $\widetilde{\UU}_{\Kb}$.  All properties (such as (\ref{E:trianng}) and (\ref{E:aaction}), finite generation etc.) extend to these settings.

\section{Summary}
Let us sum up the main results obtained in this article. To any elliptic curve
$\E$ defined over  a finite field $\kk = \mathbb{F}_q$ we have
 attached an associative algebra ${\mathbf{U}}_\E$
over the field $K= \mathbb{Q}(v)$, where  $v^{-2} = q$. Let $\sigma \in
\overline{\mathbb{Q}}$
be such that $\bar{\sigma}\sigma = q$ and $\bigl|\E(\mathbb{F}_{q^i})\bigr| =
q^i + 1 - (\sigma^i + \bar{\sigma}^i)$.  Then we have:

\medskip
\medskip

\noindent
1. The algebra ${\mathbf{U}}_\E$ is $\mathbb{Z}^2$-graded and
$K =  {\mathbf{U}}_\E\bigl[(0,0)\bigr] $ is the center of ${\mathbf{U}}_\E$.

\medskip

\noindent
2. The algebra
${\mathbf{U}}_\E$ can be described by the following generators and relations:
\begin{enumerate}
\item For $(r,d) \in \mathbb{Z}^2\setminus \{(0,0)\}$ we have a generator
$T_{(r,d)} \in  {\mathbf{U}}_\E\bigl[(r,d)\bigr]$.
\item  Let $\mathsf{gcd}(r,d) = 1$, then we defined elements $\Theta_{i(r,d)} \in
{\mathbf{U}}_\E\bigl[i(r,d)\bigr], \, i\ge 1$ using  the following equality
$$
1+  \sum\limits_{i=1}^\infty \Theta_{i(r,d)} s^i =
\exp\bigl((v^{-1} - v)\sum_{j= 1}^\infty T_{j(r,d)}s^j\bigr),
$$
where $s$ is a formal parameter.
\item\label{Rel1} If  the vectors $(r,d)$ and $(r', d')$ are collinear then  we have:
$$
\bigl[T_{(r,d)}, T_{(r',d')}\bigr] =0.
$$
\item\label{Rel2}
Assume that  $(r,d), (r',d') \in \mathbb{Z}^2\setminus \{(0,0)\}$ are such that
 $\mathsf{gcd}(r,d) = 1$ and
the triangle with the corners $(0,0), (r,d), (r',d')$
contains no interior points. Then
$$
\bigl[T_{(r,d)}, T_{(r',d')}\bigr] = \mathrm{sign}\bigl(rd' - r'd\bigr)
c_h\frac{\Theta_{(r+r', d+d')}}{v - v^{-1}},
$$
where $h = \mathsf{gcd}(r',d')$ and $c_h = 
\frac{\displaystyle v^h [h]_v}{\displaystyle h} 
\big|\E(\mathbb{F}_{q^h})\big|$.
Relations (\ref{Rel1}) and (\ref{Rel2})  form
 a complete list of relations of ${\mathbf{U}}_\E$, see Theorem \ref{T:main2}.
 The structure constants of
${\mathbf{U}}_\E$ are Laurent polynomials in $\sigma^{\pm \frac{1}{2}}$ and $\bar\sigma^{\pm \frac{1}{2}}$, so we may also introduce a generic version
${\UU}_{\Rb}$
of the Hall algebra ${\mathbf{U}}_\E$, defined  over the ring $\Rb = \C[\sigma^{\pm \frac{1}{2}}, \bar\sigma^{\pm \frac{1}{2}}]$, see  Section 6.2.
\end{enumerate}

\medskip
\noindent
3. The algebra ${\mathbf{U}}_\E$ is finitely generated and the elements
$T_{(\pm 1, 0)}, T_{(0, \pm 1)}$ generate ${\mathbf{U}}_\E$, see Corollary \ref{C:Uisfing}.

\medskip
\noindent
4. The algebra ${\mathbf{U}}_\E$ carries a natural
${SL}(2,\mathbb{Z})$--action:
for any
$\gamma \in SL(2,\mathbb{Z})$ the map
$T_{(r,d)} \mapsto T_{\gamma(r,d)}$ induces
an algebra automorphism of ${\mathbf{U}}_\E$.

\medskip
\noindent
5.
Let $\mathbf{U}_\E^\pm  = \bigl\langle T_{(r,d)}| (r,d) \in
(\mathbb{Z}^2)^\pm\bigr\rangle$,
then ${\mathbf{U}}^\pm_\E$ are graded topological bialgebras, see Lemma \ref{L:huh6}.
  This  means that  there is
a graded coassociative ring homomorphism
$$\Delta: {\mathbf{U}}_\E^\pm    \lto  {\mathbf{U}}_\E^\pm
\widehat\otimes
 {\mathbf{U}}_\E^\pm,$$
taking value in a certain completion of
$ {\mathbf{U}}_\E^\pm  \otimes
 {\mathbf{U}}_\E^\pm $ and given by the collection of linear maps for each
$\a,\b \in (\mathbb{Z}^2)^\pm$
$$
\Delta_{\a,\b}:  {\mathbf{U}}^\pm_\E\bigl[\a + \b\bigr]  \lto
 {\mathbf{U}}^\pm_\E[\a] \otimes   {\mathbf{U}}^\pm_\E[\b].
$$

\medskip
\noindent
6. The algebra ${\mathbf{U}}_\E$
is isomorphic to the  Drinfeld double
of the topological bialgebra ${\mathbf{U}}^+_\E$ and one has the decomposition
$
{\mathbf{U}}_\E = \mathbf{U}^+_\E  \otimes \mathbf{U}^-_\E,
$
where $\mathbf{U}^\pm_\E = \bigl\langle T_{(r,d)}| (r,d) \in (\mathbb{Z}^2)^+\bigr\rangle$, see Theorem
\ref{L:bial}.

\medskip
\noindent
7.  The algebra ${\mathbf{U}}_\E$ has a monomial basis
$\bigl\{T_{(r_1, d_1)} T_{(r_2, d_2)} \dots  T_{(r_n, d_n)}\bigr\}$
parameterized by the set of convex paths
$\bigl((r_1, d_1), (r_2, d_2), \dots, (r_n, d_n)\bigr)$ in  $\mathbb{Z}^2$, see Theorem \ref{L:bial}.

\medskip
\noindent
8. The algebra ${\mathbf{U}}_\E$ is a flat deformation of
the ring $$K\bigl[x_1^\pm ,x_2^\pm ,\dots, y_1^\pm ,y_2^\pm
,\dots\bigr]^{\mathfrak{S}_\infty}$$ of symmetric
Laurent series.

\medskip
\noindent
9. One can also write down some explicit formulas for the coproduct of certain generators
of ${\mathbf{U}}_\E^+$ (Proposition \ref{P:Mac}
and Lemma \ref{L:coprodun}):
$$
\Delta\bigl(T_{(0,d)}\bigr) = T_{(0,d)}\otimes 1 + 1 \otimes
T_{(0,d)} \; \; \mbox{and}  \; \; \Delta\bigl(T_{(1,d)}\bigr) =
T_{(1,d)}\otimes 1 +
\sum_{l\ge 0} \Theta_{(0,l)} \otimes T_{(1,d-l)}.
$$

\centerline{\large{Appendix A}}
\addcontentsline{toc}{section}{Appendix A}

\vspace{.1in}

In this appendix, we provide the details regarding the properties of the
Fourier-Mukai transforms on elliptic curves
defined over a \emph{finite} field $\kk$.

For a projective curve $\Y$ defined over the field $\kk$  consider the functor
${\rm Pic}^0_{\Y/\kk}: {\rm Sch}_{\kk} \lto {\rm Sets}$ given by
\begin{equation*}
\begin{split}
{\rm Pic}^0_{\Y/\kk}(\S) = \Bigl\{\kF &\in {\rm Coh}_{\Y\times\S}|\; \kF
\mbox{\,\,\rm is \,\,} \S-\mbox{\rm flat and for any closed point}\\
&s: \Spec(\ll) \lto S \,\, \; \mbox{\rm holds} \,\,\;
s_{\ll}^*(\kF) \in {\rm Pic}^0(\Y_{\ll})\Bigr\}\Big/\sim
\end{split}
\end{equation*}
where $\Y_{\ll} = \Y \times_{\Spec(\kk)} \Spec(\ll)$ and
 the map $s_{\ll}: \Y_{\ll} \lto \Y \times \S$ is induced by the base change
and the equivalence relation is
$\kF \sim \kF \otimes \pi^*_{\S}(\kL)$ for any  locally free rank one sheaf $\kL$ on
$\S$.

\medskip

In the case of an elliptic curve  $\E$ over $\kk$ with a rational
point $p_0$ the functor
${\rm Pic}^0_{\E/\kk}$ is representable by the pair $(\E,\kP)$, where
$\kP = \kO_{\E\times \E}(-\Delta+p_0\times\E +\E\times p_0)$ and
$\Delta \subset \E\times \E$ is the diagonal, see for example  \cite[example 8.9.iii]{AK1}.

\medskip

The sheaf $\kP^\vee$ is locally free on $\E\times \E$ and hence flat over
$\E$. Moreover,
for any closed  point $p: \Spec(\ll) \lto \E$ one has an isomorphism
$\kP^\vee\big|_{\E_{\ll}} \cong \bigl({\kP\big|_{\E_{\ll}}}\bigr)^\vee$.
By the universal property of $(\E,\kP)$ there exists
a unique map $i: \E \lto \E$ and a line bundle $\kL$ on  $\E$ such that
$\kP^\vee \otimes \pi_2^*\kL \cong (1 \times i)^*\kP$.
Denote by $\sigma = p_0 \times 1: \E \lto
\E \times \E$. From equalities $\sigma^*\kP^\vee  \cong \kO$
and $\sigma^*(1 \times i)^*\kP  \cong (1 \times i)^* \sigma^* \kP \cong  \kO$
we conclude that
$$\kP^\vee  \cong (1 \times i)^*\kP.$$
Moreover, the isomorphism  $\kP \cong \kP^{\vee\vee}$ and the universality of
$(\E, \kP)$ imply  $i^2 = 1$.

\vspace{.1in}

\paragraph{\textbf{Proposition A.1.}}
%\label{closure}
\textit{Let $\Y$ be a projective variety over $\kk$ and  $\bar{\kk}$ the algebraic closure
of $\kk$. For any  field extension $\kk\subset \ll$ denote by
$\Y_{\ll} = \Y \times_{\Spec(\kk)} \Spec(\ll)$ and  by  $\varphi_{\ll}:
 \Y_{\ll}\lto  \Y$ the base-change map.
Let $\kF$, $\kG
$ be two coherent sheaves, denote
$\kF_{\ll} =  \varphi_{\ll}^*(\kF)$ and $\kG_{\ll} =  \varphi_{\ll}^*(\kG)$.
Assume that  $\kF_{\bar{\kk}} \cong  \kG_{\bar{\kk}}$ then $\kF \cong \kG$.} \\

\noindent
\textit{Proof}.
Let $f: \kF_{\bar{\kk}} \lto  \kG_{\bar{\kk}}$ and $g: \kG_{\bar{\kk}} \lto
 \kF_{\bar{\kk}}$ be two maps such that
$gf = 1_{\kF_{\bar{\kk}}}$ and $fg = 1_{\kG_{\bar{\kk}}}$. From the isomorphism
$\mathrm{Hom}_\kO(\kF, \kG) \otimes_{\kk} \bar{\kk} \cong
\mathrm{Hom}_{\kO_{\bar{\kk}}}(\kF_{\bar{\kk}}, \kG_{\bar{\kk}})$
follows that
$f = \sum\limits_{i=1}^n \bar{a}_i \varphi_i$ and
$g = \sum\limits_{j=1}^m  \bar{b}_j \psi_j$,  where
$\bar{a}_1, \dots, \bar{a}_n,  \bar{b}_1,\dots, \bar{b}_m \in \bar{\kk}$,
 $\varphi_1, \dots, \varphi_n$ is a basis
of $\mathrm{Hom}_\kO(\kF, \kG)$ over $\kk$ and $\psi_1, \dots, \psi_m$ a basis of $\mathrm{Hom}_\kO(\kG, \kF)$ over $\kk$.
Let $\ll$ be the finite extension of $\kk$ generated by the elements
$\bar{a}_1, \dots, \bar{a}_n; \bar{b}_1,\dots, \bar{b}_m$,  then for the base-change  map
$\varphi_{\ll}: \Y_{\ll} \lto \E$ we have  $\varphi_{\ll}^{*}(\kF) \cong \varphi_{\ll}^{*}(\kG)$.
By the projection formula
we get $\varphi_{\ll*}(\varphi_{\ll}^{*}\kF\bigr) = \kF \otimes \varphi_{\ll*}(\kO_{\ll})$.
Let $d = {\rm deg}(\ll/\kk)$.  Since  $\Spec(\ll) =
\underbrace{\Spec(\kk) \sqcup \Spec(\kk) \sqcup \dots \sqcup \Spec(\kk)}_{\mbox{$d$ times}}$
as a scheme over
$\kk$, we have    $\Y_{\ll} = \underbrace{\Y \sqcup \Y \sqcup \dots \sqcup \Y}_{\mbox{$d$ times}}$ and
$\varphi_{\ll*}(\kO_{\ll}) = \kO^d$.
Hence   $\kF^d \cong \kG^d$ and the Krull-Schmidt theorem implies $\kF
\cong \kG$.
\qed

\medskip

\noindent
\textit{Proof of Proposition~\ref{Mukai}.}
In the case of an algebraically closed field $\bar{\kk}$ this result
 was shown by Mukai \cite{Muk}.
This isomorphism   is equivalent to the fact that
$$\mathbf{R}\pi_{13}(\pi_{12}^*\kP \otimes \pi_{23}^*\kP) \cong \kO_{i(\Delta)} [-1],$$
where $\kO_{i(\Delta)}$ is the structure sheaf of the subscheme
$i(\Delta) \subseteq \E \times \E$.
The case of a finite field $\kk$ can be derived from the corresponding result about
$\bar{\kk}$ by going into the algebraic
closure:
$\varphi_{\bar{\kk}}: \E_{\bar{\kk}} \lto \E$ and using the isomorphism
$$(\varphi_{\bar{\kk}}\times \varphi_{\bar{\kk}})^*\kO_{\E\times\E}(-\Delta + p_0\times\E +
\E\times p_0) \cong \kO_{\E_{\bar{\kk}}\times \E_{\bar{\kk}}}(-\Delta_{\bar{\kk}} +
\bar{p}_0\times \E_{\bar{\kk}} +
\E_{\bar{\kk}} \times \bar{p}_0),$$  the  flat base-change and the Proposition~A.1 above.
\qed

\vspace{.15in}

\paragraph{\textbf{Proposition A.2.}}[see  \cite[Proposition 3.8]{Muk}]
\textit{Let  $D = \mathbf{R}{\mathcal Hom}(-,\kO)$ be the dualizing functor.
Then there is an isomorphism
of functors}
$$
D \circ \Phi \cong i^* \circ [1] \circ \Phi \circ D.
$$
\textit{Proof}. This result is a corollary of
the isomorphism $\kP^\vee \cong (1 \times i)^* \kP$ and can be proven along the same
lines as in
\cite{Muk}.\qed

\vspace{.2in}

\centerline{\large{Appendix B}}
\addcontentsline{toc}{section}{Appendix B}

\vspace{.1in}

In the second appendix, we provide proofs for some technical statements
regarding the Drinfeld double construction for topological bialgebras and some properties
of Hopf algebras, which are crucial for the proof of Theorem \ref{L:bial}.

\medskip

\noindent
\textit{Proof of Lemma~\ref{L:relprod}.} For simplicity, we drop the
exponents $\pm$ in the notation. Since both  statements in the Lemma are similar, we
give a proof only of  the first one. By assumption, we have  for any $k$
\begin{equation}
\sum_{i,j} a_j^{(1)} (c_k^{(1)})^{(2)}_i \bigl((c_k^{(1)})^{(1)}_i,a_j^{(2)}\bigr)=
\sum_{i,j} (c_k^{(1)})^{(1)}_i a_j^{(2)}\bigl((c_k^{(1)})_i^{(2)},a_j^{(1)}\bigr), \tag{B.1}
\end{equation}
\begin{equation}
\sum_{i,j} b_j^{(1)} (c_k^{(2)})^{(2)}_i \bigl((c_k^{(2)})^{(1)}_i,b_j^{(2)}\bigr)=\sum_{i,j} (c_k^{(2)})^{(1)}_i b_j^{(2)}\bigl((c_k^{(2)})_i^{(2)},b_j^{(1)}\bigr). \tag{B.2}
\end{equation}
Note that all sums above are in fact finite. Now, we compute
\begin{equation}
\begin{split}
\sum_{k,j} (ab)_j^{(1)}c_k^{(2)}\bigl(c_k^{(1)},(ab)^{(2)}_j\bigr)=&\sum_{i,j,k} a_i^{(1)}b_j^{(1)}c_k^{(2)}\bigl(c_k^{(1)},a_i^{(2)}b_j^{(2)}\bigr)\\
=&\sum_{i,j,k,l} a_i^{(1)}b_j^{(1)}c_k^{(2)}\bigl((c_k^{(1)})_l^{(1)},a_i^{(2)}\bigr)\bigl((c_k^{(1)})^{(2)}_l,b_j^{(2)}\bigr),
\end{split}\tag{B.3}
\end{equation}
where we used the Hopf property of the pairing $(\; \,,\;)$ and Proposition \ref{P:topbial}.
Next, by coassociativity, we have
$\sum_{k,l} (c_k^{(1)})_l^{(1)} \otimes (c_k^{(1)})_l^{(2)} \otimes c_k^{(2)}=\sum_{k,l}c_k^{(1)} \otimes (c_k^{(2)})_l^{(1)} \otimes (c_k^{(2)})_l^{(2)}$, and substituting in (B3), we obtain
\begin{equation}
\begin{split}
\sum_{i,j,k,l} a_i^{(1)}&b_j^{(1)}c_k^{(2)}\bigl((c_k^{(1)})_l^{(1)},a_i^{(2)}\bigr)\bigl((c_k^{(1)})^{(2)}_l, b_j^{(2)}\bigr)\\
=&\sum_{i,j,k,l}a_i^{(1)}b_j^{(1)} (c_k^{(2)})^{(2)}_l (c_k^{(1)},a_i^{(2)})\bigl((c_k^{(2)})^{(1)}_l,b_j^{(2)}\bigr)\\
=&\sum_{i,j,k,l} a_i^{(1)}(c_k^{(2)})_l^{(1)}b_j^{(2)}(c_k^{(1)},a_i^{(2)})\bigl((c_k^{(2)})^{(2)}_l,b_j^{(1)}\bigr),
\end{split}\tag{B.4}
\end{equation}
where we made use of (B.2). In the same way, coassociativity and (B.1) allow us to transform the last expression into
\begin{equation*}
\begin{split}
\sum_{i,j,k,l} a_i^{(1)}&(c_k^{(1)})_l^{(2)}b_j^{(2)}\bigl((c_k^{(1)})_l^{(1)},a_i^{(2)}\bigr)(c_k^{(2)},b_j^{(1)})\\
&=\sum_{i,j,k,l}(c_k^{(1)})_l^{(1)}a_i^{(2)}b_j^{(2)}\bigl((c_k^{(1)})_l^{(2)},a_i^{(1)}\bigr)(c_k^{(2)},b_j^{(1)})\\
&=\sum_{i,j,k,l}c_k^{(1)}a_i^{(2)}b_j^{(2)}\bigl((c_k^{(2)})_l^{(1)},a_i^{(1)}\bigr)
\bigl((c_k^{(2)})^{(2)}_l,b_j^{(1)}\bigr).
\end{split}
\end{equation*}
Finally, using the Proposition \ref{P:topbial} and the Hopf property of $(\;\,,\;)$ again, we can rewrite the  last term as
$$\sum_{i,j,k}c_k^{(1)}a_i^{(2)}b_j^{(2)}\bigl(c_k^{(2)},a_i^{(1)}b_j^{(1)}\bigr) =\sum_{i,k}c_k^{(1)}(ab)^{(2)}_i\bigl(c_k^{(2)},(ab)^{(1)}_i\bigr).$$
All together, we see that $R(ab,c)$ is a consequence of relations (B.1) and (B.2). The Lemma is proved.\qed

\medskip
\noindent
The remaining part of Appendix B  is devoted to the

\medskip
\noindent
\textit{Proof of Proposition~\ref{P:DDoub}.} Recall that by the definition of $\Delta$ we have
$$
\Delta\bigl([\kF]\bigr) = \sum_{\kK \rightarrowtail \kF} v^{-\langle \kF/\kK, \kK\rangle}
\frac{{P}^\kF_{\kF/\kK, \kK}}{a_\kF}[\kF/\kK] \otimes \kK.
$$
Iterating this formula we have
$$
\Delta^2\bigl([\kF]\bigr) = (1 \otimes \Delta)\Delta\bigl([\kF]\bigr) =
\sum_{\kL  \rightarrowtail \kK \rightarrowtail \kF}
c^\kF_{\kF/\kK, \kK/\kL, \kL}
[\kF/\kK] \otimes [\kK/\kL]
\otimes [\kL],
$$
where $c^\kF_{\kF/\kK, \kK/\kL, \kL} =
{P}^{\kF}_{\kF/\kK, \kK} {P}^{\kK}_{\kK/\kL, \kL} \frac{\displaystyle v^{-\langle \kK/\kK, \kK\rangle - \langle\kK/\kL, \kL \rangle}}{\displaystyle a_{\kF}}$.
So, in general we can write
\begin{equation}
\Delta^n\bigl([\kF]\bigr) =
\sum_{\kF_n \rightarrowtail \kF_{n-1} \rightarrowtail \dots \rightarrowtail \kF_1 \rightarrowtail \kF}
c^\kF_{\kA_1,\kA_2,\dots,\kA_{n+1}}
[A_1] \otimes [\kA_2] \otimes \dots \otimes
[\kA_{n+1}],\tag{B.5}
\end{equation}
where $\kA_i = \kF_{i-1}/\kF_i$, $\kA_{n+1} = \kF_n$ and
\begin{equation*}
c^\kF_{\kA_1,\kA_2,\dots,\kA_{n+1}} =
{P}^{\kF}_{\kA_1, \kF_1} {P}^{\kF_1}_{\kA_2, \kF_2} \dots
{P}^{\kF_{n-1}}_{\kA_n, \kF_n} \frac{v^{-\sum_{i=1}^n \langle \kA_i, \kF_i\rangle}}{a_{\kF}}.
\end{equation*}
We can also write in a dual  way:
$$
\Delta^n\bigl([\kF]\bigr) =
\sum_{\kF \twoheadrightarrow \kF_1 \twoheadrightarrow \kF_2
\twoheadrightarrow \dots \twoheadrightarrow \kF_n}
d^\kF_{\kB_1,\kB_2,\dots,\kB_{n+1}}
 [\kB_{n+1}] \otimes
[\kB_n] \otimes \dots \otimes
[\kB_{1}],
$$
where $\kB_{n+1} = \kF_n$, $\kB_i = \mathrm{ker}(\kF_{i-1} \lto \kF_i)$ and
\begin{equation}
d^\kF_{\kB_1,\kB_2,\dots,\kB_{n+1}} =
{P}^{\kF}_{\kF_1, \kB_1} {P}^{\kF_1}_{\kF_2, \kB_2} \dots {P}^{\kF_{n-1}}_{\kF_n, \kB_n}
\frac{v^{-\sum_{i=1}^n \langle \kF_i, \kB_i\rangle}}{a_{\kF}}.\tag{B.6}
\end{equation}

\medskip

\noindent
\textbf{Definition B.1}. For $\alpha \in (\mathbb{Z}^2)^+$ define an operator
$T: \H_\E[\alpha] \lto \widehat{\H}_\E[\alpha]$ by the following formulas:
\begin{enumerate}
\item $T\bigl([0]\bigr) = T(1) = 1$.
\item Let $\alpha \ne (0,0)$ and $\kF$ be a coherent sheaf of class $\alpha$. Then
\begin{equation}
T\bigl([\kF]\bigr) =
 - [\kF] + \sum_{n=1}^\infty (-1)^n \sum_{
\kF_n \stackrel{\ne}\rightarrowtail   \dots \stackrel{\ne}\rightarrowtail\kF_1 \stackrel{\ne}\rightarrowtail  \kF}
c^\kF_{\kA_1,\kA_2,\dots,\kA_{n+1}}
  [\kA_{n+1}] \otimes \dots  \otimes
[\kA_{1}],\tag{B.7}
\end{equation}
where $\kA_1, \kA_2,\dots, \kA_{n+1}$ and $c^\kF_{\kA_1,\kA_2,\dots,\kA_{n+1}}$  are the same as in
(B.5).
\end{enumerate}

\medskip
\noindent
In order to see that the operator $T$ is well-defined we introduce one more  definition:

\medskip
\noindent
\textbf{Definition B.2}. We call two coherent sheaves $\kF$ and $\kG$  of the same rank and degree
\emph{swept-equivalent}, if  there exist two filtrations
$0 = \kF_{n+1}\subset \kF_n \subset \kF_{n-1} \subset \dots \subset \kF_1 \subset  \kF_0 = \kF$  and
$0  = \kG_{n+1} \subset \kG_n \subset \kG_{n-1} \subset \dots \subset \kG_1 \subset \kG_0 = \kG$
with quotients
$\kK_i :=  \kF_{i-1}/\kF_i \cong \kG_{n-i+1}/\kG_{n-i+2}$, $1 \le i \le n+1$.
Two such  filtrations are called \emph{admissible} filtrations associated to the  swept-equivalent pair
$(\kF, \kG)$.

\vspace{.15in}

\paragraph{\textbf{Lemma B.3.}}\textit{ For given two coherent sheaves $\kF$ and $\kG$
\begin{itemize}
\item there are finitely many  swept-equivalent pairs of coherent sheaves
$(\kF', \kG')$ such that $\kF' \rightarrowtail \kF$ and  $\kG \twoheadrightarrow \kG'$.
\item If $\kF$ and $\kG$ are themselves swept-equivalent, then
there exist only finitely many admissible filtrations associated with
$(\kF, \kG)$.
\end{itemize}}

\vspace{.1in}

\noindent
\textit{Proof}.  Let us first deal with  the second part.
Denote by  $\tau(\mathcal{H})$ the torsion part of the sheaf $\mathcal{H}$.
We argue by induction on the pair
$\bigl(\rank(\mathcal{G}), \deg(\tau(\mathcal{F}))\bigr)$, where the order is lexicographic.
The Lemma is obvious if $\rank(\mathcal{G})=0$.
Now we fix $\mathcal{F},\mathcal{G} \in Coh(\E)$ and assume we have  an
admissible filtration associated with the pair  $(\kF,\kG)$
and having the  quotients $\kK_1$, $\kK_2, \dots, \kK_{n+1}$.

Note that $\mathcal{K}_{n+1}$ is both a
subsheaf of $\mathcal{F}$ and a quotient of
$\mathcal{G}$. Hence there are only finitely many possibilities for $\mathcal{K}_{n+1}$,
and for each such $\mathcal{K}_{n+1}$, only finitely
many embeddings  $\phi: \mathcal{K}_{n+1} \rightarrowtail
 \mathcal{F}$ and quotients
$\psi: \mathcal{G}  \twoheadrightarrow  \mathcal{K}_{n+1}$.
For fixed $\phi$ and $\psi$ there is a bijection between admissible filtrations
of $(\kF, \kG)$ with quotients $\kK_1, \dots, \kK_{n+1}$
and admissible filtrations of $(\mathrm{coker}(\phi), \mathrm{ker}(\psi))$
with quotients $\kK_1, \dots, \kK_{n}$.
But $\bigl(\rank(\mathrm{ker}(\psi)),
\deg(\tau(\mathrm{coker}(\phi)))\bigr) <
\bigl(\rank(\mathcal{G}), \deg(\tau(\mathcal{F}))\bigr)$, so the induction hypothesis allows
us to conclude.

To prove the first part note that there is a bijection between the sequences of inclusions
$0 = \kH'_{n+1} \rightarrowtail \kH'_n \rightarrowtail  \kH'_{n-1} \rightarrowtail \dots \rightarrowtail  \kH'_0 = \kH$
and sequences of projections
$\kH = \kH''_0 \twoheadrightarrow \kH''_1 \twoheadrightarrow \dots \kH''_n \twoheadrightarrow
\kH''_{n+1} = 0
$
such that $\mathrm{coker}(\kH'_{i+1} \rightarrowtail \kH'_{i}) \cong \mathrm{ker}(\kH_{n-i}'' \twoheadrightarrow \kH_{n-i+1}'')$
(we can simply put $\kH''_i := \kH/\kH'_{n-i+1})$.

Therefore, existence of a swept-equivalent pair $(\kF', \kG')$, where $\kF' \rightarrowtail
 \kF$,  $\kG \twoheadrightarrow \kG'$
is equivalent to existence of a sequence of inclusions
$0  \rightarrowtail \kF'_n  \rightarrowtail  \dots \rightarrowtail
 \kF'_1 \rightarrowtail  \kF'_0 = \kF' \rightarrowtail
\kF$
and a sequence of surjections
$
\kG \twoheadrightarrow \kG' = \kG'_0  \twoheadrightarrow \kG'_1 \twoheadrightarrow \dots \twoheadrightarrow
\kG'_n \twoheadrightarrow 0
$
such that $\mathrm{coker}(\kF'_{i+1} \lto \kF'_i) \cong  \mathrm{ker}(\kG_{i}' \lto \kG_{i+1}')$.
But obviously, such sequences stand in a bijection with sequences associated with
the pair $\bigl(\kF/\kF'_n, \mathrm{ker}(\kG \lto \kG'_n)\bigr)$. This implies the first part. The lemma is proved.
\qed

\medskip

From this lemma it follows  that the operator $T: \H_\E[\alpha] \lto \widehat{\H}_\E[\alpha]$ is well-defined,
i.e.~the series (B.7)  for $T\bigl([\kF]\bigr)$ is convergent.
Indeed, for any coherent sheaf $\kG$ of class $\alpha$ there exist only finitely many admissible filtrations
associated with $(\kF,\kG)$, what means that each term $[\kG] \in \H_\E[\alpha]$ appears
the expansion of $T\bigl([\kF]\bigr)$ finitely many times.

\medskip
\noindent
For  $\alpha \in (\mathbb{Z}^2)^+$ and  let  $\Delta^n_{*}$ be the
composition of $\Delta^n$ and the canonical projection
$$
\prod_{\substack{\a_1 + \dots + \a_n = \alpha \\ \a_i \in (\mathbb{Z}^2)^+}} \H_\E[\a_1]
\widehat\otimes
\dots \widehat\otimes \H_\E[\a_n] \lto
\prod_{\substack{\a_1 + \dots +\a_n = \a \\ \a_i \in (\mathbb{Z}^2)^+, \alpha_i \ne 0}}
 \H_\E[\a_1]  \widehat\otimes
\dots \widehat\otimes \H_\E[\a_n]$$
then we denote
$$
\Delta_*(a) = \sum_{\overline{a_i^{(1)}} \neq 0, \cdots,
\overline{a_i^{(l+1)}} \neq 0} a_i^{(1)} \otimes  a_i^{(2)} \otimes \dots \otimes a_i^{(n)}.
$$
Using this notations, we may  write the operator $T: \H_\E[\alpha] \lto \widehat{\H}_\E[\alpha]$:
$$T(a)=
\left(-a+\sum_{l=1}^{\infty} (-1)^l
\sum_{\overline{a_i^{(1)}} \neq 0, \cdots,
\overline{a_i^{(l+1)}} \neq 0} a_i^{(l+1)} \cdots a_i^{(1)}\right).
$$
Note that in the case of the Hall algebra of the category of representations of a finite
quiver the map $T$ is the inverse of the antipode.

\vspace{.1in}

\paragraph{\textbf{Lemma B.4.}}\textit{
Let $\kF$ be a coherent sheaf and $\Delta\bigl([\kF]\bigr) = \sum_{i} \kF_i^{(1)} \otimes \kF_i^{(2)} \in
\H_\E \widehat\otimes \H_\E$. Then
we have   $\sum_{i}  [\kF_i^{(2)}] T\bigl([\kF_i^{(1)}]\bigr) =
\varepsilon\bigl([\kF]\bigr)1$ and
 $\sum_{i}  T\bigl([\kF_i^{(2)}]\bigr)[\kF_i^{(1)}] = \varepsilon\bigl([\kF]\bigr)1$,
where both equalities are taken in $\widehat\H_\E$.}

\vspace{.1in}

\noindent
\textit{Proof}.  We shall prove, following Theorem 1.6.3 in \cite{Kap} only the first statement, the proof of the second is  dual.
Since  the assertion trivially holds for $\kF = 0$, assume
we have a coherent sheaf
 $\kF$ of class $\alpha \ne 0$.
First of all let us check  the convergence of the series
$\sum_{i}  [\kF_i^{(2)}] T\bigl([\kF_i^{(1)}]\bigr)$ in $\widehat{\H}_\E[\a]$.
It is clear that each $[\kG] \in \H_\E[\alpha]$ gets non-zero contributions only from  finitely many summands
$[\kF_i^{(2)}] T\bigl([\kF_i^{(1)}]\bigr)$. Indeed, by
Lemma~B.3 there are only finitely many exact sequences
$0 \to \kF_i^{(2)} \to \kF \to \kF_i^{(1)} \to 0$ and
$0 \to \kG' \to \kG \to \kF_i^{(2)} \to 0$ such that $\kG'$ and $\kF_i^{(1)}$ are swept-equivalent.
Now note that
$$\sum_{i}  [\kF_i^{(2)}] T\bigl([\kF_i^{(1)}]\bigr)  =
\sum_{n=1}^\infty (-1)^n
\sum_{\kF \stackrel{\varphi_1}\twoheadrightarrow \kF_1 \stackrel{\varphi_2}\twoheadrightarrow \kF_2
\twoheadrightarrow \dots \stackrel{\varphi_n}\twoheadrightarrow \kF_n}
d^\kF_{\kB_1,\dots,\kB_{n+1}}
 [\kB_{n+1}] \otimes  \dots \otimes
[\kB_{1}],
$$
where
 $\kB_{n+1} = \kF_n$, $\kB_i =
\mathrm{ker}(\kF_{i-1} \stackrel{\varphi_i}\lto \kF_i)$
and the
sum is taken in such a way that
the epimorphisms $\varphi_2, \dots \varphi_n$ are strict  and $\varphi_1$ is arbitrary.
Now note that each term $[\kB_{n+1}] \otimes
 \dots \otimes
[\kB_{1}]$ occurs exactly twice in the sum with two different signs: one comes from the sequence
$\kF \stackrel{\varphi_1}\twoheadrightarrow \kF_1 \stackrel{\varphi_2}\twoheadrightarrow \kF_2
\twoheadrightarrow \dots \stackrel{\varphi_n}\twoheadrightarrow \kF_n$  where all epimorphisms
 $\varphi_i$
are strict and the second comes from
$\kF \stackrel{id}\lto  \kF \stackrel{\varphi_1}\twoheadrightarrow \kF_1
\stackrel{\varphi_2}\twoheadrightarrow \kF_2
\twoheadrightarrow \dots \stackrel{\varphi_n}\twoheadrightarrow \kF_n$.
This shows the lemma.
\qed

\medskip

\paragraph{\textbf{Lemma B.5.}}\textit{
Let $\a,\b \in (\mathbb{Z}^2)^+$, $a \in \H_\E[\a]$ and  $b\in \H_\E[\b]$. Then we have
$T(ab) = T(b)T(a)$ in $\widehat{\H}_\E[\a + \b]$. }

\vspace{.1in}

\noindent
\textit{Proof}. Let $\Delta^2(a) = \sum_i a_i^{(1)} \otimes a_i^{(2)} \otimes a_i^{(3)}$ and
$\Delta^2(b) = \sum_j b_j^{(1)} \otimes b_j^{(2)} \otimes b_j^{(3)}$. Consider an expression
$$
c = \sum_{i,j} T\bigl(b_j^{(3)}\bigr)  T(a_i^{(3)}) a_i^{(2)} b_j^{(2)} 
T\bigl(a_i^{(1)}b_j^{(1)}\bigr)
$$
in $\widehat{\H}_\E[\a+\b]$. To see that this sum converges, assume that both $a$ and $b$ are classes of
coherent sheaves.  A coherent sheaf $\kF$ of class $\a + \b$ enters in the sum $c$ if and only if
there is a filtration $ 0 \rightarrowtail \kF_4  \stackrel{\varphi_3}\rightarrowtail  \kF_3
\stackrel{\varphi_2}\rightarrowtail
 \kF_2  \stackrel{\varphi_1}\rightarrowtail \kF_1  \stackrel{\varphi_0}\rightarrowtail  \kF_0 = \kF$
such that $\mathrm{coker}(\varphi_0)$ is swept-equivalent to $b_j^{(3)}$,  $\mathrm{coker}(\varphi_1)$ is swept-equivalent to
$a_i^{(3)}$, $\mathrm{coker}(\varphi_2)$ is isomorphic to $a_i^{(2)}$,
$\mathrm{coker}(\varphi_3)$ is isomorphic to $b_j^{(2)}$ and  finally
$\kF_4$ is swept-equivalent to $a_i^{(1)} b_j^{(1)}$.

Since $b_j^{(3)}$ is a subsheaf of $b$ and is swept-equivalent to $\mathrm{coker}(\varphi_0)$, by Lemma~B.3
 there are finitely many
choices for $b_j^{(3)}$ and $\varphi_0$, and therefore  finitely many contributions of  $T(b_j^{(3)})$ to $\kF$.
Assuming that $b_j^{(3)}$ and $\kF_1 \stackrel{\varphi_0}\rightarrowtail \kF$
are fixed, by the same argument we see that there are finitely many subobjects $a_i^{(3)}$ of $a$ and finitely many inclusions
$\kF_2 \stackrel{\varphi_1}\rightarrowtail \kF_1$ such that $a_i^{(3)}$ and $\mathrm{coker}(\varphi_1)$ are swept-equivalent.
Next, there are finitely many inclusions $\kF_3 \stackrel{\varphi_2}\rightarrowtail \kF_2$ such that there is a subobject
$a_i^{(2)}$ of $a/a_i^{(3)}$ isomorphic to  $\mathrm{coker}(\varphi_2)$. In the same way, we have only finitely many
inclusions  $\kF_4 \stackrel{\varphi_3}\rightarrowtail \kF_3$ such that $\mathrm{coker}(\varphi_3)$ is isomorphic to a subobject
of $b/b_j^{(3)}$.  But choices of $a_i^{(2)}$ and $a_i^{(3)}$ also determine $a_i^{(1)}$, the same holds for
$b_j^{(2)}$ and $b_j^{(3)}$, hence there are finitely many subobjects of $\kF_4$ swept-equivalent to
some summand of $a_i^{(1)}b_j^{(1)}$.  Gathering all together we conclude, that the element $c$ is correctly defined in
$\widehat{\H}_\E[\alpha +\beta]$.

Using Lemma~B.4 we can transform the series $c$ in two different ways.
From the one side  we have
\begin{equation*}
\begin{split}
c =
\sum_{j} T(b_j^{(3)})  \varepsilon(a_i^{(2)})  b_j^{(2)} T(a_i^{(1)}b_j^{(1)}) &=
\sum_{j} T(b_j^{(3)})  b_j^{(2)} T(a b_j^{(1)}) \\
&=\sum_{j} \varepsilon(b_j^{(2)})  T(a b_j^{(1)}) =   T(ab)
\end{split}
\end{equation*}
and from another side,
\begin{equation*}
\begin{split}
c  &=
\sum_{i,j,k,l} T(b_j^{(2)})  T(a_i^{(2)}) (a_i^{(1)})_l^{(2)} (b_j^{(1)})_k^{(2)}
T\bigl((a_i^{(1)})_l^{(1)} (b_{j}^{(1)})_k^{(1)})\\
&=\sum_{i,j,k} T(b_j^{(2)})  T(a_i^{(2)}) (a_i^{(1)} b_j^{(1)})_k^{(2)} T\bigl((a_i^{(1)} b_j^{(1)})_k^{(1)}\bigr)\\
&= \sum_{i,j} T(b_j^{(2)})  T(a_i^{(2)}) \varepsilon(a_i^{(1)} b_j^{(1)}) = T(b) T(a).
\end{split}
\end{equation*}
\qed

\paragraph{\textbf{Lemma B.6.}}\textit{
Let $a \in \H_\E^-$ and $b\in \H_\E^+$, then the following equation holds in the Drinfeld double
$\mathbf{D}\H_\E$:}
\begin{equation}
ab=\sum_{i,j} \bigl(a_i^{(1)},b_j^{(3)}\bigr) b_j^{(2)}a_i^{(2)}\bigl(T(a_i^{(3)}),b_j^{(1)}\bigr).\tag{B.8}
\end{equation}
\textit{Proof}. First of all note, that the right-hand side of the equation (B.8)
is finite by Lemma~B.3.
The relation $R(a,b)$ in the Drinfeld double implies
$$
ab = - \sum_{\underset{\overline{a_i^{(2)}} \ne 0}{i,j}} a_i^{(1)} b_j^{(2)}
\bigl(a_i^{(2)}, b_j^{(1)}\bigr) +
\sum_{i,j} b_j^{(1)}a_i^{(2)}\bigl(a_i^{(1)}, b_j^{(2)}\bigr).
$$
Now use this equality to rewrite each term $a_i^{(1)} b_j^{(2)}$:
\begin{equation*}
\begin{split}
a_i^{(1)} b_j^{(2)} = -
\sum_{\underset{\overline{(a_{i}^{(2)})_k^{(2)}} \ne 0}{k,l}}
(a_i^{(1)})_k^{(1)} (b_j^{(2)})_l^{(2)} & \bigl((a_i^{(1)})_k^{(2)}, (b_j^{(2)})_l^{(1)}\bigr) + \\&
+\sum_{k,l} (b_j^{(2)})_l^{(1)} (a_i^{(1)})_k^{(2)} \bigl((a_i^{(1)})_{k}^{(1)}, (b_j^{(2)})_l^{(2)}\bigr).
\end{split}
\end{equation*}
Combining these two equations we obtain
\begin{equation*}
\begin{split}
ab =  \sum_{\underset{\overline{a_i^{(3)}} \ne 0, \overline{a_i^{(2)}}\ne 0}{i,j}}
a_i^{(1)} b_j^{(3)}(a_i^{(2)}, b_j^{(2)})(a_i^{(3)}, b_j^{(1)}) - &
\sum_{\underset{\overline{a_i^{(3)}} \ne 0}{i,j}} b_j^{(2)}a_i^{(2)} (a_i^{(1)}, b_j^{(3)})(a_i^{(3)}, b_j^{(1)}) +
\\& + \sum_{i,j} b_j^{(1)}a_i^{(2)} (a_i^{(1)}, b_j^{(2)}).
\end{split}
\end{equation*}
But note that
$$
\sum_{\underset{\overline{a_i^{(3)}} \ne 0, \overline{a_i^{(2)}} \ne 0}{i,j}}
a_i^{(1)} b_j^{(3)}(a_i^{(2)}, b_j^{(2)})(a_i^{(3)}, b_j^{(1)}) =
\sum_{\underset{\overline{a_i^{(3)}} \ne 0, \overline{a_i^{(2)}} \ne 0}{i,j}} a_i^{(1)} b_j^{(2)} (a_i^{(3)}
a_i^{(2)}, b_j^{(1)}).
$$
Moreover, we have
$$
\sum_{i,j} b_j^{(1)}a_i^{(2)} (a_i^{(1)}, b_j^{(2)}) =
\sum_{\underset{\overline{a_i^{(3)}} = 0}{i,j}} (a_i^{(1)}, b_j^{(3)}) b_j^{(2)} a_i^{(2)}
(a_i^{(3)}, b_j^{(1)}).
$$
Summing everything up, we get
\begin{equation*}
\begin{split}
ab = \sum_{\underset{\overline{a_i^{(3)}} = 0}{i,j}} (a_i^{(1)}, b_j^{(3)}) b_j^{(2)} a_i^{(2)}
(a_i^{(3)}, b_j^{(1)}) -  &
\sum_{\underset{\overline{a_i^{(3)}} \ne 0}{i,j}} b_j^{(2)}a_i^{(2)} (a_i^{(1)}, b_j^{(3)})(a_i^{(3)}, b_j^{(1)}) + \\ &
+ \sum_{\underset{\overline{a_i^{(3)}} \ne 0, \overline{a_i^{(2)}} \ne 0}{i,j}} a_i^{(1)} b_j^{(2)} (a_i^{(3)}
a_i^{(2)}, b_j^{(1)}).
\end{split}
\end{equation*}
Iterating this procedure, we get, for each $k > 0$,
\begin{equation}
\begin{split}
ab=\sum_{i,j} (a_i^{(1)},b_j^{(3)}) &b_j^{(2)}a_i^{(2)}\bigl(T^k(a_i^{(3)}),b_j^{(1)}\bigr) +\\
&+(-1)^k
\sum_{\underset{\overline{a_i^{(2)}} \neq 0, \ldots, \overline{a_i^{(k+1)}} \neq 0}{i,j}}
a_i^{(1)}b_j^{(2)}(a_i^{(k+1)} \cdots a_{i}^{(2)}, b_j^{(1)}),
\end{split}\tag{B.9}
\end{equation}
where
$$
T_k(a)=\left(-a+\sum_{l=1}^{k-2} (-1)^l
\sum_{\overline{a_i^{(1)}} \neq 0, \cdots,
\overline{a_i^{(l+1)}} \neq 0} a_i^{(l+1)}
\cdots a_i^{(1)}\right).
$$
It follows from the first part of
Lemma~B.3 that the second term in (B.9) vanishes for
$k \gg 0$ and the operators $T^k$ converge pointwise as $k \to \infty$ to operators
$T:  \H_{\E}[\a] \to \widehat{\H}_{\E}[\a]$ yielding the equation
(B.8) as wanted.
\qed

\medskip

This lemma shows that the  map $\H_\E^+ \otimes \H_\E^- \stackrel{m}\lto
\mathbf{D}\H_\E$ is surjective.
Next, we define an associative algebra structure on $\H^+_\E \otimes \H^-_\E$ by setting
$$
(a\otimes a')\cdot(b\otimes b') = (m\otimes m)(a \otimes L(a',b) \otimes b'),
$$
where
$$
L(x,y) = \sum_{i,j}\bigl(x_1^{(i)}, y_j^{(3)}\bigr) y_j^{(2)} x_{i}^{(2)}\bigl(T(x_i^{(3)}), y_j^{(1)}\bigr).
$$
A proof  of associativity  of this product is based on Proposition \ref{P:topbial}, Remark \ref{R:hopfp}
and Lemma~B.5 and can be shown along the same lines as in
in \cite[3.2.4]{Joseph} using similar calculations as in the proof of
Lemma~B.6.

Now we  can  construct the inverse
map $n: \mathbf{D}\H_\E \lto \H_\E^+ \otimes \H_\E^-$
by putting $n(a) = 1 \otimes a$, $n(b) = b \otimes 1$ for $ a \in \H_\E^- \subset
\mathbf{D}\H_\E$ and  $b \in \H_\E^+ \subset
\mathbf{D}\H_\E$.
To see that  we get  a well-defined map we need to check that all relations $R(a,b)$ are
preserved in the Drinfeld double. But indeed,

\begin{equation*}
\begin{split}
\sum_{i,j} a_i^{(1)} b_j^{(2)}(a_i^{(2)}, b_j^{(1)}) =
\sum_{i,j}(a_i^{(1)}, b_j^{(4)}) a_i^{(2)}  b_j^{(3)}
\bigl(T(a_i^{(3)}), b_j^{(2)}\bigr)(a_i^{(4)}, b_j^{(1)}) = & \\
\sum_{i,j}(a_i^{(1)}, b_j^{(3)}) a_i^{(2)}  b_j^{(2)} \bigl(T(a_i^{(3)}) a_i^{(4)}, b_j^{(1)}\bigr) =
\sum_{i,j}(a_i^{(1)}, b_j^{(3)}) a_i^{(2)}  b_j^{(2)} \bigl(\varepsilon(a_i^{(3)})1, b_j^{(1)}\bigr) = &\\
= \sum_{i,j} (a_i^{(1)}, b_j^{(2)}) a_i^{(2)}b_j^{(1)}.
\end{split}
\end{equation*}
This concludes the proof of injectivity and surjectivity of the linear map \linebreak
$m: \H_\E^+ \otimes \H_\E^- \lto \mathbf{D}\H_\E$. Proposition~\ref{P:DDoub} is proven. \qed

\vspace{.2in}

\centerline{\textbf{Acknowledgments}}

\vspace{.1in}

We would like to thank Yu.~Drozd, P.~Etingof, S.~Ovsienko, E.~Vasserot and especially I. Gordon
for useful discussions on the results of this article.
Parts of this work
were conducted at ENS Paris, Max-Plank-Institut in Bonn,
Mathematical  Research Institute  Oberwolfach (RIP Programme), Johannes Gutenberg Universit\"at Mainz
 and Yale University. We would like to thank
these institutions for their hospitality. The  first-named author was supported by the DFG grants
Bu 1866/1-1 and  Bu 1866/1-2.

\small{}

\end{document}